\documentclass[12pt,twoside]{article}  

\usepackage{amssymb}

\usepackage{amsfonts}
\usepackage{amssymb}
\usepackage{mathrsfs}
\usepackage{amsmath}
\usepackage{graphicx}
\usepackage{multirow}
\usepackage{booktabs}
\usepackage{bm,amsthm,stmaryrd}
\usepackage[colorlinks,bookmarksopen,bookmarksnumbered,citecolor=blue, linkcolor=blue, urlcolor=blue]{hyperref}
\usepackage[ruled,vlined]{algorithm2e}
\numberwithin{equation}{section}
\allowdisplaybreaks[4]
\usepackage{cases} 

\usepackage[left=2cm, right=2cm, top=1.0in, bottom=1.0in]{geometry}
\usepackage[marginal]{footmisc}

\newtheorem{theorem}{Theorem}[section]
\newtheorem{lemma}{Lemma}[section]
\newtheorem{example}{Example}[section]

\begin{document}

\date{}
\author{Min Cai, Changpin Li\thanks{Corresponding author. E-mail: lcp@shu.edu.cn}\,, Yu Wang
\\
\small \textit{Department of Mathematics, Shanghai University, Shanghai 200444, China}
}   
\title{Two kinds of numerical algorithms for ultra-slow diffusion equations\thanks{This work was funded by the National Natural Science Foundation of China under Grant Nos.~12201391 and 12271339.}}
\maketitle
 \hrulefill
 
\footnote{Each of the authors equally contributes to this work.}

\begin{abstract}
In this article, two kinds of numerical algorithms are derived for the ultra-slow (or superslow) diffusion equation in one and two space dimensions, where the ultra-slow diffusion is characterized by the Caputo-Hadamard fractional derivative of order $\alpha \in (0,1)$. To describe the spatial interaction, the Riesz fractional derivative and the fractional Laplacian are used in one and two space dimensions, respectively. The Caputo-Hadamard derivative is discretized by two typical approximate formulae, i.e., L2-1$_{\sigma}$ and L1-2 methods.
The spatial fractional derivatives are discretized by the 2-nd order finite difference methods. When L2-1$_{\sigma}$ discretization is used, the derived numerical scheme is unconditionally stable with error estimate $\mathcal{O}(\tau^{2}+h^{2})$  for all $\alpha \in (0, 1)$, in which $\tau$ and $h$ are temporal and spatial stepsizes, respectively. When L1-2 discretization is used, the derived numerical scheme is stable with error estimate $\mathcal{O}(\tau^{3-\alpha}+h^{2})$ for $\alpha \in (0, 0.3738)$. The illustrative examples displayed are in line with the theoretical analysis.
\end{abstract}


\vskip 5mm


\hrulefill

\section{Introduction}
\label{sec:Introduction}
Scientific explorations on ultra-slow (or superslow) process originates from the study on Lomnitz's logarithmic creep law of igneous \cite{Lom1956}, elastic crack or fatigue fracture of viscoelastic solid materials \cite{Lom1957}, and stress wave attenuation in the solid earth \cite{Lom1962}.  It is also reported that ultra-slow process appears in geochemistry and geophysics \cite{Bre2021,War2016}, hydrothermal fluid \cite{Tao2020}, and neuroscience \cite{Sil2013}. Ultra-slow diffusion is one of the most common phenomena among ultra-slow process \cite{Ou2022,Yang2022,Zaky2022}. From the perspective of continuous time random walk, the ultra-slow diffusion is caused by the super-heavy tailed distribution of waiting time in  particle motion \cite{Den2010}, which can be characterized by a logarithmic growth of the mean squared displacement as a function of the temporal variable \cite{Lom2016}. 

In the existing time-space fractional models for ultra-slow diffusion which take the temporal logarithmic decay and the spatial non-locality into account, the temporal derivative can be governed by the Caputo-Hadamard (fractional) derivative \cite{Diethelm2010,LiL2020}
\begin{equation*}
    _{CH}{\rm D}_{\tilde{a},t}^\alpha \varphi(t)
    =\frac{1}{\Gamma(n-\alpha)}
    \int^{t}_{\tilde{a}}\left(\log\frac{t}{s}\right)^{n-\alpha-1}\delta^{n}\varphi(s)\frac{{\rm d}s}{s},\quad 0<\tilde{a}<t,
\end{equation*}
with
\begin{equation*}
    \delta^{n}\varphi(s) 
    = \left(s\frac{{\rm d}}{{\rm d}  s}\right)^{n}\varphi(s),
    \  
    n-1 < \alpha< n \in \mathbb{Z}^+.  
\end{equation*}
The spatial fractional derivative is usually given by the Riesz (fractional) derivative in one dimension \cite{Bai22021,DingL2017,Lin2017}
\begin{equation*}
    _{RZ}{\rm D}_{x}^\beta \psi(x) 
    = -\frac{1}{2\cos(\pi \beta /2)}
    \left[{}_{RL}{\rm D}_{a,x}^\beta \psi(x)+ {}_{RL}{\rm D}_{x,b}^\beta \psi(x)\right],\, m-1< \beta < m \in \mathbb{Z}^+,
\end{equation*}
or fractional Laplacian in higher dimensions \cite{Jiao2021,LiandZeng2015,Wan2022} which is defined by a singular or hyper-singular integral, 
\begin{equation}\label{eq:lap} 
    \left(-\Delta\right)^{\frac{\beta}{2}}\psi(\mathbf{x}) 
    = c_{\beta}\, \mathrm{P.V.} \int_{\mathbb{R}^d}
    \frac{\psi(\mathbf{x})-\psi(\mathbf{z})}{\lvert\mathbf{x}-\mathbf{z}\rvert^{d+\alpha}}{\rm d}\mathbf{z},
    \, c_{\beta}=\frac{2^\beta \Gamma(d/2+\beta/2)}{\pi^{d/2}\lvert\Gamma(-\beta/2)\rvert}, 
\end{equation}
with $\mathrm{P.V.}$ standing for the Cauchy principal value. Here 
\begin{equation*}
    _{RL}{\rm D}_{a,x}^\beta \psi(x) 
    = \frac{1}{\Gamma(m - \beta)} \frac{{\rm d}^{m}}{{\rm d} x^{m}}\int_{a}^{x}(x-s)^{m-\beta-1}\psi(s){\rm d}s
\end{equation*}
and
\begin{equation*}
    _{RL}{\rm D}_{x,b}^\beta \psi(x) 
    = \frac{(-1)^m}{\Gamma(m - \beta)} \frac{{\rm d}^{m}}{{\rm d} x^{m}}\int_{x}^{b}(s-x)^{m-\beta-1}\psi(s){\rm d}s
\end{equation*}
are the left- and right-side Riemann-Liouville (fractional) derivatives. The sufficient condition for the existence of Caputo-Hadamard derivative ${}_{CH}{\rm D}_{\tilde{a},t}^\alpha \varphi(t)$ is $\varphi(t)\in AC^{n}_{\delta}[\tilde{a},T]=\left\{\varphi: \delta^{n-1}\varphi\in AC[\tilde{a},T] \right\}$ with $AC(\Omega)$ denoting the space of absolute continuous functions defined on $\Omega$. Sufficient condition for the existence of ${}_{RZ}{\rm D}_{x}^\beta \psi(x)$ is that
$$
\psi(x)\in AC^{\lceil\beta\rceil}[a,b]=\left\{\psi: \psi^{(k)}\in AC[a,b], k=0,1,\cdots,\lfloor\beta\rfloor\right\}.
$$
One sufficient condition for the existence of the fractional Laplacian $\left(-\Delta\right)^{\frac{\beta}{2}}\psi(\mathbf{x})$ is that $\psi(\mathbf{x})$ belongs to the following Schwartz space  \cite{LiCai2019}
\begin{equation*}
\mathcal{S}=\left\{\psi\in C^{\infty}(\mathbb{R}^{d}): \sup\limits_{\mathbf{x}\in\mathbb{R}^{d}}(1+\lvert\mathbf{x}\rvert)^{N}\sum\limits_{k=0}^{N}\lvert D^{k}\psi(\mathbf{x})\rvert<+\infty, 
N = 0,1,2,\cdots\right\}.
\end{equation*}

From Lemma 3.1 in \cite{LiL2021}, one knows that the Caputo-Hadamard derivative can be changed into the Caputo derivative. Afterwards, some other works, e.g., \cite{Zheng2021} and \cite{Alikhanov2022}, also display this point. But it is inconvenient for analysis and computation if we do like this. So it is completely unnecessary to change the Caputo-Hadamard derivative into the Caputo one whenever we meet the Caputo-Hadamard type differential equations.

In this paper, we first study the following one-dimensional ultra-slow diffusion equation,
\begin{equation}\label{eq:FAE}
    \left\{
    \begin{array}{ll}
    {}_{CH}{\rm D}^{\alpha}_{\tilde{a},t}u(x,t)
    -{}_{RZ}{\rm D}^{\beta}_{x}u(x,t)=f(x,t), 
    & (x,t)\in\Omega\times(\tilde{a},T], 
    \\
    u(a,t) = u(b,t) = 0, 
    & t\in(\tilde{a},T], 
    \\
    u(x,\tilde{a}) = u_{0}(x), 
    & x\in\Omega, 
    \end{array}
    \right.
\end{equation} 
with $\alpha\in(0,1)$, $\beta\in(1,2)$, $\tilde{a}>0$, $\Omega=(a,b)\subset\mathbb{R}$ being a bounded interval, $u_{0}(x)$ being a known function, and $f(x,t)$ being a source term. 

Next, we study the two-dimensional ultra-slow diffusion equation of the form 
\begin{equation}\label{eq:FAE2}
    \left\{
    \begin{array}{ll}
    {}_{CH}{\rm D}^{\alpha}_{\tilde{a},t}u(x,y,t)
    +\left(-\Delta\right)^{\frac{\beta}{2}}u(x,y,t)=f(x,y,t), 
    & (x,y)\in\widetilde{\Omega},\, t \in (\tilde{a},T], 
    \\
    u(x,y,t) = 0, 
    & (x,y)\in \mathbb{R}^2\backslash \widetilde{\Omega},\, t \in (\tilde{a},T], 
    \\
    u(x,y,\tilde{a}) = u_{0}(x,y), 
    & (x,y)\in\widetilde{\Omega}, 
    \end{array}
    \right.
\end{equation}
with $\alpha\in(0,1)$, $\beta\in(1,2)$, $\tilde{a}>0$,  $\widetilde{\Omega}=(-L,L)^2\subset\mathbb{R}^2$ being a bounded domain, $u_{0}(x,y)$ being a given  function, and $f(x,y,t)$ being a source term. 

Models \eqref{eq:FAE} and \eqref{eq:FAE2} have been considered in \cite{WangCai2023}, where the Caputo-Hadamard derivative is approximated by the L1 method. In this paper, we continue to study Eqs. \eqref{eq:FAE} and \eqref{eq:FAE2} in order to improve the accuracy in time direction. To this end, we adopt two approximations of $(3-\alpha)$-th order for Caputo-Hadamard derivative of order $\alpha\in(0,1)$. Say, L2-1$_{\sigma}$ and L1-2 formulae. In the spatial discretization of this paper, we choose the weighted and shifted Gr\"{u}nwald-Letnikov formula for the Riesz derivative in one space dimension, and the fractional centered difference formula for the fractional Laplacian in two space dimensions. Although both L2-1$_{\sigma}$ and L1-2 formulae are of $(3-\alpha)$-th order when evaluating Caputo-Hadamard derivative, the resulting fully discrete schemes are of different accuracy in temporal direction. For the numerical schemes based on L2-1$_{\sigma}$ formula which is used to evaluate Caputo-Hadamard derivative at the non-integer temporal grid point $t_{k+\sigma}=\tilde{a}+(k+\sigma)\tau$ with $\sigma=1-\frac{\alpha}{2}$ via numerical combinations of the values at integer grid point $t_j \,(0\leq j\leq k+1)$, necessary pre-process is needed to achieve the unification at the temporal level. The interpolation technique is adopted to achieve this goal and the resulting temporal accuracy is of 2-nd order. And the corresponding numerical schemes are proved to be unconditionally stable. On the other side, the numerical schemes based on L1-2 formula do not have this issue and their temporal accuracy remains $(3-\alpha)$-th order. However, for the time being, they are proved to be numerically stable only when $\alpha$ belongs to a proper subset of $(0,1)$.

The remaining part of this paper is organized as follows. Sec. \ref{sec:Preliminaries} introduces numerical approximations of Caputo-Hadamard derivative, Riesz derivative, and fractional Laplacian those are adopted in the present paper, along with their properties. Fully discrete schemes for Eqs. \eqref{eq:FAE} and \eqref{eq:FAE2} with L2-1$_{\sigma}$ formula in temporal discretization are constructed and numerically analyzed in Sec. \ref{sec:Scheme}. Construction and numerical analysis of the fully discrete schemes for Eqs. \eqref{eq:FAE} and \eqref{eq:FAE2} based on L1-2 approximation in temporal discretization are shown in Sec. \ref{sec:Scheme2}. Numerical examples in Sec. \ref{sec:Experiments} verify the theoretical results given in Secs. \ref{sec:Scheme} and \ref{sec:Scheme2}. For the sake of simplicity, the constants $C$ and $\widetilde{C}$ appeared in the present paper may be different for distinct cases.

\section{Preliminaries}
\label{sec:Preliminaries}
The current paper focuses on numerical algorithms for the ultra-slow diffusion equations in both one-dimensional and two-dimensional cases. In this preliminary section, we introduce the numerical approximations which are cornerstones of the numerical schemes to be proposed in the coming sections, along with their properties.  

\subsection{L2-1$_{\sigma}$ formula for Caputo-Hadamard derivative}
Let $t_k = \tilde{a} + k\tau$ with $k = 0,1,\cdots ,N \in\mathbb{Z}^{+}$ and define $t_{k+\sigma} = t_{k} + \sigma \tau $ with $\sigma = 1 - \frac{\alpha}{2}$, where $\tau = (T-\tilde{a})/N$ is the temporal stepsize. In this setting, if $\varphi(t)\in C^3[\tilde{a},T]$, the Caputo-Hadamard derivative ${}_{CH}{\rm D}_{\tilde{a},t}^\alpha \varphi(t)$ at $t = t_{k+\sigma}\ (0\leq k \leq N-1)$ can be evaluated by the following L2-1$_{\sigma}$ approximation {\cite{Fan2022}}, 
\begin{equation}\label{eq:L1_CH}
    \begin{split}
      \left.{}_{CH}{\rm D}_{\tilde{a},t}^\alpha \varphi(t)\right|_{t=t_{k+\sigma}}
      \approx
      \sum_{i=1}^{k+1} c^{(\alpha,\sigma)}_{i,k}\left[\varphi(t_i)-\varphi(t_{i-1})\right]  
    \end{split}
\end{equation}
with the numerical error being $\mathcal{O}(\tau^{3-\alpha})$.  When $k=0, 1$, the coefficients are given by 
\begin{equation}\label{eq:c_{j,k}_1}
\begin{aligned}
    &c_{1,0}^{(\alpha,\sigma)}
    =\frac{1}{\Gamma(2-\alpha)}\frac{1}{\log\frac{t_1}{t_0}}\left(\log\frac{t_{\sigma}}{t_0}\right)^{1-\alpha}, 
    \qquad\qquad\quad 
    k = 0, 
    \\
    &
    \left\{
    \begin{array}{ll}
    c_{1,1}^{(\alpha,\sigma)}
    ={\displaystyle \frac{1}{\Gamma(2-\alpha)}\frac{1}{\log\frac{t_1}{t_0}}}
    \left(a_{1,1}^{(\alpha,\sigma)}-b_{1,1}^{(\alpha,\sigma)}\right),
    \\
    c_{2,1}^{(\alpha,\sigma)}
    ={\displaystyle \frac{1}{\Gamma(2-\alpha)}\frac{1}{\log\frac{t_2}{t_1}}}
    \left(b_{1,1}^{(\alpha,\sigma)}+\left(\log\frac{t_{1+\sigma}}{t_1}\right)^{1-\alpha}\right),
    \end{array}
    \right. k =1. 
\end{aligned}
\end{equation}
For the case with $k\geq 2$, the coefficients are defined by  
\begin{equation}\label{eq:c_{j,k}}
\begin{split}
    c^{(\alpha,\sigma)}_{i,k}=&
    \left\{
    \begin{array}{ll}
    {\displaystyle \frac{1}{\Gamma(2-\alpha)}\frac{1}{\log\frac{t_1}{t_0}}}
    \left(a_{1,k}^{(\alpha,\sigma)}-b_{1,k}^{(\alpha,\sigma)} \right), &i = 1,
    \\
    {\displaystyle \frac{1}{\Gamma(2-\alpha)}\frac{1}{\log\frac{t_{i}}{t_{i-1}}}}
    \left(a_{i,k}^{(\alpha,\sigma)}+b_{i-1,k}^{(\alpha,\sigma)}-b_{i,k}^{(\alpha,\sigma)} \right), &2\leq i \leq k,
    \\
    {\displaystyle \frac{1}{\Gamma(2-\alpha)}\frac{1}{\log\frac{t_{k+1}}{t_{k}}}}
    \left(b_{k,k}^{(\alpha,\sigma)}+\left(\log\frac{t_{k+\sigma}}{t_k}\right)^{1-\alpha} \right), &i = k+1. 
    \end{array}
    \right.
\end{split}
\end{equation}
Here 
\begin{equation*} 
\begin{split}
    a_{i,k}^{(\alpha,\sigma)}=&\left(\log\frac{t_{k+\sigma}}{t_{i-1}}\right)^{1-\alpha} - \left(\log\frac{t_{k+\sigma}}{t_{i}}\right)^{1-\alpha},
    \\
    b_{i,k}^{(\alpha,\sigma)}
    =\,&\Bigg\{
     \frac{2}{2-\alpha}\left[\left(\log\frac{t_{k+\sigma}}{t_{i-1}}\right)^{2-\alpha}-\left(\log\frac{t_{k+\sigma}}{t_i}\right)^{2-\alpha}\right] 
    \\
    &-\log\frac{t_i}{t_{i-1}}\left[\left(\log\frac{t_{k+\sigma}}{t_{i}}\right)^{1-\alpha} + \left(\log\frac{t_{k+\sigma}}{t_{i-1}}\right)^{1-\alpha} \right]
    \Bigg\}
    \frac{1}{\log\frac{t_{i+1}}{t_{i-1}}}. 
\end{split}
\end{equation*} 

For properties of the coefficients $c_{i,k}^{(\alpha,\sigma)}$, we have the following results. 
\begin{lemma}{\rm\cite{Fan2022}}\label{lemma:lc_kk}
For $\alpha \in (0,1)$, $\sigma = 1-\frac{\alpha}{2}$, and sufficiently small temporal stepsize $\tau$, the coefficients $c_{i,k}^{(\alpha,\sigma)}$ given by Eqs. \eqref{eq:c_{j,k}_1} and \eqref{eq:c_{j,k}} satisfy:
\begin{equation*}
    c^{(\alpha,\sigma)}_{k+1,k}
    >c^{(\alpha,\sigma)}_{k,k}
    >\cdots
    >c^{(\alpha,\sigma)}_{i,k}
    >c^{(\alpha,\sigma)}_{i-1,k}
    >\cdots
    >c^{(\alpha,\sigma)}_{1,k}>0. 
\end{equation*}
\end{lemma}

\begin{lemma}\label{lemma:lc_kksigma}
For $\alpha \in (0,1)$, $\sigma = 1-\frac{\alpha}{2}$, and sufficiently small  temporal stepsize $\tau$, the following inequality holds for the coefficients $c_{i,k}^{(\alpha,\sigma)}$ defined in Eqs. \eqref{eq:c_{j,k}_1} and \eqref{eq:c_{j,k}},  
\begin{equation*}
    (2\sigma-1)c^{(\alpha,\sigma)}_{k+1,k}
    >\sigma c^{(\alpha,\sigma)}_{k,k}, 
    \ 
    k\geq 1.
\end{equation*}
\begin{proof}  
We first show the case with $k=1$. Note that when $\tau \rightarrow 0$, 
\begin{equation*}
\begin{aligned}
     a_{1,1}^{(\alpha,\sigma)}
    =&\left(\log\frac{t_{1+\sigma}}{t_{0}}\right)^{1-\alpha} - \left(\log\frac{t_{1+\sigma}}{t_{1}}\right)^{1-\alpha}
    \\ 
    \approx
    &\left(\frac{(1+\sigma)\tau}{t_0}\right)^{1-\alpha} - \left(\frac{\sigma\tau}{t_1}\right)^{1-\alpha}
    \\
    =&\frac{\tau^{1-\alpha}}{t_0^{1-\alpha}}\left[\left(1+\sigma\right)^{1-\alpha} - \sigma^{1-\alpha}\left(\frac{t_0}{t_1}\right)^{1-\alpha}\right], 
\end{aligned}
\end{equation*}
and 
\begin{equation*} 
\begin{split} 
    b_{1,1}^{(\alpha,\sigma)}
    &=\frac{1}{{\displaystyle \log\frac{t_{2}}{t_{0}}}}\frac{2}{2-\alpha}
    \left[\left(\log\frac{t_{1+\sigma}}{t_{0}}\right)^{2-\alpha}-\left(\log\frac{t_{1+\sigma}}{t_1}\right)^{2-\alpha}\right] 
    \\
    &-\frac{1}{{\displaystyle\log\frac{t_{2}}{t_{0}}}}\log\frac{t_1}{t_{0}}\left[\left(\log\frac{t_{1+\sigma}}{t_{1}}\right)^{1-\alpha} + \left(\log\frac{t_{1+\sigma}}{t_{0}}\right)^{1-\alpha}\right]
    \\ 
    &\approx
    \frac{t_0}{2\tau}\frac{2}{2-\alpha}\left[\left(\frac{(1+\sigma)\tau}{t_0}\right)^{2-\alpha}-\left(\frac{\sigma\tau}{t_1}\right)^{2-\alpha}\right] 
    -\frac{t_0}{2\tau}\frac{\tau}{t_0}\left[\left(\frac{\sigma\tau}{t_1}\right)^{1-\alpha} + \left(\frac{(1+\sigma)\tau}{t_0}\right)^{1-\alpha}\right]
    \\
    &=\frac{\tau^{1-\alpha}}{t_0^{1-\alpha}}\left\{\frac{1}{2-\alpha}\left[\left(1+\sigma\right)^{2-\alpha}-\sigma^{2-\alpha}\left(\frac{t_0}{t_1}\right)^{2-\alpha}\right] 
    -\frac{1}{2}\left[\sigma^{1-\alpha}\left(\frac{t_0}{t_1}\right)^{1-\alpha} + \left(1+\sigma\right)^{1-\alpha}\right]\right\}. 
\end{split}
\end{equation*}
It follows from the definition of the coefficients $c_{i,k}^{(\alpha,\sigma)}$ that  
\begin{equation*}
\begin{aligned}
    \frac{c_{2,1}^{(\alpha,\sigma)}}{ c_{1,1}^{(\alpha,\sigma)}} 
    =&\frac{\log\frac{t_1}{t_0}}{\log\frac{t_2}{t_1}} 
    \,
    \frac{b_{1,1}^{(\alpha,\sigma)}+\left(\log\frac{t_{1+\sigma}}{t_1}\right)^{1-\alpha}}{a_{1,1}^{(\alpha,\sigma)}-b_{1,1}^{(\alpha,\sigma)}}
    \\ 
    \approx
    &\frac{\frac{1}{2-\alpha}\left[\left(1+\sigma\right)^{2-\alpha}-\sigma^{2-\alpha}\right] 
    -\frac{1}{2}\left[\left(1+\sigma\right)^{1-\alpha}-\sigma^{1-\alpha} \right]}
    {\left(1+\sigma\right)^{1-\alpha} - \frac{1}{2-\alpha}\left[\left(1+\sigma\right)^{2-\alpha}-\sigma^{2-\alpha}\right] 
    +\frac{1}{2}\left[\left(1+\sigma\right)^{1-\alpha}-\sigma^{1-\alpha} \right]} 
    \\
    \triangleq
    &\frac{I_1}{I_2-I_1},
    \ 
    \tau\rightarrow 0,  
\end{aligned}
\end{equation*} 
where $I_1 = \frac{1}{2-\alpha}[(\sigma+1)^{2-\alpha}-\sigma^{2-\alpha}]-\frac{1}{2}[(\sigma+1)^{1-\alpha}-\sigma^{1-\alpha}]$ and $I_2 = (\sigma+1)^{1-\alpha}$. It is evident that $I_1>0$ and $I_2>0$. Note that $2\sigma = 2-\alpha$. We have $I_1 = \frac{1}{2\sigma}(\sigma+1)^{1-\alpha}$ and $I_2 - I_1 = \frac{2\sigma-1}{2\sigma}(\sigma+1)^{1-\alpha}$. Since $\sigma\in(\frac{1}{2}, 1)$, there holds  
\begin{equation*}
    \displaystyle\lim_{\tau \rightarrow 0}\frac{(2\sigma-1)c^{(\alpha,\sigma)}_{2,1}}{\sigma c^{(\alpha,\sigma)}_{1,1}}=\frac{(2\sigma-1)I_1}{\sigma(I_2-I_1)}=\frac{1}{\sigma}>1,
\end{equation*}
which indicates $(2\sigma-1)c^{(\alpha,\sigma)}_{2,1}> \sigma c^{(\alpha,\sigma)}_{1,1}$ for sufficiently small $\tau$.

Next, we prove the case with $k\geq 2$. Using almost the same technique as that for $k=1$ yields that  
\begin{equation*} 
     \frac{c_{k+1,k}^{(\alpha,\sigma)}}{ c_{k,k}^{(\alpha,\sigma)}} 
    =\frac{\log\frac{t_k}{t_{k-1}}}{\log\frac{t_{k+1}}{t_k}} \frac{b_{k,k}^{(\alpha,\sigma)}+\left(\log\frac{t_{k+\sigma}}{t_k}\right)^{1-\alpha}}{a_{k,k}^{(\alpha,\sigma)}+b_{k-1,k}^{(\alpha,\sigma)}-b_{k,k}^{(\alpha,\sigma)}} 
    \approx
    \frac{I_1}
    {I_3-I_1},
    \ 
    \tau\rightarrow 0,  
\end{equation*} 
where $I_3 = \frac{1}{2-\alpha}[(\sigma+2)^{2-\alpha}-(\sigma+1)^{2-\sigma}]-\frac{1}{2}[(\sigma+2)^{1-\alpha}-(\sigma+1)^{1-\alpha}]>0$. Using $2\sigma = 2 - \alpha$, we have $(3\sigma-1)I_1-\sigma I_3 = \frac{1}{2\sigma}[(4\sigma-1)(\sigma+1)^{2\sigma-1}-2\sigma(\sigma+2)^{2\sigma-1}]>0$ via preliminary computations. Therefore, 
\begin{equation*}
    \displaystyle\lim_{\tau \rightarrow 0}\frac{(2\sigma-1)c^{(\alpha,\sigma)}_{k+1,k}}{\sigma c^{(\alpha,\sigma)}_{k,k}}=\frac{(2\sigma-1)I_1}{\sigma(I_3-I_1)}>1, 
    \ 
    k\geq 2. 
\end{equation*}
That is, $(2\sigma-1)c^{(\alpha,\sigma)}_{k+1,k}> \sigma c^{(\alpha,\sigma)}_{k,k}$ for sufficiently small $\tau$. The proof is thus completed. 
\end{proof}
\end{lemma} 

\begin{lemma}\label{lemma:1/c_kk}
Let $\alpha \in (0,1)$, $\sigma = 1-\frac{\alpha}{2}$, and the coefficients $c^{(\alpha,\sigma)}_{i,k}$ be defined by Eqs. \eqref{eq:c_{j,k}_1} and \eqref{eq:c_{j,k}}. There holds ${\displaystyle \frac{1}{c^{(\alpha,\sigma)}_{1,k}} < \frac{2\Gamma(1-\alpha)}{\tilde{a}^{\alpha}}} \left((k+\sigma)\tau\right)^{\alpha}$ for $k \geq 0$.
\begin{proof} 
The definition of $c_{1,0}^{(\alpha,\sigma)}$ yields that
\begin{equation*}
\begin{split}
    c_{1,0}^{(\alpha,\sigma)}
    =\frac{\left(\log\frac{t_{\sigma}}{t_0}\right)^{1-\alpha}}{\Gamma(2-\alpha)\log\frac{t_1}{t_0}}    
    >\frac{\left(\log\frac{t_{\sigma}}{t_0}\right)^{-\alpha}}{2\Gamma(2-\alpha)}
    >\frac{\tilde{a}^{\alpha}}{2\Gamma(2-\alpha)\left(\sigma\tau\right)^{\alpha}} 
    >\frac{\tilde{a}^{\alpha}}{2\Gamma(1-\alpha)\left(\sigma\tau\right)^{\alpha}}
    >0, 
\end{split}
\end{equation*}
where $\displaystyle2\log\frac{t_{\sigma}}{t_0}>\log\frac{t_1}{t_0}$ is used. Therefore, we have
\begin{equation*}
    \frac{1}{c^{(\alpha,\sigma)}_{1,0}} < \frac{2\Gamma(1-\alpha)\left(\sigma\tau\right)^{\alpha}}{\tilde{a}^{\alpha}}.
\end{equation*}

The definition of $c_{1,k}^{(\alpha,\sigma)}\,(k\geq 1)$ yields that 
\begin{equation}\label{eq:cjk>0}
\begin{split}
    c^{(\alpha,\sigma)}_{1,k}  
    =& \frac{1}{\Gamma(2-\alpha)}\frac{1}{\log\frac{t_1}{t_{0}}}\Bigg\{
    \left(\log\frac{t_{k+\sigma}}{t_{0}}\right)^{1-\alpha} - \left(\log\frac{t_{k+\sigma}}{t_{1}}\right)^{1-\alpha} 
    \\
    &- \left\{
     \frac{2}{2-\alpha}\left[\left(\log\frac{t_{k+\sigma}}{t_{0}}\right)^{2-\alpha}-\left(\log\frac{t_{k+\sigma}}{t_1}\right)^{2-\alpha}\right]\right. 
    \\
    &-\left.\left.\log\frac{t_1}{t_{0}}\left[\left(\log\frac{t_{k+\sigma}}{t_{1}}\right)^{1-\alpha} + \left(\log\frac{t_{k+\sigma}}{t_{0}}\right)^{1-\alpha} \right]
    \right\}
    \frac{1}{\log\frac{t_{2}}{t_{0}}}\right\} 
    \\ 
    =&\frac{1}{\Gamma(2-\alpha)}\frac{1}{\log\frac{t_1}{t_{0}}}\frac{\log\frac{t_{2}}{t_1}}{\log\frac{t_2}{t_0}}\left\{\left(\log\frac{t_{k+\sigma}}{t_0}\right)^{1-\alpha} -\left(\log\frac{t_{k+\sigma}}{t_1}\right)^{1-\alpha} \right\}
    \\
    &+\frac{1}{\Gamma(2-\alpha)}\frac{1}{\log\frac{t_1}{t_{0}}}\frac{1-\alpha}{\log\frac{t_{2}}{t_{0}}}
    \xi^{-\alpha}\left(\log\frac{t_{1}}{t_0}\right)^{2} 
    \\
    > & \frac{1}{\Gamma(1-\alpha)}\frac{\log\frac{t_{1}}{t_0}}{\log\frac{t_{2}}{t_{0}}}\xi^{-\alpha}
    > \frac{1}{2\Gamma(1-\alpha)}\frac{\log\left(\frac{t_{1}}{t_0}\right)^2}{\log\frac{t_{2}}{t_{0}}}\left(\log\frac{t_{k+\sigma}}{t_{0}}\right)^{-\alpha}
    \\
    >& \frac{\left(\frac{(k+\sigma)\tau}{\tilde{a}}\right)^{-\alpha}}{2\Gamma(1-\alpha)}, 
    \ 
    \xi\in\left(\log\frac{t_{k+\sigma}}{t_{1}},\log\frac{t_{k+\sigma}}{t_0}\right),  
\end{split}
\end{equation}
which indicates that  
\begin{equation*}
    \frac{1}{c^{(\alpha,\sigma)}_{1,k}} < \frac{2\Gamma(1-\alpha)\left((k+\sigma)\tau\right)^{\alpha}}{\tilde{a}^{\alpha}},
    \ 
    k\geq 1. 
\end{equation*}
All this ends the proof. 
\end{proof}
\end{lemma}

\subsection{L1-2 formula for Caputo-Hadamard derivative}
Let $t_k = \tilde{a} + k\tau$ with $k = 0,1,\cdots ,N\in\mathbb{Z}^{+}$, where $\tau = (T-\tilde{a})/N$ is the temporal stepsize. For $\varphi(t)\in C^3[\tilde{a},T]$, its Caputo-Hadamard derivative of order $\alpha\in(0,1)$ at the integer grid point $t = t_k$ can be evaluated by the following L1-2 approximation \cite{Fan2022},
\begin{equation}\label{eq:L12_CH}
    \begin{split}
       \left.{}_{CH}{\rm D}_{\tilde{a},t}^\alpha \varphi(t)\right\rvert_{t=t_k}
      \approx
      \sum\limits_{i=1}^k c^{(\alpha)}_{i,k}\left[\varphi(t_i)-\varphi(t_{i-1})\right] 
    \end{split}
\end{equation}
with the truncated error being of $\mathcal{O}(\tau^{3-\alpha})$. When $k=1, 2$, the coefficients are given by 
\begin{equation}\label{eq:c_{j,k}_2}
\begin{split}
&c_{1,1}^{(\alpha)}
=\frac{1}{\Gamma(2-\alpha)}\frac{1}{\log\frac{t_{1}}{t_{0}}}
a_{1,1}^{(\alpha)},\ \  k=1,
\\
&c_{1,2}^{(\alpha)}
=\frac{1}{\Gamma(2-\alpha)}\frac{1}{\log\frac{t_{1}}{t_{0}}}
\left(a_{1,2}^{(\alpha)}+b_{2,2}^{(\alpha)}\right),\ \
c_{2,2}^{(\alpha)}
=\frac{1}{\Gamma(2-\alpha)}\frac{1}{\log\frac{t_{2}}{t_{1}}}
\left(a_{2,2}^{(\alpha)}-b_{2,2}^{(\alpha)}\right),\ \ k=2.
\end{split}
\end{equation}
For the case with $k\geq 3$, the coefficients are defined by 
\begin{equation}\label{eq:c_{j,k}2} 
c^{(\alpha)}_{i,k}
=\left\{
\begin{array}{ll}
{\displaystyle\frac{1}{\Gamma(2-\alpha)}\frac{1}{\log\frac{t_{1}}{t_{0}}}}
\left(a^{(\alpha)}_{1,k}+b^{(\alpha)}_{2,k}\right),
&
i=1,
\\
{\displaystyle\frac{1}{\Gamma(2-\alpha)}\frac{1}{\log\frac{t_{i}}{t_{i-1}}}}
\left(a^{(\alpha)}_{i,k}-b^{(\alpha)}_{i,k}+b^{(\alpha)}_{i+1,k}\right), 
&
2\leq i\leq k-1,
\\
{\displaystyle\frac{1}{\Gamma(2-\alpha)}\frac{1}{\log\frac{t_{k}}{t_{k-1}}}}
\left(a^{(\alpha)}_{k,k}-b^{(\alpha)}_{k,k}\right),
&
i=k.
\end{array}
\right. 
\end{equation}
Here
\begin{equation*}
\begin{split}
    a^{(\alpha)}_{i,k}
    =&\left(\log\frac{t_{k}}{t_{i-1}}\right)^{1-\alpha}
-\left(\log\frac{t_{k}}{t_{i}}\right)^{1-\alpha},
\end{split}
\end{equation*}
and 
\begin{equation*}
\begin{split}
    b^{(\alpha)}_{i,k}
=\,&\Bigg\{
\log\frac{t_{i}}{t_{i-1}}
\left[\left(\log\frac{t_{k}}{t_{i}}\right)^{1-\alpha}
+\left(\log\frac{t_{k}}{t_{i-1}}\right)^{1-\alpha}\right]
\\
&+\frac{2}{2-\alpha}
\left[\left(\log\frac{t_{k}}{t_{i}}\right)^{2-\alpha}
-\left(\log\frac{t_{k}}{t_{i-1}}\right)^{2-\alpha}\right]
\Bigg\}
\frac{1}{\log\frac{t_i}{t_{i-2}}}.
\end{split}
\end{equation*} 

In the following, we introduce the important properties of coefficients $c_{i,k}^{(\alpha)}$ which will be used later on.  
\begin{lemma}{\rm\cite{Fan2022}}\label{lemma:lc_kk2}
For $\alpha \in (0,1)$ and sufficiently small temporal stepsize $\tau$, the coefficients $c_{i,k}^{(\alpha)}$ in Eqs. \eqref{eq:c_{j,k}_2} and \eqref{eq:c_{j,k}2} satisfy
\begin{description}
    \item[(1)] $c_{2,2}^{(\alpha)}>0$, 
    \item[(2)] $c_{3,3}^{(\alpha)}>c_{1,3}^{(\alpha)}>0$, 
    \item[(3)] $c_{k,k}^{(\alpha)}>c_{k-2,k}^{(\alpha)}>c_{k-3,k}^{(\alpha)}>\cdots>c_{1,k}^{(\alpha)}>0, \,k\geq 4$,
    \item[(4)] $c_{k,k}^{(\alpha)}>\left|c_{k-1,k}^{(\alpha)}\right|, \, k\geq 2$.
\end{description} 
\end{lemma}

\begin{lemma}\label{lemma:lc_kk_2}
For $\alpha \in (0,0.3738)$ and sufficiently small temporal stepsize $\tau$, the coefficients $c_{i,k}^{(\alpha)}$ given by Eqs. \eqref{eq:c_{j,k}_2} and \eqref{eq:c_{j,k}2} satisfy 
\begin{equation*}
    c^{(\alpha)}_{k,k}
    >c^{(\alpha)}_{k-1,k}
    >\cdots
    >c^{(\alpha)}_{i,k}
    >c^{(\alpha)}_{i-1,k}
    >\cdots
    >c^{(\alpha)}_{1,k}>0, 
    \, 
    k\geq 2.
\end{equation*}
\begin{proof}
We first show that $c_{1,2}^{(\alpha)}>0$ if $\alpha\in(0,0.4701)$. It follows from the definition that  
\begin{equation*}
\begin{split}
    c_{1,2}^{(\alpha)}  
    =&\frac{1}{\Gamma(3-\alpha)}\frac{1}{\log\frac{t_{1}}{t_{0}}\log\frac{t_{2}}{t_{0}}} 
    \Bigg\{\log\frac{t_{2}}{t_{0}} \left[(2-2\alpha)\left(\log\frac{t_{2}}{t_{0}}\right)^{1-\alpha}  -(2-\alpha)\left(\log\frac{t_{2}}{t_{1}}\right)^{1-\alpha} \right] 
    \\
    &+\alpha \left[\left(\log\frac{t_{2}}{t_{0}}\right)^{2-\alpha} -\left(\log\frac{t_{2}}{t_{1}}\right)^{2-\alpha}\right]  \Bigg\}
    \\
    >&\frac{1}{\Gamma(3-\alpha)}\frac{1}{\log\frac{t_{1}}{t_{0}}}\left[(2-2\alpha)\left(\log\frac{t_{2}}{t_{0}}\right)^{1-\alpha}  -(2-\alpha)\left(\log\frac{t_{2}}{t_{1}}\right)^{1-\alpha} \right] .
\end{split}
\end{equation*}
Using the equivalent infinitesimal substitution technique gives
\begin{equation*}
\begin{split}
    &{(2-2\alpha)\left(\log\frac{t_{2}}{t_{0}}\right)^{1-\alpha}}-{(2-\alpha)\left(\log\frac{t_{2}}{t_{1}}\right)^{1-\alpha}}  
    \approx  
    \left(\frac{\tau}{t_{0}}\right)^{1-\alpha}\left[{(2-2\alpha)2^{1-\alpha}}-{(2-\alpha)}\right],
    \ 
    \tau\rightarrow 0.
\end{split}
\end{equation*}
Through numerical computation, one has ${(2-2\alpha)2^{1-\alpha}}-{(2-\alpha)} > 0$ if $\alpha\in(0,0.4701)$. So we have $c_{1,2}^{(\alpha)}>0$ with sufficiently small temporal stepsize $\tau$ and $\alpha\in(0,0.4701)$. 

Then we show that $c_{2,3}^{(\alpha)}>c_{1,3}^{(\alpha)}$ if $\alpha\in(0,0.3909)$ and $\tau$ is sufficiently small. 
Note that
\begin{equation*}
\begin{split}
    a^{(\alpha)}_{1,3} 
    \approx 
    \left(\frac{\tau}{t_1}\right)^{1-\alpha}\left(3^{1-\alpha} -2^{1-\alpha}\right),
    \ 
    \tau\rightarrow 0,
\end{split}
\end{equation*}
\begin{equation*} 
    a^{(\alpha)}_{2,3} 
    \approx
    \left(\frac{\tau}{t_1}\right)^{1-\alpha}\left(2^{1-\alpha} - 1\right),
    \ 
    \tau\rightarrow 0, 
\end{equation*}
\begin{equation*} 
    b^{(\alpha)}_{2,3}
    \approx
    \left(\frac{\tau}{t_1}\right)^{1-\alpha}
    \left(\frac{1+2^{1-\alpha}}{2}
    +\frac{1-2^{2-\alpha}}{2-\alpha}\right),
    \ 
    \tau\rightarrow 0, 
\end{equation*}
and 
\begin{equation*} 
    b^{(\alpha)}_{3,3} 
    \approx 
    -\left(\frac{\tau}{t_1}\right)^{1-\alpha}\frac{\alpha}{2(2-\alpha)},
    \ 
    \tau\rightarrow 0.  
\end{equation*}
We arrive at 
\begin{equation*} 
    c^{(\alpha)}_{2,3}-c^{(\alpha)}_{1,3} 
    \approx
    \left(\frac{\tau}{t_1}\right)^{2-\alpha}
    \left( 
    -3^{1-\alpha}+2^{1-\alpha}+\frac{2^{3-\alpha}-3}{2-\alpha} -\frac{3}{2}\right),
    \ 
    \tau\rightarrow 0.  
\end{equation*}
Numerical computation yields that $\displaystyle -3^{1-\alpha}+2^{1-\alpha}+\frac{2^{3-\alpha}-3}{2-\alpha} -\frac{3}{2} > 0$ if $\alpha\in(0,0.3909)$. Hence $c_{2,3}^{(\alpha)}>c_{1,3}^{(\alpha)}$ holds for sufficiently small temporal stepsize $\tau$ and $\alpha\in(0,0.3909)$. 

We finally show that $c_{k-1,k}^{(\alpha)}>c_{k-2,k}^{(\alpha)}\, (k\geq 4)$ if $\alpha\in(0,0.3738)$ and $\tau$ is sufficiently small. Note that  
\begin{equation*} 
    a^{(\alpha)}_{k-2,k} 
    \approx
    \left(\frac{\tau}{t_{k-2}}\right)^{1-\alpha}\left(3^{1-\alpha} - 2^{1-\alpha}\right),
    \ 
    \tau\rightarrow 0, 
\end{equation*}
\begin{equation*} 
    a^{(\alpha)}_{k-1,k} 
    \approx
    \left(\frac{\tau}{t_{k-2}}\right)^{1-\alpha}\left(2^{1-\alpha} - 1\right),
    \ 
    \tau\rightarrow 0, 
\end{equation*}
\begin{equation*} 
    b^{(\alpha)}_{k-2,k} 
    \approx
    \left(\frac{ \tau}{t_{k-2}}\right)^{1-\alpha}
    \left(\frac{2^{1-\alpha}+3^{1-\alpha}}{2}
    +\frac{2^{2-\alpha}-3^{2-\alpha}}{2-\alpha}\right),
    \ 
    \tau\rightarrow 0, 
\end{equation*}
\begin{equation*} 
    b^{(\alpha)}_{k-1,k} 
    \approx
    \left(\frac{ \tau}{t_{k-2}}\right)^{1-\alpha}
    \left(\frac{1 +2^{1-\alpha}}{2}
    + \frac{1-2^{2-\alpha}}{2-\alpha}\right),
    \ 
    \tau\rightarrow 0, 
\end{equation*}
and 
\begin{equation*} 
    b^{(\alpha)}_{k,k} 
    \approx
    -\left(\frac{\tau}{t_{k-2}}\right)^{1-\alpha}\frac{\alpha}{2(2-\alpha)},
    \ 
    \tau\rightarrow 0. 
\end{equation*}
Therefore, 
\begin{equation*} 
    c^{(\alpha)}_{k-1,k}-c^{(\alpha)}_{k-2,k} 
    \approx
    \left(\frac{\tau}{t_{k-2}}\right)^{2-\alpha}
    \left(\frac{-3^{1-\alpha}+3\cdot2^{1-\alpha}-3}{2}+\frac{-3^{2-\alpha}+3\cdot2^{2-\alpha}-3}{2-\alpha}\right),
    \ 
    \tau\rightarrow 0.  
\end{equation*}
Applying the numerical computation, one gets 
$$\displaystyle \frac{-3^{1-\alpha}+3\cdot2^{1-\alpha}-3}{2}+\frac{-3^{2-\alpha}+3\cdot2^{2-\alpha}-3}{2-\alpha} > 0$$
if $\alpha\in(0,0.3738)$. Hence $c_{k-1,k}^{(\alpha)}>c_{k-2,k}^{(\alpha)}$ is valid for sufficiently small $\tau$ and $\alpha\in(0,0.3738)$.  

Taking Lemma \ref{lemma:lc_kk2} into account, we can find that 
\begin{equation*}
    c^{(\alpha)}_{k,k}
    >c^{(\alpha)}_{k-1,k}
    >\cdots
    >c^{(\alpha)}_{i,k}
    >c^{(\alpha)}_{i-1,k}
    >\cdots
    >c^{(\alpha)}_{1,k}>0, 
    \ 
    k\geq 2,
\end{equation*}
if $\alpha \in (0,0.3738)$ and the temporal stepsize $\tau$ is sufficiently small, 
\end{proof}
\end{lemma}

\begin{lemma}\label{lemma:1/c_kk_alpha}
For $\alpha\in(0,0.4701)$ and sufficiently small temporal stepsize $\tau$, the coefficients $c^{(\alpha)}_{i,k}$ defined by Eqs. \eqref{eq:c_{j,k}_2} and \eqref{eq:c_{j,k}2} satisfy the following estimate, 
\begin{equation}\label{eq:1/c_kk_alpha}
    \frac{1}{c^{(\alpha)}_{1,k}} <C^{(\alpha)} \left(k\tau\right)^{\alpha},\ k\geq 1, 
\end{equation}
where $C^{(\alpha)} =  \max\left\{\displaystyle\frac{6\Gamma(1-\alpha)}{\tilde{a}^{\alpha}}, \frac{3\Gamma(2-\alpha)}{\alpha\tilde{a}^{\alpha}}\right\}$. 
\begin{proof}  
For the case with $k=1$,  
\begin{equation*}
    c_{1,1}^{(\alpha)}=\frac{1}{\Gamma(2-\alpha)}\frac{\left(\log\frac{t_{1}}{t_{0}}\right)^{1-\alpha}}{\log\frac{t_{1}}{t_{0}}}=
    \frac{\left(\log\left(1+\frac{\tau}{\tilde{a}}\right)\right)^{-\alpha}}{\Gamma(2-\alpha)}>\frac{\left(\frac{\tau}{\tilde{a}}\right)^{-\alpha}}{\Gamma(2-\alpha)}, 
\end{equation*}
which yields that  
\begin{equation*}
    \frac{1}{c_{1,1}^{(\alpha)}}<\frac{\Gamma(2-\alpha)}{\tilde{a}^{\alpha}}\tau^{\alpha}.
\end{equation*}

Now we consider the case with $k=2$. It follows from the proof of Lemma \ref{lemma:lc_kk_2} and the mean value theorem that if $\alpha\in (0,0.4701)$ and the temporal stepsize $\tau$ is sufficiently small, 
\begin{equation*} 
    c_{1,2}^{(\alpha)} 
    > \frac{\alpha}{\Gamma(3-\alpha)}\frac{1}{\log\frac{t_{1}}{t_{0}}\log\frac{t_{2}}{t_{0}}}  \left[\left(\log\frac{t_{2}}{t_{0}}\right)^{2-\alpha} -\left(\log\frac{t_{2}}{t_{1}}\right)^{2-\alpha}\right] 
    \\ 
    = \frac{\alpha}{\Gamma(2-\alpha)}\frac{1}{\log\frac{t_2}{t_{0}}}
    \xi^{1-\alpha}
\end{equation*} 
holds, where $\xi\in\left(\log\frac{t_2}{t_{1}}, \log\frac{t_2}{t_{0}}\right)$.  Recall that $t_0 = \tilde{a}$. We have 
\begin{equation*} 
\begin{split}
    c^{(\alpha)}_{1,2}
    >& \frac{\alpha}{\Gamma(2-\alpha)}\frac{\left(\log\frac{t_2}{t_{1}}\right)^{1-\alpha}}{\log\frac{t_2}{t_{0}}}
    >\frac{\alpha\left(\log\frac{t_2}{t_{1}}\right)^{-\alpha}}{3\Gamma(2-\alpha)}
    =\frac{\alpha\left(\log(1+\frac{\tau}{\tilde{a}+\tau})\right)^{-\alpha}}{3\Gamma(2-\alpha)}
    > \frac{\alpha\left(\frac{\tau}{\tilde{a}}\right)^{-\alpha}}{3\Gamma(2-\alpha)}
\end{split}
\end{equation*}
because of $\displaystyle3\log\frac{t_2}{t_{1}}>\log\frac{t_2}{t_{0}}$
with $t_0=\tilde{a}>\tau$. In other words, 
\begin{equation*}
    \frac{1}{c^{(\alpha)}_{1,2}} < \frac{3\Gamma(2-\alpha)\tau^{\alpha}}{\alpha\tilde{a}^{\alpha}}. 
\end{equation*} 

We finally consider the case with $k\geq 3$. The Taylor formula yields that 
\begin{equation*}
\begin{split}
    c^{(\alpha)}_{1,k}=\frac{1}{\Gamma(2-\alpha)}\frac{1}{\log\frac{t_{1}}{t_{0}}}\left[\frac{1}{\log\frac{t_2}{t_{0}}} J_1
    +\frac{(1-\alpha)\left(\log\frac{t_2}{t_1}\right)^2}{\log\frac{t_2}{t_0}}\eta^{-\alpha}\right], 
\end{split}
\end{equation*}
with $\eta\in\left(\log\frac{t_k}{t_{2}}, \log\frac{t_k}{t_{1}}\right)$ and 
\begin{equation*}
    J_1 = \log\frac{t_2}{t_{0}}\left[\left(\log\frac{t_k}{t_0}\right)^{1-\alpha}-\left(\log\frac{t_k}{t_1}\right)^{1-\alpha}\right]+\log\frac{t_2}{t_{1}}\left[\left(\log\frac{t_k}{t_2}\right)^{1-\alpha}-\left(\log\frac{t_k}{t_1}\right)^{1-\alpha}\right].
\end{equation*}
From \cite{Fan2022}, we know that $J_1>0$ when $k\geq 3$. Therefore, 
\begin{equation*}
    c^{(\alpha)}_{1,k}
    >\frac{1}{\Gamma(1-\alpha)}\frac{\left(\log\frac{t_2}{t_1}\right)^2}{\log\frac{t_{1}}{t_{0}}\log\frac{t_2}{t_0}}\eta^{-\alpha}, 
    \ 
    \eta\in\left(\log\frac{t_k}{t_{2}},\, \log\frac{t_k}{t_{1}}\right).
\end{equation*}
Due to the inequalities $2\log\frac{t_2}{t_{1}}>\log\frac{t_1}{t_{0}}$ and $3\log\frac{t_2}{t_{1}}>\log\frac{t_2}{t_{0}}$
as $t_0=\tilde{a}>\tau$ when $\tau$ is sufficiently small, we have 
\begin{equation*}
    c^{(\alpha)}_{1,k}
    > \frac{\left(\log\frac{t_k}{t_{1}}\right)^{-\alpha}}{6\Gamma(1-\alpha)} 
    = \frac{\left(\log\left(1+\frac{(k-1)\tau}{\tilde{a}+\tau}\right)\right)^{-\alpha}}{6\Gamma(1-\alpha)}
    > \frac{\left(\frac{k\tau}{\tilde{a}}\right)^{-\alpha}}{6\Gamma(1-\alpha)}.
\end{equation*}
That is to say, 
\begin{equation*}
    \frac{1}{c^{(\alpha)}_{1,k}} < \frac{6\Gamma(1-\alpha)}{\tilde{a}^{\alpha}} \left(k\tau\right)^{\alpha},\ k\geq 3. 
\end{equation*}

Based on the above analysis, the desired result is proved. 
\end{proof}
\end{lemma}
Next, we discuss the discretizations for the Riesz derivative and the fractional Laplacian used in space.
\subsection{Weighted and shifted Gr\"{u}nwald-Letnikov formula for Riesz derivative}
Let $x_j = a + jh\ (j = 0,1,\cdots,M)$ with $h = (b-a)/M$ being the  spatial stepsize. Define the spaces of grid functions as $\mathcal{U}_h = \{w\ \lvert\ w = (w_0,w_1,\cdots,w_M)\}$ and $\mathcal{U}_h^{\circ} = \{w\ \lvert\ w\in \mathcal{U}_h,w_0 = w_M =0 \}$. 

For any grid functions $w, v\in\mathcal{U}_h^{\circ}$, a discrete inner product and the associated norm are defined as usual, i.e.,
\begin{equation*}
    (w,v)_{h} = h \sum_{j=1}^{M-1}w_{j}v_{j}, \ \ 
    \Vert w \Vert_{h}^2 = (w,w)_{h}.    
\end{equation*}

If $\psi(x)$, $_{RL}{\rm D}_{a,x}^\beta \psi(x)$, $_{RL}{\rm D}_{x,b}^\beta \psi(x)$, and their Fourier transforms belong to $L^{1}(\mathbb{R})$, then the Riesz derivative of $\psi(x)$ at $x=x_j\,(1\leq j \leq M-1)$ can be approximated by the following weighted and shifted Gr\"{u}nwald-Letnikov formula with 2-nd order accuracy \cite{Tian2015},
\begin{equation}\label{eq:RZ1}
\begin{split}
    \left.{}_{RZ}{\rm D}_{x}^\beta\psi\left(x\right)\right|_{x_{j}}
    &=-\frac{\Psi_{\beta}}{h^{\beta}}
    \bigg [v_{1}\sum_{k=0}^{j+l_{1}}g^{(\beta)}_{k}\psi\left(x_{j-k+l_{1}}\right)
    +v_{2}\sum_{k=0}^{j+l_{2}}g^{(\beta)}_{k}\psi\left(x_{j-k+l_{2}}\right)
    \\
    +&v_{1}\sum_{k=0}^{M-j+l_{1}}g^{(\beta)}_{k}\psi(x_{j+k-l_{1}})
    +v_{2}\sum_{k=0}^{M-j+l_{2}}g^{(\beta)}_{k}\psi(x_{j+k-l_{2}})\bigg]
    +\mathcal{O}(h^2),
\end{split}
\end{equation}
where $\displaystyle \Psi_{\beta} = \frac{1}{2\cos(\frac{\pi\beta}{2})},\   g^{(\beta)}_{k}=(-1)^{k}\binom{\beta}{k}$, $\displaystyle v_{1}=\frac{\beta-2l_{2}}{2(l_{1}-l_{2})}$, and $\displaystyle v_{2}=\frac{2l_{1}-\beta}{2(l_{1}-l_{2})}$ with $l_{1}\neq l_{2}$. In particular, the coefficients $g^{(\beta)}_{k}$ can be computed via the following recursive formula, 
\begin{equation}\label{eq:gkgk_1}
    g^{(\beta)}_{0}=1,\ 
    g^{(\beta)}_{k}=(1-\frac{\beta+1}{k})g^{(\beta)}_{k-1}, \ 
    k = 1,2,\cdots.
\end{equation}

\begin{lemma}{\rm\cite{Meerschaert2006,Podlubny1999,Tian2015}}\label{lemma:g_k}
For  $1< \beta <2$, the coefficients $g^{(\beta)}_{k}$ in Eq. {\rm (\ref{eq:RZ1})} satisfy:
\begin{equation}\label{eq:g_k}
    \left\{
    \begin{array}{ll}
    g^{(\beta)}_{0} = 1, g^{(\beta)}_{1} = -\beta,
    \\
    1\geq g^{(\beta)}_{2}\geq g^{(\beta)}_{3}\geq \cdots\geq 0,
    \\
    \sum\limits_{k=0}\limits^{\infty}g^{(\beta)}_{k} = 0,  \sum\limits_{k=0}\limits^{m}g^{(\beta)}_{k} < 0,
    \ 
    m \geq 1.
    \end{array}
    \right.
\end{equation}
\end{lemma} 

As a matter of fact, if $l_1 = 1$ and $l_2 =0$, Eq. \eqref{eq:RZ1} can be rewritten as \cite{WangCai2023}
\begin{equation}\label{eq:RZ3}
\begin{split}
    _{RZ}{\rm D}_{x}^\beta \psi\left(x_{j}\right)
    &=-\frac{1}{h^{\beta}}
    \left(\sum_{k=0}^{j}r^{(\beta)}_{j-k}\psi\left(x_{k}\right)
    +\sum_{k=j+1}^{M}r^{(\beta)}_{k-j}\psi\left(x_{k}\right) \right)
    +\mathcal{O}(h^2)
    \\
    &=-\frac{1}{h^{\beta}}\sum_{k=0}^{M}r^{(\beta)}_{j-k}\psi\left(x_{k}\right) 
    +\mathcal{O}(h^2),
\end{split}
\end{equation}
where 
\begin{equation}\label{eq:r_k}
    \left\{
    \begin{array}{ll}
    \displaystyle
    r^{(\beta)}_{0} 
    = 2\Psi_{\beta}\left(\frac{\beta}{2}g_{1}^{(\beta)}+\frac{2-\beta}{2}g_0^{(\beta)}\right)
    = \frac{1}{\cos(\frac{\pi\beta}{2})}
    \left(1-\frac{\beta + \beta^2}{2}\right)>0,
    \\
    \displaystyle
    r^{(\beta)}_{1} 
    = \Psi_{\beta}\left(\frac{\beta}{2}g_{0}^{(\beta)}+\frac{2-\beta}{2}g_{1}^{(\beta)}+\frac{\beta}{2}g_{2}^{(\beta)}\right)
    = \frac{1}{2\cos(\frac{\pi\beta}{2})}\frac{\beta(\beta^2 + \beta -2)}{4}<0 ,
    \\
    \displaystyle
    r^{(\beta)}_{k} 
    = \Psi_{\beta}\left(\frac{\beta}{2}g_{k+1}^{(\beta)}+\frac{2-\beta}{2}g_k^{(\beta)}\right) <0,
    \ 
    k\geq 2,
    \\
    r^{(\beta)}_{k} = r^{(\beta)}_{-k}<0,
    \ 
    k\geq 1.
    \end{array}
    \right.
\end{equation}

The following lemma is fundamental and important.
\begin{lemma}\label{lemma:lr_k} 
Let $1< \beta <2$. The coefficients $r^{(\beta)}_{k}$ given by Eq. \eqref{eq:r_k} satisfy 
\begin{equation*}
    \left\{
    \begin{array}{ll}
    \sum\limits_{k=0}\limits^{m}r^{(\beta)}_{j-k} > 0,\ \ m-1\geq j\geq 1,
    \, 
    m\geq 2,
    \\ 
    \sum\limits_{k=1-m }\limits^{m-1}r_{k}^{(\beta)} > \displaystyle\frac{c_{\ast}^{(\beta)}}{m^{\beta}}>0,\ \  m \geq 2,
    \end{array}
    \right.
\end{equation*}
where $\displaystyle c_{\ast}^{(\beta)} = \frac{(1-\beta)(2-\beta)(3-\beta)4^{\beta}e^{-\frac{9}{4}}\Psi_{\beta}}{3}>0$. 
\begin{proof}
The above assertions have been proved in Lemma 6 of \cite{WangCai2023} except for the positive lower bound of $\sum\limits_{k=1-m }\limits^{m-1}r_{k}^{(\beta)}$. In the following, we shall determine this lower bound which will be used in the corresponding numerical stability analysis. 

Note that $1-x>e^{-x-x^2}$ when $x\in(0, \frac{2}{3})$. For $k\geq 5$, we have 
\begin{equation}\label{eq:gkgk4} 
    g^{(\beta)}_{k}
    =\prod\limits_{m=5}^{k}\left(1-\frac{\beta+1}{m}\right)g_{4}^{(\beta)}
    \geq e^{-(\beta+1)\sum\limits_{m=5}\limits^{k}\frac{1}{m}-(\beta+1)^2\sum\limits_{m=5}\limits^{k}\frac{1}{m^2}}g^{(\beta)}_{4}>0. 
\end{equation}
Since $\displaystyle \sum\limits_{m=5}\limits^{k}\frac{1}{m} \leq {\log}\left(\frac{k}{4}\right)$ 
and $\displaystyle \sum\limits_{m=5}\limits^{k}\frac{1}{m^2} < \frac{1}{4}$, there holds 
\begin{equation*}
\begin{split}
    g^{(\beta)}_{k} \geq \left(\frac{4}{k}\right)^{\beta+1}e^{-\frac{9}{4}}g^{(\beta)}_{4} = c^{(\beta)}\left(\frac{1}{k}\right)^{\beta+1},\,\, k\geq 5,
\end{split}
\end{equation*}
where $\displaystyle c^{(\beta)} = 4^{\beta+1}e^{-\frac{9}{4}}g^{(\beta)}_{4}=\frac{-\beta(1-\beta)(2-\beta)(3-\beta)4^{\beta}e^{-\frac{9}{4}}}{6} >0$. As a matter of fact, direct calculation indicates that the above inequality is also valid for $k=2,3,4$. Therefore, we have  
\begin{equation*}
    g^{(\beta)}_{k} \geq c^{(\beta)}\left(\frac{1}{k}\right)^{\beta+1},\,\, k\geq 2, 
\end{equation*}
which implies 
\begin{equation*}
    \sum\limits_{k=l}\limits^{\infty} g^{(\beta)}_{k} \geq c^{(\beta)}\sum\limits_{k=l}\limits^{\infty}\left(\frac{1}{k}\right)^{\beta+1} 
    \geq  c^{(\beta)}\int_{l}^{\infty}\frac{1}{x^{\beta+1}}{\rm d}x 
    = \frac{c^{(\beta)}}{\beta l^{\beta}},\ \  l\geq 2.
\end{equation*} 
Recall that $\sum\limits_{k=0}\limits^{\infty}g^{(\beta)}_{k} = 0$ (see Eq. \eqref{eq:g_k}). It is evident that 
\begin{equation*}
    \sum\limits_{k=0}\limits^{m}g^{(\beta)}_{k} = \sum\limits_{k=0}\limits^{\infty}g^{(\beta)}_{k} - \sum\limits_{k=m+1}\limits^{\infty} g^{(\beta)}_{k} \leq -\frac{c^{(\beta)}}{\beta (m+1)^{\beta}} ,\ \  m\geq 1.
\end{equation*}
Combining the above inequality with Eq. \eqref{eq:g_k} yields 
\begin{equation*}
\begin{split}
    \sum\limits_{k=1-m }\limits^{m-1}r_{k}^{(\beta)} 
    &= 2\Psi_{\beta}\left(\frac{\beta}{2}\sum\limits_{k=0}\limits^{m}g^{(\beta)}_{k} +\frac{2-\beta}{2}\sum\limits_{k=0}\limits^{m-1}g^{(\beta)}_{k} \right) 
    \\
    &\geq 2\Psi_{\beta}\sum\limits_{k=0}\limits^{m-1}g^{(\beta)}_{k}
    \geq - 2\Psi_{\beta}\frac{c^{(\beta)}}{\beta{m}^{\beta}}
    =\frac{c_{\ast}^{(\beta)}}{m^{\beta}}>0 ,    \, \,m \geq 2, 
\end{split}
\end{equation*}
with $c_{\ast}^{(\beta)}= - 2\Psi_{\beta}\frac{c^{(\beta)}}{\beta}>0$. All this completes the proof. 
\end{proof}
\end{lemma} 

\subsection{Fractional centered difference formula for fractional Laplacian}
In  two dimensional case, let $h = 2L/M$ with $M\in\mathbb{Z}^{+}$, $x_j = -L + jh$ with $0\leq j \leq M$, and $y_k = -L + kh$ with $0\leq k \leq M$. Denote $\widetilde{\Omega}_h = \{(x_j,y_k)\ \lvert\ 0\leq j, k \leq M\}$,
$\Omega_h = \widetilde{\Omega}_h \cap \widetilde{\Omega}$ and $\partial\Omega_h = \widetilde{\Omega}_h \cap \partial\widetilde{\Omega}$.

For any grid function $w, v$ in the grid function space $S_h^{\circ} = \{w\ \lvert\  w=\{w_{jk}\}, w_{jk} = 0$ with $ (x_j,y_k) \in \partial\Omega_h$\}, a discrete inner product and the associated norm are defined as 
\begin{equation*}
    (w,v)_{h^2} = h^2 \sum_{j=1}^{M-1}\sum_{k=1}^{M-1}w_{jk}v_{jk}, \ \ 
    \Vert w \Vert_{L_h^2}^2 = (w,w)_{h^2}.    
\end{equation*}
Set $L_h^2 = \left\{w\ \lvert\  w= \{w_{jk}\}, \Vert w \Vert_{L_h^2}^2 < +\infty  \right\} $. For $w \in L_h^2$, we define the semi-discrete Fourier transform $\widehat{w}: [-\frac{\pi}{h},\frac{\pi}{h}]^2 \rightarrow \mathbb{C} $ as
\begin{equation*}
    \widehat{w}(\eta_1,\eta_2) = h^2 \sum_{j=1}^{M-1}\sum_{k=1}^{M-1}w_{jk}e^{-\mathbf{i}(\eta_1jh+\eta_2kh)},
\end{equation*}
and the inverse semi-discrete Fourier transform
\begin{equation*}
    w_{jk}=\frac{1}{4\pi^2}\int_{-\pi/h}^{\pi/h}\int_{-\pi/h}^{\pi/h}\widehat{w}(\eta_1,\eta_2)e^{-\mathbf{i}(\eta_1jh+\eta_2kh)}{\rm d}\eta_1{\rm d}\eta_2, 
\end{equation*}
with $\mathbf{i}$ being the imaginary unit. It follows from Parseval's identity that  the continuous definition of the inner product takes the form 
\begin{equation*}
\begin{split}
    (w,v)_{h^2}
    &=\frac{1}{4\pi^2}\int_{-\pi/h}^{\pi/h}\int_{-\pi/h}^{\pi/h}\widehat{w}(\eta_1,\eta_2)\widehat{v}(\eta_1,\eta_2){\rm d}\eta_1{\rm d}\eta_2,
\end{split}
\end{equation*}
with the norm given by 
\begin{equation*}
\begin{split} 
    \Vert w \Vert_{L_h^2}^2 
    &= \frac{1}{4\pi^2}\int_{-\pi/h}^{\pi/h}\int_{-\pi/h}^{\pi/h}
    \left\lvert\widehat{w}(\eta_1,\eta_2)\right\rvert^2{\rm d}\eta_1{\rm d}\eta_2.
\end{split}
\end{equation*}
For an arbitrary positive constant $s$, define the fractional Sobolev semi-norm $\lvert\cdot\rvert_{H_h^s}$ 
as 
\begin{equation*}
\begin{split}
    \lvert w\rvert_{H_h^s}^2
    &=\frac{1}{4\pi^2}\int_{-\pi/h}^{\pi/h}\int_{-\pi/h}^{\pi/h}\left(\eta_1^2+\eta_2^2 \right)^s\left\lvert\widehat{w}(\eta_1,\eta_2)\right\rvert^2{\rm d}\eta_1{\rm d}\eta_2. 
\end{split}
\end{equation*} 

In the above settings, we introduce the following two-dimensional fractional centered difference formula 
\begin{equation}\label{eq:Deltah}
    \left(-\Delta_h\right)^{\frac{\beta}{2}}\psi(x,y) = \frac{1}{h^{\beta}}\sum_{j,k\in \mathbb{Z}}a_{j,k}^{(\beta)}\psi(x+jh,y+kh),
\end{equation}
with 
\begin{equation*} 
    a_{j,k}^{(\beta)} = \frac{1}{4\pi^2}\int_{-\pi}^{\pi}\int_{-\pi}^{\pi}\left[4{\rm sin}^2(\frac{\eta_1}{2})+4{\rm sin}^2(\frac{\eta_2}{2}) \right]^{\frac{\beta}{2}}e^{-\mathbf{i}(\eta_1j+\eta_2k)} {\rm d}\eta_1{\rm d}\eta_2. 
\end{equation*}
In fact, the coefficients $a_{j,k}^{(\beta)}$ can be numerically computed by the built-in function ``fft2" in Matlab \cite{HaoZhang2021}. Throughout the numerical simulations in this paper, we compute $a_{j,k}^{(\beta)}$ as the way given by \cite{HaoZhang2021,WangCai2023}.  

Before introducing properties of the fractional central finite difference formula, the following space 
\begin{equation*}
    \mathcal{B}^s(\mathbb{R}^2)
    =\left\{\psi\left\lvert\ \psi\in L^1(\mathbb{R}^2), \int_{\mathbb{R}^2}\left[1+\lvert\eta \rvert^2\right]^s\lvert\widehat{\psi}(\eta_1,\eta_2)\rvert {\rm d}\eta_1{\rm d}\eta_2 < \infty \right.\right\}
\end{equation*}
will be used, where $\lvert\eta \rvert^2 = \eta_1^2+\eta_2^2$.

\begin{lemma}{\rm\cite{HaoZhang2021}}\label{lemma:ErrorRL}
Let $\psi(x,y)\in\mathcal{B}^{2+\beta}(\mathbb{R}^2)$. Its fractional Laplacian defined in Eq. \eqref{eq:lap} can be approximated by the fractional centered difference operator, i.e., 
\begin{equation*} 
    \left(-\Delta\right)^{\frac{\beta}{2}}\psi(x,y)
    =\left(-\Delta_h\right)^{\frac{\beta}{2}}\psi(x,y)+R_{L}(x,y),\ \ 0< \beta \leq 2,
\end{equation*}
where the truncation error satisfies 
\begin{equation}\label{eq:Rlap}
    \left\lvert R_{L}(x,y)\right\rvert 
    \leq 
    C_0 h^2\int_{\mathbb{R}^2}\left(1+\lvert\eta\rvert\right)^{\beta + 2}
    \left\lvert\widehat{\psi}(\eta_1,\eta_2)\right\rvert{\rm d}\eta_1{\rm d}\eta_2 = Ch^2 
\end{equation}
with $C$ being a constant independent of h. 
\end{lemma}
 
\begin{lemma}\label{lemma:psi} 
Let $\psi\in H_{h}^{\frac{\beta}{2}}(\mathbb{R}^2)$ with $\beta \in (1,2)$. It holds that 
\begin{equation*}
    C\left(\frac{2}{\pi}\right)^{\beta}\frac{2^{\frac{\beta-2}{2}}}{4\pi^2} \Vert \psi \Vert_{L_h^2}^2 
    \leq \left( \left(-\Delta_h\right)^{\frac{\beta}{2}}\psi,\psi \right)_{h^2}
    \leq 2^{\frac{\beta}{2}}\frac{\pi^{\beta}}{h^{\beta}} \Vert \psi \Vert_{L_h^2}^2,
\end{equation*}
where $C$ is a positive constant independent of $h$. 
\begin{proof}
It follows from the proof of Lemma 3.3 in \cite{HaoZhang2021} that 
\begin{equation}\label{eq:Parseval}
    \left( \left(-\Delta_h\right)^{\frac{\beta}{2}}\psi,\psi \right)_{h^2} 
    =\frac{1}{4\pi^2}\frac{1}{h^{\beta}} \int_{-\pi/h}^{\pi/h}\int_{-\pi/h}^{\pi/h}\left[4{\rm sin}^2(\frac{\eta_1h}{2})+4{\rm sin}^2(\frac{\eta_2h}{2}) \right]^{\frac{\beta}{2}}
    \left\lvert\widehat{\psi}(\eta_1,\eta_2)\right\rvert^2{\rm d}\eta_1{\rm d}\eta_2, 
\end{equation}
and 
\begin{equation}\label{eq:eta1eta2}
    \left(\frac{2}{\pi}\right)^{\beta}h^{\beta}\left(\eta_1^2 +\eta_2^2 \right)^{\frac{\beta}{2}}
    \leq \left[4{\rm sin}^2(\frac{\eta_1h}{2})+4{\rm sin}^2(\frac{\eta_2h}{2}) \right]^{\frac{\beta}{2}}
    \leq h^{\beta}\left(\eta_1^2 +\eta_2^2 \right)^{\frac{\beta}{2}}. 
\end{equation} 
Therefore,  
\begin{equation*}
\begin{aligned}
    \left( \left(-\Delta_h\right)^{\frac{\beta}{2}}\psi,\psi \right)_{h^2} 
    \leq 
    &\frac{1}{4\pi^2}\frac{1}{h^{\beta}} 
    \int_{-\pi/h}^{\pi/h}\int_{-\pi/h}^{\pi/h}
    h^{\beta}\left(\eta_1^2 +\eta_2^2 \right)^{\frac{\beta}{2}}
    \left\lvert\widehat{\psi}(\eta_1,\eta_2)\right\rvert^2{\rm d}\eta_1{\rm d}\eta_2
    \\
    \leq 
    &\frac{1}{4\pi^2}\frac{1}{h^{\beta}}\frac{h^{\beta}(2\pi^2)^{\frac{\beta}{2}}}{h^{\beta}}
    \int_{-\pi/h}^{\pi/h}\int_{-\pi/h}^{\pi/h}
    \left\lvert\widehat{\psi}(\eta_1,\eta_2)\right\rvert^2{\rm d}\eta_1{\rm d}\eta_2
    \\
    =&2^{\frac{\beta}{2}}\frac{\pi^{\beta}}{h^{\beta}}
    \frac{1}{4\pi^2}\int_{-\pi/h}^{\pi/h}\int_{-\pi/h}^{\pi/h}
    \left\lvert\widehat{\psi}(\eta_1,\eta_2)\right\rvert^2{\rm d}\eta_1{\rm d}\eta_2
    = 2^{\frac{\beta}{2}}\frac{\pi^{\beta}}{h^{\beta}} \Vert \psi \Vert_{L_h^2}^2.
\end{aligned}
\end{equation*}

Note also that 
\begin{equation*}
    {\rm min}\{1,2^{1-s}\}\left(x^2+y^2\right)^s \leq (x^{2 s}+y^{2 s}) \leq {\rm max}\{1,2^{1-s}\}(x^2+y^2)^s,\,\, s>0,  
\end{equation*}
from which we have 
\begin{equation*}
    2^{\frac{\beta-2}{2}}(|\eta_1|^{2\cdot \frac{\beta}{2}}+|\eta_2|^{2\cdot \frac{\beta}{2}})
    \leq \left(|\eta_1|^2 +|\eta_2|^2 \right)^{\frac{\beta}{2}}.
\end{equation*}
Then Eq. \eqref{eq:Parseval} and inequality \eqref{eq:eta1eta2} imply that  
\begin{equation*}
\begin{split}
    \left(\frac{2}{\pi}\right)^{\beta}\frac{2^{\frac{\beta-2}{2}}}{4\pi^2}
    \lvert \psi\rvert^2_{H_h^{\frac{\beta}{2},\frac{\beta}{2}}(\mathbb{R}^2)}
    =&  
    \left(\frac{2}{\pi}\right)^{\beta}\frac{2^{\frac{\beta-2}{2}}}{4\pi^2}
    \int_{-\pi/h}^{\pi/h}\int_{-\pi/h}^{\pi/h}
    (|\eta_1|^{2\cdot \frac{\beta}{2}}+|\eta_2|^{2\cdot \frac{\beta}{2}})
    \left\lvert\widehat{\psi}(\eta_1,\eta_2)\right\rvert^2{\rm d}\eta_1{\rm d}\eta_2
    \\
    \leq& 
    \left(\frac{2}{\pi}\right)^{\beta}
    \frac{1}{4\pi^2}
    \int_{-\pi/h}^{\pi/h}\int_{-\pi/h}^{\pi/h}
    \left(|\eta_1|^2 +|\eta_2|^2 \right)^{\frac{\beta}{2}}
    \left\lvert\widehat{\psi}(\eta_1,\eta_2)\right\rvert^2{\rm d}\eta_1{\rm d}\eta_2
    \\
    \leq&  
    \frac{1}{4\pi^2}\frac{1}{h^{\beta}}
    \int_{-\pi/h}^{\pi/h}\int_{-\pi/h}^{\pi/h} 
    \left[4{\rm sin}^2(\frac{\eta_1h}{2})+4{\rm sin}^2(\frac{\eta_2h}{2}) \right]^{\frac{\beta}{2}}
    \left\lvert\widehat{\psi}(\eta_1,\eta_2)\right\rvert^2{\rm d}\eta_1{\rm d}\eta_2
    \\
    =& 
    \left( \left(-\Delta_h\right)^{\frac{\beta}{2}}\psi,\psi \right)_{h^2}.   
 \end{split}
\end{equation*}
Here $|\cdot|_{H_h^{\frac{\beta}{2},\frac{\beta}{2}}(\mathbb{R}^2)}$ denotes the fractional Sobolev semi-norm defined by  
\begin{equation*}
    \lvert \psi\rvert^2_{H_h^{\frac{\beta}{2},\frac{\beta}{2}}(\mathbb{R}^2)}
    = \int_{-\pi/h}^{\pi/h}\int_{-\pi/h}^{\pi/h}
    (|\eta_1|^{2\cdot \frac{\beta}{2}}+|\eta_2|^{2\cdot \frac{\beta}{2}})
    \left\lvert\widehat{\psi}(\eta_1,\eta_2)\right\rvert^2{\rm d}\eta_1{\rm d}\eta_2.
\end{equation*}
By the fractional Poincare inequality \cite{HaoCao2021} 
\begin{equation*}
    C \Vert \psi \Vert_{L_h^2}^2 \leq \lvert \psi\rvert^2_{H_h^{\frac{\beta}{2},\frac{\beta}{2}}},
\end{equation*}
there holds 
\begin{equation*}
    C\left(\frac{2}{\pi}\right)^{\beta}\frac{2^{\frac{\beta-2}{2}}}{4\pi^2} \Vert \psi \Vert_{L_h^2}^2 
    \leq \left( \left(-\Delta_h\right)^{\frac{\beta}{2}}\psi,\psi \right)_{h^2}, 
\end{equation*}
with $C$ being a positive constant independent of $h$.  
\end{proof}
\end{lemma}

\section{Fully discrete schemes with L2-1$_{\sigma}$ formula in temporal discretization}
\label{sec:Scheme}
In this section, we derive fully discrete schemes for Eqs. \eqref{eq:FAE} and \eqref{eq:FAE2} based on L2-1$_{\sigma}$ formula in temporal discretization. The corresponding numerical stability analysis and error estimates are investigated as well. 

In order to apply the L2-1$_{\sigma}$ formula which evaluates the temporal Caputo-Hadamard derivative at the non-integer point $t= t_{k+\sigma}$ via a linear combination of the values at integer points $t_i$, it is necessary to pre-process Eqs. \eqref{eq:FAE} and \eqref{eq:FAE2} using the following assertion.  

\begin{lemma}\label{lemma:tk+sigma}
Assume that $\varphi(t)\in C^2[\tilde{a},T]$, $0< \alpha <1$, and $\gamma\in(0,1)$. Then 
\begin{equation*}
    \varphi(t_{k+\gamma}) = (1-\gamma)\varphi(t_{k}) + \gamma\varphi(t_{k+1}) + \mathcal{O}(\tau^2), 
    \ 
    0<k<N-1. 
\end{equation*}
\begin{proof}
The above equality can be readily derived via Taylor's expansion and so is omitted here. 
\end{proof}    
\end{lemma}

In view of Lemma \ref{lemma:tk+sigma}, Eqs. \eqref{eq:FAE} and \eqref{eq:FAE2} can be approximated by the following systems, 
\begin{equation}\label{eq:FAE_preprocess}
    \left\{
    \begin{array}{ll}
    \left.{}_{CH}{\rm D}^{\alpha}_{\tilde{a},t}u(x,t)\right|_{t_{k+\sigma}}
    -(1-\sigma){}_{RZ}{\rm D}^{\beta}_{x}u(x,t_k)
    -\sigma{}_{RZ}{\rm D}^{\beta}_{x}u(x,t_{k+1})
    \\
    \hfill
    =f(x,t_{k+\sigma}) + \mathcal{O}(\tau^2), 
    \ 
    (x,t)\in\Omega\times(\tilde{a},T], 
    \\
    u(a,t) = u(b,t) = 0, 
    \ 
    \hfill t\in(\tilde{a},T], 
    \\
    u(x,\tilde{a}) = u_{0}(x), 
    \ 
    \hfill x\in\Omega, 
    \end{array}
    \right.
\end{equation} 
and 
\begin{equation}\label{eq:FAE2_preprocess}
    \left\{
    \begin{array}{ll}
    \left.{}_{CH}{\rm D}^{\alpha}_{\tilde{a},t}u(x,y,t)\right|_{t_{k+\sigma}}
    +(1-\sigma)\left(-\Delta\right)^{\frac{\beta}{2}}u(x,y,t_k)
    +\sigma\left(-\Delta\right)^{\frac{\beta}{2}}u(x,y,t_{k+1})
    \\
    \hfill 
    =f(x,y,t_{k+\sigma}) + \mathcal{O}(\tau^2), 
    \ 
    (x,y)\in\widetilde{\Omega},t \in (\tilde{a},T], 
    \\
    u(x,y,t) = 0, 
    \hfill 
    (x,y)\in \mathbb{R}^2\backslash \widetilde{\Omega}, t \in (\tilde{a},T], 
    \\
    u(x,y,\tilde{a}) = u_{0}(x,y), 
    \hfill 
    (x,y)\in\widetilde{\Omega}, 
    \end{array}
    \right.
\end{equation}
respectively. 

Once the approximation in Lemma \ref{lemma:tk+sigma} is applied, the following lemma which is useful for numerical analysis of the coming algorithms can be derived.  
\begin{lemma}\label{lemma:cikv}
For $\alpha \in (0,1)$, $\sigma = 1 - \frac{\alpha}{2}$, and the coefficients $c_{i,k}^{(\alpha,\sigma)}$ defined by Eqs. \eqref{eq:c_{j,k}_1} and \eqref{eq:c_{j,k}}{\rm :} \\
{\rm (1)} For $v^0, v^1, \cdots, v^k, v^{k+1} \in \mathcal{U}_h$, the following inequality holds when the temporal stepsize $\tau$ is sufficiently small, 
\begin{equation}\label{eq:cikv_1D}
    \sum_{i=1}^{k+1}c^{(\alpha,\sigma)}_{i,k}\left(v^i-v^{i-1},\sigma v^{k+1} + (1-\sigma)v^{k}\right)_h \geq 
    \frac{1}{2}\sum_{i=1}^{k+1}c^{(\alpha,\sigma)}_{i,k}\left(\Vert v^i \Vert_{h}^2  
    -\Vert v^{i-1} \Vert_{h}^2  \right), 
\end{equation}
which corresponds to Eq. \eqref{eq:FAE_preprocess}. \\
{\rm (2)} For $w^0, w^1, \cdots, w^k, w^{k+1} \in \mathcal{S}_{h}^\circ$, the following inequality holds if $\tau$ is sufficiently small,
\begin{equation}\label{eq:cikv_2D}
    \sum_{i=1}^{k+1}c^{(\alpha,\sigma)}_{i,k}\left(w^i-w^{i-1},\sigma w^{k+1} + (1-\sigma)w^{k}\right)_{h^2} \geq 
    \frac{1}{2}\sum_{i=1}^{k+1}c^{(\alpha,\sigma)}_{i,k}\left(\Vert w^i \Vert_{L_h^2}^2  
    -\Vert w^{i-1} \Vert_{L_h^2}^2  \right),  
\end{equation}
which corresponds to Eq. \eqref{eq:FAE2_preprocess}. 
\begin{proof}
In the following, we only prove inequality \eqref{eq:cikv_1D}. The proof of inequality \eqref{eq:cikv_2D} can be shown using the same techniques and so is omitted. 

We first show that for arbitrary $k=0,1,2,\cdots$, there holds 
\begin{equation}\label{eq:ckk_vk1}
\begin{split}
    &\sum_{i=1}^{k+1}c^{(\alpha,\sigma)}_{i,k}(v^i-v^{i-1},v^{k+1})_h 
    \\
    \geq &\frac{1}{2}\sum_{i=1}^{k+1}c^{(\alpha,\sigma)}_{i,k}\left(\Vert v^i \Vert_{h}^2 -\Vert v^{i-1} \Vert_{h}^2  \right) + \frac{1}{2c^{(\alpha,\sigma)}_{k+1,k}} 
    \left\Vert \sum_{i=1}^{k+1}c^{(\alpha,\sigma)}_{i,k}(v^i-v^{i-1}) \right\Vert_{h}^2.
\end{split}
\end{equation} 
The case with $k=0$ is evident because of
\begin{equation*}
    (v^{1}-v^0, v^{1})_h = \frac{1}{2}(\Vert v^{1} \Vert_{h}^2-\Vert v^{0} \Vert_{h}^2) + \frac{1}{2} \Vert v^{1} - v^{0} \Vert_{h}^2.
\end{equation*} 
For the the case with $k\geq 1$, we have  
\begin{equation*}
\begin{split}
    A 
    &\triangleq 
    \sum_{i=1}^{k+1}c^{(\alpha,\sigma)}_{i,k}(v^i-v^{i-1},v^{k+1})_h
    -\frac{1}{2}\sum_{i=1}^{k+1}c^{(\alpha,\sigma)}_{i,k}\left(\Vert v^i \Vert_{h}^2 -\Vert v^{i-1} \Vert_{h}^2  \right)
    \\
    & = \sum_{i=1}^{k+1}c^{(\alpha,\sigma)}_{i,k}\left(v^i-v^{i-1},v^{k+1}-\frac{v^i+v^{i-1}}{2}\right)_h
    \\
    & = \sum_{i=1}^{k+1}c^{(\alpha,\sigma)}_{i,k}\left(v^i-v^{i-1},\frac{v^i-v^{i-1}}{2}+ (v^{k+1}-v^{i})\right)_h
    \\
    & = \frac{1}{2}\sum_{i=1}^{k+1}c^{(\alpha,\sigma)}_{i,k}\left\Vert v^i-v^{i-1}\right\Vert_{h}^2 + \sum_{i=1}^{k}c^{(\alpha,\sigma)}_{i,k}\left(v^i-v^{i-1},\sum_{s=i+1}^{k+1}(v^{s}-v^{s-1})\right)_h
    \\
    & = \frac{1}{2}\sum_{i=1}^{k+1}c^{(\alpha,\sigma)}_{i,k}\left\Vert v^i-v^{i-1}\right\Vert_{h}^2 + \sum_{s=2}^{k+1}\left(v^s-v^{s-1},\sum_{i=1}^{s-1}c^{(\alpha,\sigma)}_{i,k}(v^{i}-v^{i-1})\right)_h.
\end{split}
\end{equation*}
Define
\begin{equation*}
    w_{s-1}
    =
    \left\{
    \begin{array}{ll}
    0, 
    & 
    s=1, 
    \\
    \sum\limits_{i=1}^{s-1}c^{(\alpha,\sigma)}_{i,k}(v^{i}-v^{i-1}), 
    & 
    s\geq 2. 
    \end{array}
    \right.
\end{equation*} 
Then 
\begin{equation*} 
\begin{split} 
     v^{s}-v^{s-1} =\frac{1}{c^{(\alpha,\sigma)}_{s,k}}(w_{s}-w_{s-1}),
     \, 
     s\geq 1, 
\end{split} 
\end{equation*}
which yields that 
\begin{equation*}
\begin{split}
    A 
    &= \frac{1}{2}\sum_{i=1}^{k+1}\frac{1}{c^{(\alpha,\sigma)}_{i,k}}\left\Vert w_i-w_{i-1}\right\Vert_{h}^2 + \sum_{s=2}^{k+1}\left(\frac{1}{c^{(\alpha,\sigma)}_{s,k}}(w_{s}-w_{s-1}),w_{s-1}\right)_h
    \\
    &= \frac{1}{2c^{(\alpha,\sigma)}_{1,k}}\left\Vert w_{1}\right\Vert_{h}^2 + \frac{1}{2}\sum_{s=2}^{k+1}\frac{1}{c^{(\alpha,\sigma)}_{s,k}} \left(\left\Vert w_s\right\Vert_{h}^2-\left\Vert w_{s-1}\right\Vert_{h}^2 \right)
    \\
    &= \frac{1}{2c^{(\alpha,\sigma)}_{k+1,k}}\left\Vert w_{k
    +1}\right\Vert_{h}^2 + \frac{1}{2}\sum_{s=1}^{k} 
    \left(\frac{1}{c^{(\alpha,\sigma)}_{s,k}}-\frac{1}{c^{(\alpha,\sigma)}_{s+1,k}} \right) \left\Vert w_s\right\Vert_{h}^2
    \\
    &\geq \frac{1}{2c^{(\alpha,\sigma)}_{k+1,k}}\left\Vert w_{k
    +1}\right\Vert_{h}^2 = \frac{1}{2c^{(\alpha,\sigma)}_{k+1,k}} 
    \left\Vert \sum_{i=1}^{k+1}c^{(\alpha,\sigma)}_{i,k}(v^i-v^{i-1}) \right\Vert_{h}^2.
\end{split}
\end{equation*}
Here Lemma \ref{lemma:lc_kk} is used. Therefore, inequality \eqref{eq:ckk_vk1} holds for an arbitrary non-negative integer $k$. 

Now we prove that for $k=0,1,2,\cdots$, there holds 
\begin{equation}\label{eq:ckk_vk2}
\begin{split}
    &\sum_{i=1}^{k+1}c^{(\alpha,\sigma)}_{i,k}(v^i-v^{i-1},v^{k})_h 
    \\
    \geq &\frac{1}{2}\sum_{i=1}^{k+1}c^{(\alpha,\sigma)}_{i,k}\left(\Vert v^i \Vert_{h}^2 -\Vert v^{i-1} \Vert_{h}^2 \right) - \frac{1}{2(c^{(\alpha,\sigma)}_{k+1,k}-c^{(\alpha,\sigma)}_{k,k})} 
    \left\Vert \sum_{i=1}^{k+1}c^{(\alpha,\sigma)}_{i,k}(v^i-v^{i-1}) \right\Vert_{h}^2.
\end{split}
\end{equation}
In fact, the techniques for proving \eqref{eq:ckk_vk1} yield that 
\begin{equation*}
\begin{split}
    &\sum_{i=1}^{k+1}c^{(\alpha,\sigma)}_{i,k}(v^i-v^{i-1},v^{k})_h 
    - \frac{1}{2}\sum_{i=1}^{k+1}c^{(\alpha,\sigma)}_{i,k}\left(\Vert v^i \Vert_{h}^2 -\Vert v^{i-1} \Vert_{h}^2 \right) 
    \\
    &+ \frac{1}{2(c^{(\alpha,\sigma)}_{k+1,k}-c^{(\alpha,\sigma)}_{k,k})} 
    \left\Vert \sum_{i=1}^{k+1}c^{(\alpha,\sigma)}_{i,k}(v^i-v^{i-1}) \right\Vert_{h}^2
    \\
    =& \sum_{i=1}^{k+1}c^{(\alpha,\sigma)}_{i,k}(v^i-v^{i-1},v^{k+1})_h 
    - \frac{1}{2}\sum_{i=1}^{k+1}c^{(\alpha,\sigma)}_{i,k}\left(\Vert v^i \Vert_{h}^2 -\Vert v^{i-1} \Vert_{h}^2 \right)
    \\
    &-(\sum_{i=1}^{k+1}c^{(\alpha,\sigma)}_{i,k}(v^i-v^{i-1}),v^{k+1}-v^{k})_h + \frac{1}{2(c^{(\alpha,\sigma)}_{k+1,k}-c^{(\alpha,\sigma)}_{k,k})} 
    \left\Vert \sum_{i=1}^{k+1}c^{(\alpha,\sigma)}_{i,k}(v^i-v^{i-1}) \right\Vert_{h}^2
    \\
    =&\frac{1}{2c^{(\alpha,\sigma)}_{k+1,k}}\left\Vert w_{k
    +1}\right\Vert_{h}^2 + \frac{1}{2}\sum_{s=1}^{k} 
    \left(\frac{1}{c^{(\alpha,\sigma)}_{s,k}}-\frac{1}{c^{(\alpha,\sigma)}_{s+1,k}} \right) \left\Vert w_s\right\Vert_{h}^2 
    \\
    &- (w_{k+1}, \frac{1}{c^{(\alpha,\sigma)}_{k+1,k}}(w_{k+1}-w_{k}))_h+\frac{1}{2(c^{(\alpha,\sigma)}_{k+1,k}-c^{(\alpha,\sigma)}_{k,k})} 
    \left\Vert w_{k+1} \right\Vert_{h}^2
    \\
    =& \frac{1}{2}\sum_{s=1}^{k-1} 
    \left(\frac{1}{c^{(\alpha,\sigma)}_{s,k}}-\frac{1}{c^{(\alpha,\sigma)}_{s+1,k}} \right) \left\Vert w_s\right\Vert_{h}^2 + \frac{c^{(\alpha,\sigma)}_{k,k}}{2c^{(\alpha,\sigma)}_{k+1,k}(c^{(\alpha,\sigma)}_{k+1,k}-c^{(\alpha,\sigma)}_{k,k})}\left\Vert w_{k+1}+\frac{c^{(\alpha,\sigma)}_{k+1,k}-c^{(\alpha,\sigma)}_{k,k}}{c^{(\alpha,\sigma)}_{k,k}}w_{k} \right\Vert_{h}^2>0. 
\end{split}    
\end{equation*} 
Combining inequalities \eqref{eq:ckk_vk1} and \eqref{eq:ckk_vk2} with Lemma \ref{lemma:lc_kksigma} shows that \eqref{eq:cikv_1D} is valid. 
\end{proof}
\end{lemma}

\subsection{Fully discrete scheme for Eq. \eqref{eq:FAE}}
In this part, we introduce a finite difference scheme for Eq. \eqref{eq:FAE}, where the L2-1$_{\sigma}$ formula and the weighted and shifted Gr\"{u}nwald-Letnikov formula are applied. Numerical stability and error estimate are discussed as well. 

Denote the exact solution and the numerical solution at the grid point $(x_j,t_n)$ by $u_j^k$ and $U_j^n$, respectively. Denote also $f_j^{n+\sigma}=f(x_j,t_{n+\sigma})$. Substituting Eqs. \eqref{eq:L1_CH} and \eqref{eq:RZ3} into Eq. \eqref{eq:FAE_preprocess} which approximates Eq. \eqref{eq:FAE} with 2-nd order accuracy in time, we obtain the following implicit difference scheme after omitting the high-order terms, 
\begin{equation}\label{eq:CH_RZ}
    \left\{
    \begin{array}{ll} 
    \sum\limits_{i=1}\limits^{n+1} c^{(\alpha,\sigma)}_{i,n}\left(U_j^i-U_j^{i-1}\right)
    + h^{-\beta}\sum\limits_{k=0}\limits^{M}r_{j-k}^{(\beta)}(\sigma U^{n+1}_{k}+(1-\sigma)U^{n}_{k})
    =f^{n+\sigma}_{j},
    \\ 
    \hfill 
    0\leq n\leq N-1,\,1\leq j \leq M-1,
    \\
    U^{0}_{j}=u_0(x_j),
    \hfill
    1\leq j \leq M-1,
    \\
    U^{n}_{0}=U^{n}_{M}=0,
    \hfill 
    0\leq n\leq N. 
    \end{array}
    \right.
\end{equation} 

In the following, we give numerical stability analysis and error estimate for the fully discrete scheme \eqref{eq:CH_RZ}. To this end, we first show a fundamental lemma as follows.  

\begin{lemma}\label{lemma:RZh>0} 
For any grid function $v = (v_0,v_1,\cdots,v_M) \in \mathcal{U}_h^{\circ}$ and $M\geq 2$,  
\begin{equation*}
    \frac{c_{\ast}^{(\beta)}}{(b-a)^{\beta}}h\sum\limits_{j=1}\limits^{M-1}v_j^2 \leq h^{1-\beta}\sum\limits_{j=1}\limits^{M-1}\left(\sum\limits_{k=0}\limits^{M}r_{j-k}^{(\beta)}v_{k}\right)v_j 
    \leq 2r_{0}^{(\beta)}h^{1-\beta}\sum\limits_{j=1}\limits^{M-1}v_j^2, 
\end{equation*}
where $r_{k}^{(\beta)}$ is defined by Eq. \eqref{eq:r_k} with $\beta\in(1 ,2)$ and $c_{\ast}^{(\beta)} = \displaystyle\frac{(1-\beta)(2-\beta)(3-\beta)4^{\beta}e^{-\frac{9}{4}}\Psi_{\beta}}{3}$. In other words, 
\begin{equation*}\label{eq:RZh>0} 
    \frac{c_{\ast}^{(\beta)}}{(b-a)^{\beta}} \Vert v \Vert_{h}^2
    \leq  
    h^{1-\beta}\sum\limits_{j=1}^{M-1}\sum\limits_{k=0}\limits^{M}r_{j-k}^{(\beta)}v_{k}v_{j}
    \leq 2r_{0}^{(\beta)}h^{-\beta}\Vert v \Vert_{h}^2.  
\end{equation*}
\begin{proof}
It follows from Lemma \ref{lemma:lr_k} and Eq. \eqref{eq:r_k} that 
\begin{equation}\label{eq:A123}
\begin{split}
    B\triangleq &h^{-\beta}h\sum\limits_{j=1}\limits^{M-1}\left(\sum\limits_{k=0}\limits^{M}r_{j-k}^{(\beta)}v_{k}\right)v_j
    \\
    =& h^{-\beta}\left[h\sum\limits_{j=1}\limits^{M-1}r_{0}^{(\beta)}v_j^2 + h\sum\limits_{j=1}\limits^{M-1}\left(\sum\limits_{k=0, k\neq j  }\limits^{M}r_{j-k}^{(\beta)}v_{k}v_j\right) \right]
    \\
    \geq &h^{-\beta}\left[h\sum\limits_{j=1}\limits^{M-1}r_{0}^{(\beta)}v_j^2 + \frac{h}{2}\sum\limits_{j=1}\limits^{M-1}\sum\limits_{k=0, k\neq j  }\limits^{M}r_{j-k}^{(\beta)}(v_{k}^2 + v_j^2) \right]
    \\
    =&h^{-\beta}\left(h\sum\limits_{j=1}\limits^{M-1}r_{0}^{(\beta)}v_j^2 
    + \frac{h}{2}\sum\limits_{j=1}\limits^{M-1}\sum\limits_{k=0, k\neq j  }\limits^{M}r_{j-k}^{(\beta)}v_{j}^2 
    + \frac{h}{2}\sum\limits_{j=1}\limits^{M-1}\sum\limits_{k=0, k\neq j  }\limits^{M}r_{j-k}^{(\beta)}v_{k}^2 \right)
    \\
    \triangleq & h^{-\beta}(B_1+B_2+B_3).
\end{split}
\end{equation} 
Note that $r_{k}^{(\beta)}<0$ and $r_{-k}^{(\beta)}<0$ when $k\neq 0$. Therefore, 
\begin{equation}\label{eq:A2}
\begin{split}
    B_2 
    &= \frac{h}{2}\sum\limits_{j=1}\limits^{M-1}\sum\limits_{k=0, k\neq j  }\limits^{M}r_{j-k}^{(\beta)}v_{j}^2 
    = \frac{h}{2}\sum\limits_{j=1}\limits^{M-1}v_{j}^2 \sum\limits_{k=0, k\neq j  }\limits^{M}r_{j-k}^{(\beta)}
    \\
    &= \frac{h}{2}\sum\limits_{j=1}\limits^{M-1}v_{j}^2 \sum\limits_{k=j-M,k\neq 0 }\limits^{j}r_{k}^{(\beta)}
    \geq \frac{h}{2}\sum\limits_{j=1}\limits^{M-1}v_{j}^2 \sum\limits_{\lvert k\rvert=1 }\limits^{M-1}r_{k}^{(\beta)}, 
\end{split}
\end{equation}
and 
\begin{equation}\label{eq:A3}
\begin{split}
    B_3 
    &= \frac{h}{2}\sum\limits_{j=1}\limits^{M-1}\sum\limits_{k=0, k\neq j  }\limits^{M}r_{j-k}^{(\beta)}v_{k}^2
    \\
    &= \frac{h}{2}\left(\sum\limits_{j=1}\limits^{M-1}\sum\limits_{k=j-M  }\limits^{-1}r_{k}^{(\beta)}v_{j-k}^2 +\sum\limits_{j=1}\limits^{M-1} \sum\limits_{k=1}\limits^{j}r_{k}^{(\beta)}v_{j-k}^2  \right)
    \\
    &= \frac{h}{2}\left(\sum\limits_{k=1-M}\limits^{-1}r_{k}^{(\beta)}\sum\limits_{j=1}\limits^{k+M}v_{j-k}^2 
    +\sum\limits_{k=1}\limits^{M-1}r_{k}^{(\beta)}\sum\limits_{j=k}\limits^{M-1}v_{j-k}^2  \right)
    \\
    &\geq \frac{h}{2}\left(\sum\limits_{k=1-M}\limits^{-1}r_{k}^{(\beta)} 
    +\sum\limits_{k=1}\limits^{M-1}r_{k}^{(\beta)}\right)\sum\limits_{j=1}\limits^{M-1}v_{j}^2
    \\
    &= \frac{h}{2}\sum\limits_{j=1}\limits^{M-1}v_{j}^2 \sum\limits_{\lvert k\rvert=1 }\limits^{M-1}r_{k}^{(\beta)}.
\end{split}
\end{equation}
Substituting \eqref{eq:A2} and \eqref{eq:A3} into \eqref{eq:A123} gives 
\begin{equation*}
\begin{split}
    B 
    &\geq h^{-\beta}\left(h\sum\limits_{j=1}\limits^{M-1}r_{0}^{(\beta)}v_j^2 
    +h\sum\limits_{j=1}\limits^{M-1}v_{j}^2 \sum\limits_{\lvert k\rvert=1 }\limits^{M-1}r_{k}^{(\beta)}
    \right)
    \\
    &= h^{-\beta}h\left(r_{0}^{(\beta)}+\sum\limits_{\lvert k\rvert=1 }\limits^{M-1}r_{k}^{(\beta)}
    \right)\sum\limits_{j=1}\limits^{M-1}v_j^2 
    \\
    &=h^{-\beta}h\sum\limits_{k=1-M }\limits^{M-1}r_{k}^{(\beta)} \sum\limits_{j=1}\limits^{M-1}v_j^2
    \\
    &\geq 
    \frac{c_{\ast}^{(\beta)}}{M^{\beta}}h^{-\beta}h\sum\limits_{j=1}\limits^{M-1}v_j^2
    >
    \frac{c_{\ast}^{(\beta)}}{(b-a)^{\beta}}h\sum\limits_{j=1}\limits^{M-1}v_j^2 >0,\,\,M\geq 2, 
\end{split}
\end{equation*}
where Lemma \ref{lemma:lr_k} is applied.

In the meanwhile, 
\begin{equation}\label{eq:A123_2}
\begin{split}
    B 
    \leq &h^{-\beta}\left[h\sum\limits_{j=1}\limits^{M-1}r_{0}^{(\beta)}v_j^2 - \frac{h}{2}\sum\limits_{j=1}\limits^{M-1}\sum\limits_{k=0, k\neq j  }\limits^{M}r_{j-k}^{(\beta)}(v_{k}^2 + v_j^2) \right]
    \\ 
    = & h^{-\beta}(B_1-B_2-B_3)
    \\
    \leq & h^{-\beta}h\left(r_{0}^{(\beta)}-\sum\limits_{\lvert k\rvert=1 }\limits^{M-1}r_{k}^{(\beta)}
    \right)\sum\limits_{j=1}\limits^{M-1}v_j^2 
    \\
    = & h^{-\beta}h\left(2r_{0}^{(\beta)}-\sum\limits_{k=1-M}\limits^{M-1}r_{k}^{(\beta)}
    \right)\sum\limits_{j=1}\limits^{M-1}v_j^2 
     \\
    \leq & 2r_{0}^{(\beta)}h^{-\beta}h
    \sum\limits_{j=1}\limits^{M-1}v_j^2,\,\,M\geq 2.
\end{split}
\end{equation}
The proof is thus completed. 
\end{proof}
\end{lemma} 

Now we are ready for  numerical stability analysis of the fully discrete scheme \eqref{eq:CH_RZ}. 
\begin{theorem}\label{TH:3.1}
Let $\alpha\in(0,1)$, $\beta\in(1,2)$, and $\sigma = 1-\frac{\alpha}{2}$. The fully discrete scheme \eqref{eq:CH_RZ} is unconditionally stable. Moreover, the numerical solution satisfies the following priori estimate for sufficiently small temporal stepsize $\tau$ and spatial stepsize $h$,
\begin{equation}\label{eq:Un}
    \Vert U^{n} \Vert_{h}^2 \leq \Vert U^0 \Vert_{h}^2+\frac{(b-a)^{\beta}\Gamma(1-\alpha)}{c_{\ast}^{(\beta)}\tilde{a}^{\alpha}} \max\limits_{0\leq k \leq n-1} \{((k+\sigma)\tau)^{\alpha}\Vert f^{k+\sigma} \Vert_{h}^2\},
    \, 
    1\leq n \leq N,
\end{equation}
where $U^n = (U^n_0,U^n_1,\cdots,U^n_M)$ and $\displaystyle c_{\ast}^{(\beta)} = \frac{(1-\beta)(2-\beta)(3-\beta)4^{\beta}e^{-\frac{9}{4}}\Psi_{\beta}}{3}$.
\begin{proof} 
For $n = 1$, taking the discrete inner product of the first equation in \eqref{eq:CH_RZ} with $\sigma U^{1}+(1-\sigma)U^{0}$ gives
\begin{equation*}
\begin{split}
    &c^{(\alpha,\sigma)}_{1,0}\left(U^1-U^0, \sigma U^{1}+(1-\sigma)U^{0}\right)_h 
    \\
    =& - h^{-\beta}h\sum\limits_{j=1}^{M-1}\sum\limits_{k=0}\limits^{M}r_{j-k}^{(\beta)}\left(\sigma U^{1}_{k}+(1-\sigma)U^{0}_{k}\right)\left(\sigma U_j^{1}+(1-\sigma)U_j^{0}\right)
    + \left(f^{\sigma}, \sigma U^{1}+(1-\sigma)U^{0}\right)_h.
\end{split}
\end{equation*}
Applying Lemmas \ref{lemma:cikv} and \ref{lemma:RZh>0}, and Cauchy-Schwarz inequality to the above equation, we have
\begin{equation*}
\begin{split}
    \frac{1}{2}c^{(\alpha,\sigma)}_{1,0}\left(\Vert U^1 \Vert_{h}^2  
    -\Vert U^{0} \Vert_{h}^2 \right)
    &\leq
    - \frac{c_{\ast}^{(\beta)}}{(b-a)^{\beta}}\Vert \sigma U^{1}+(1-\sigma)U^{0} \Vert_{h}^2
    \\
    &+\frac{c_{\ast}^{(\beta)}}{(b-a)^{\beta}} \Vert \sigma U^{1}+(1-\sigma)U^{0} \Vert_{h}^2 + \frac{(b-a)^{\beta}}{4c_{\ast}^{(\beta)}} \Vert f^{\sigma} \Vert_{h}^2, 
\end{split}
\end{equation*}
which yields  
\begin{equation}\label{eq:U1}
    \Vert U^1 \Vert_{h}^2 
    \leq\Vert U^0 \Vert_{h}^2 + \frac{2}{c^{(\alpha,\sigma)}_{1,0}}\frac{(b-a)^{\beta}}{4c_{\ast}^{(\beta)}} \Vert f^{\sigma} \Vert_{h}^2
    \leq\Vert U^0 \Vert_{h}^2 + \frac{(b-a)^{\beta}\Gamma(1-\alpha)}{c_{\ast}^{(\beta)}\tilde{a}^{\alpha}}
    \left(\sigma\tau\right)^{\alpha} \Vert f^{\sigma} \Vert_{h}^2. 
\end{equation}
Here Lemma \ref{lemma:1/c_kk} is used in the last inequality.  

Assume that \eqref{eq:Un} holds for $1 \leq i \leq n$. For the case with $i=n+1$, taking the inner product of the first equation in \eqref{eq:CH_RZ} with $\sigma U^{n+1}+(1-\sigma)U^{n}$ leads to
\begin{equation*}
\begin{split}
    &\sum\limits_{i=1}\limits^{n+1} c^{(\alpha,\sigma)}_{i,n}\left(U^i-U^{i-1},\sigma U^{n+1}+(1-\sigma)U^{n}\right)_h 
    \\
    =& -h^{-\beta}h\sum\limits_{j=1}^{M-1}\sum\limits_{k=0}\limits^{M}r_{j-k}^{(\beta)}\left(\sigma U^{n+1}_{k}+(1-\sigma)U^{n}_{k}\right)\left(U_j^{n+1}+(1-\sigma)U_j^{n}\right) 
    \\
    & + \left(f^{n+\sigma}, \sigma U^{n+1}+(1-\sigma)U^{n}\right)_h.
\end{split}
\end{equation*}
It follows from Lemmas \ref{lemma:lc_kk}, \ref{lemma:cikv}, \ref{lemma:RZh>0}, and Cauchy-Schwarz inequality that 
\begin{equation*}
\begin{split}
    &\frac{1}{2}\sum_{i=1}^{n+1}c^{(\alpha,\sigma)}_{i,n}\left(\Vert U^i \Vert_{h}^2  
    -\Vert U^{i-1} \Vert_{h}^2  \right)
    \\
    \leq
    &- \frac{c_{\ast}^{(\beta)}\Vert \sigma U^{n+1}+(1-\sigma)U^{n} \Vert_{h}^2}{(b-a)^{\beta}} 
    +\frac{c_{\ast}^{(\beta)}\Vert \sigma U^{n+1}+(1-\sigma)U^{n}\Vert_{h}^2}{(b-a)^{\beta}}  
    +\frac{(b-a)^{\beta}\Vert f^{n+\sigma} \Vert_{h}^2}{4c_{\ast}^{(\beta)}} 
    \\
    =&\frac{(b-a)^{\beta}\Vert f^{n+\sigma} \Vert_{h}^2}{4c_{\ast}^{(\beta)}}. 
\end{split}
\end{equation*}
Hence,  
\begin{equation*}
\begin{split}
    &c^{(\alpha,\sigma)}_{n+1,n}\Vert U^{n+1} \Vert_{h}^2 
    \leq 
    \sum_{i=1}^{n}(c^{(\alpha,\sigma)}_{i+1,n}-c^{(\alpha,\sigma)}_{i,n})\Vert U^i \Vert_{h}^2  + c^{(\alpha,\sigma)}_{1,n}\left(\Vert U^0 \Vert_{h}^2 + \frac{2}{c^{(\alpha,\sigma)}_{1,n}}\frac{(b-a)^{\beta}}{4c_{\ast}^{(\beta)}} \Vert f^{n+\sigma} \Vert_{h}^2\right)
    \\
    \leq&\sum_{i=1}^{n}(c^{(\alpha,\sigma)}_{i+1,n}-c^{(\alpha,\sigma)}_{i,n})\Vert U^i \Vert_{h}^2  
    +c^{(\alpha,\sigma)}_{1,n}
    \left(\Vert U^0 \Vert_{h}^2 + \frac{(b-a)^{\beta}\Gamma(1-\alpha)}{c_{\ast}^{(\beta)}\tilde{a}^{\alpha}}
    \left((n+\sigma)\tau\right)^{\alpha} \Vert f^{n+\sigma} \Vert_{h}^2\right) 
    \\
    \leq & \sum_{i=1}^{n}(c^{(\alpha,\sigma)}_{i+1,n}-c^{(\alpha,\sigma)}_{i,n})\left(\Vert U^0 \Vert_{h}^2+\frac{(b-a)^{\beta}\Gamma(1-\alpha)}{c_{\ast}^{(\beta)}\tilde{a}^{\alpha}} \max\limits_{0\leq k \leq i-1} \left\{((k+\sigma)\tau)^{\alpha}\Vert f^{k+\sigma} \Vert_{h}^2\right\}\right)
    \\
    +&c^{(\alpha,\sigma)}_{1,n}
    \left(\Vert U^0 \Vert_{h}^2 + \frac{(b-a)^{\beta}\Gamma(1-\alpha)}{c_{\ast}^{(\beta)}\tilde{a}^{\alpha}}\left((n+\sigma)\tau\right)^{\alpha} \Vert f^{n+\sigma} \Vert_{h}^2\right).
    \\
    \leq& 
    \left(\sum_{i=1}^{n}(c^{(\alpha,\sigma)}_{i+1,n}-c^{(\alpha,\sigma)}_{i,n})+c^{(\alpha,\sigma)}_{1,n}\right)\left(\Vert U^0 \Vert_{h}^2+\frac{(b-a)^{\beta}\Gamma(1-\alpha)}{c_{\ast}^{(\beta)}\tilde{a}^{\alpha}} \max\limits_{0\leq k \leq n} \left\{((k+\sigma)\tau)^{\alpha}\Vert f^{k+\sigma} \Vert_{h}^2\right\}\right)
    \\
    =&
    c_{n+1,n}^{(\alpha,\sigma)}\left(\Vert U^0 \Vert_{h}^2+\frac{(b-a)^{\beta}\Gamma(1-\alpha)}{c_{\ast}^{(\beta)}\tilde{a}^{\alpha}} \max\limits_{0\leq k \leq n} \left\{((k+\sigma)\tau)^{\alpha}\Vert f^{k+\sigma} \Vert_{h}^2\right\}\right).
\end{split}
\end{equation*}
Consequently,
\begin{equation*}
    \Vert U^{n} \Vert_{h}^2 \leq \Vert U^0 \Vert_{h}^2+\frac{(b-a)^{\beta}\Gamma(1-\alpha)}{c_{\ast}^{(\beta)}\tilde{a}^{\alpha}} \max\limits_{0\leq k \leq n-1} \left\{((k+\sigma)\tau)^{\alpha}\Vert f^{k+\sigma} \Vert_{h}^2\right\}, 1\leq n \leq N.
\end{equation*}
This completes the proof.
\end{proof}
\end{theorem}

\begin{theorem}
Let $\alpha\in(0,1)$, $\beta\in(1,2)$, and $\sigma = 1-\frac{\alpha}{2}$. Assume that $u(x,t)$, $_{RL}{\rm D}_{a,x}^\beta u(x,t)$, $_{RL}{\rm D}_{x,b}^\beta u(x,t)$ and their Fourier transforms belong to $C^3\left([\tilde{a},T], L^1(\mathbb{R})\right)$. Let $u_j^n$ be the exact solution to Eq. \eqref{eq:FAE} at $(x_j,t_n)$ and $U_j^n$ be the numerical solution solved by scheme \eqref{eq:CH_RZ}. Then the error function given by  
\begin{equation*}
    \varepsilon_j^n = U_j^n-u_j^n, 
    \ 
    0\leq j\leq M, 
    \ 
    0\leq n \leq N, 
\end{equation*}
has the following estimate when the temporal stepsize $\tau$ and the spatial stepsize $h$ are sufficiently small,
\begin{equation}\label{eq:varepsC}
    \Vert \varepsilon^n \Vert _h \leq C\left(\tau^{2}+h^2\right), 
    \ 
    0\leq n \leq N, 
\end{equation}
with $C > 0$ being a constant and $\varepsilon^n = (\varepsilon^n_0,\varepsilon^n_1,\cdots,\varepsilon^n_M)$.
\begin{proof}
The case with $n=0$ is trivial. Now we prove the results for $n=1, 2,3, \cdots$. It follows from Eqs. \eqref{eq:FAE_preprocess} and \eqref{eq:CH_RZ} that
\begin{equation}
    \left\{
    \begin{array}{ll} 
    \sum\limits_{i=1}\limits^{n+1} c^{(\alpha,\sigma)}_{i,n}\left(\varepsilon_j^i-\varepsilon_j^{i-1}\right)
    + h^{-\beta}\sum\limits_{k=0}\limits^{M}r_{j-k}^{(\beta)}(\sigma \varepsilon^{n+1}_{k}+(1-\sigma)\varepsilon^{n+1}_{k})
    =R^{n+\sigma}_{j},
    \\ 
    \hfill 
    0\leq n\leq N-1,\,1\leq j \leq M-1,
    \\
    \varepsilon^{0}_{j}=0,
    \hfill
    1\leq j \leq M-1,
    \\
    \varepsilon^{n}_{0}=\varepsilon^{n}_{M}=0,
    \hfill 
    0\leq n\leq N. 
    \end{array}
    \right.
\end{equation}
By Eqs. \eqref{eq:L1_CH} and \eqref{eq:RZ1}, and Lemma \ref{lemma:tk+sigma}, the truncation error $R_j^{n+\sigma}$ satisfies
\begin{equation*}
    \left\lvert R_j^{n+\sigma}\right\rvert \leq 
    \widetilde C\left(\tau^2+h^2\right),
    \ 
    0\leq n \leq N,
    \, 
    1\leq j \leq M-1.  
\end{equation*}
In view of Theorem \ref{TH:3.1}, there holds 
\begin{equation*}
\begin{split}
     \Vert \varepsilon^{n} \Vert_{h}^2 
     &\leq \Vert \varepsilon^0 \Vert_{h}^2+\frac{(b-a)^{\beta}\Gamma(1-\alpha)}{c_{\ast}^{(\beta)}\tilde{a}^{\alpha}} \max\limits_{0\leq k \leq n-1} \left\{((k+\sigma)\tau)^{\alpha}\Vert R^{k+\sigma} \Vert_{h}^2\right\}
     \\
     &\leq \frac{(b-a)^{\beta}\Gamma(1-\alpha)} {c_{\ast}^{(\beta)}\tilde{a}^{\alpha}}T^{\alpha}\left(\widetilde C\left(\tau^2+h^2\right)\right)^2,\,\,1\leq n \leq N,  
\end{split}
\end{equation*}
which yields that 
\begin{equation*}
    \Vert \varepsilon^n \Vert _{h} \leq {C}\left(\tau^{2}+h^2\right), 
    \ 
    1\leq n \leq N.
\end{equation*}
 All this ends the proof. 
\end{proof}
\end{theorem}

\subsection{Fully discrete scheme for Eq. \eqref{eq:FAE2}}\label{Sec:L2_1Sigma_FAE2}
In the current subsection, we establish a numerical scheme for Eq. \eqref{eq:FAE2} and show the corresponding numerical analysis. Adopting the L2-1$_{\sigma}$  approximation in Eq. \eqref{eq:L1_CH} for the temporal Caputo-Hadamard derivative and the fractional centred difference formula in Eq. \eqref{eq:Deltah} for the fractional Laplacian, we obtain the following fully discrete scheme after removing the high-order terms,  
\begin{equation}\label{eq:twodimen}
    \left\{
    \begin{array}{ll} 
    \sum\limits_{i=1}\limits^{n+1} c^{(\alpha,\sigma)}_{i,n}\left(U_{jk}^i-U_{jk}^{i-1}\right)+\left(-\Delta_h\right)^{\frac{\beta}{2}}\left(\sigma U^{n+1}_{jk}-(1-\sigma)U^{n}_{jk}\right)=f^{n+\sigma}_{jk},
    \\ 
    \hfill 
    0\leq n\leq N-1,\,1\leq j, k \leq M-1,
    \\
    U^{0}_{jk}=u_0(x_j,y_k),
    \hfill 
    1\leq j, k \leq M-1,
    \\
    U^{n}_{jk}=0,
    \hfill 
    (x_j,y_k)\in \partial\Omega_h, 0\leq n\leq N.
    \end{array}
    \right.
\end{equation}
Here $U_{jk}^n$ denotes the approximation of $u(x,y,t)$ at $(x_j,y_k,t_n)$. This finite difference system can be written into the following matrix form, 
\begin{equation}\label{eq:twodimen2} 
    \sum\limits_{i=1}\limits^{n+1} c^{(\alpha,\sigma)}_{i,n}\left(\mathbf{U}^{i}-\mathbf{U}^{i-1}\right) 
    +\frac{1}{h^{\beta}}\mathbf{A}(\sigma\mathbf{U}^{n+1}+(1-\sigma)\mathbf{U}^n)=\mathbf{F}^{n+\sigma}, 
    \ 
    0\leq n\leq N-1, 
\end{equation}
where 
$$\mathbf{F}^{n+\sigma} = (f_{11}^{n+\sigma},\cdots,f_{(M-1)1}^{n+\sigma},f_{12}^{n+\sigma},\cdots,f_{(M-1)2}^{n+\sigma},\cdots,f_{1(M-1)}^{n+\sigma},\cdots,f_{(M-1)(M-1)}^{n+\sigma})^{\rm T}$$
and 
$\mathbf{U}^n = \left(\mathbf{U}_1^n,\mathbf{U}_2^n,\cdots,\mathbf{U}_{M-1}^n\right)$
with $\mathbf{U}_j^n = \left(U_{1j}^n,U_{2j}^n, \cdots, U_{(M-1)j}^n\right)^{\rm T}$. 
The coefficient matrix $\mathbf{A}$ is a real block symmetric matrix with Toeplitz blocks, say, 
\begin{equation}\label{eq:AA}
\mathbf{A} = 
\left (\begin{array}{ccccccc}
A_0 &A_1 &A_2 &\cdots &A_{M-3} &A_{M-2}  \\
A_1 &A_0 &A_1 &\cdots &A_{M-4} &A_{M-3}  \\
A_2 &A_1 &A_0 &\cdots &A_{M-5} &A_{M-4}  \\
\vdots & \vdots &\vdots &\ddots &\vdots &\vdots \\
A_{M-3} &A_{M-4} &A_{M-5} &\cdots &A_0 &A_1 \\
A_{M-2} &A_{M-3} &A_{M-4} &\cdots &A_1 &A_0
\end{array}\right),
\end{equation}
with  
\begin{equation*}
A_j = 
\left (\begin{array}{ccccccc}
a_{0,j}^{(\beta)} &a_{1,j}^{(\beta)} &a_{2,j}^{(\beta)} &\cdots &a_{M-3,j}^{(\beta)} &a_{M-2,j}^{(\beta)}  \\
a_{1,j}^{(\beta)} &a_{0,j}^{(\beta)} &a_{1,j}^{(\beta)} &\cdots &a_{M-4,j}^{(\beta)} &a_{M-3,j}^{(\beta)}  \\
a_{2,j}^{(\beta)} &a_{1,j}^{(\beta)} &a_{0,j}^{(\beta)} &\cdots &a_{M-5,j}^{(\beta)} &a_{M-4,j}^{(\beta)}  \\
\vdots & \vdots &\vdots &\ddots &\vdots &\vdots \\
a_{M-3,j}^{(\beta)} &a_{M-4,j}^{(\beta)} &a_{M-5,j}^{(\beta)} &\cdots &a_{0,j}^{(\beta)} &a_{1,j}^{(\beta)}  \\
a_{M-2,j}^{(\beta)} &a_{M-3,j}^{(\beta)} &a_{M-4,j}^{(\beta)} &\cdots &a_{1,j}^{(\beta)} &a_{0,j}^{(\beta)}  
\end{array}\right).
\end{equation*}
It follows from Lemma \ref{lemma:psi} that the matrix $\textbf{A}$ is a real symmetric positive definite matrix \cite{Wan2022}.

Next, we show the unconditional stability and the error estimate for the fully discrete scheme \eqref{eq:twodimen}.

\begin{theorem}\label{TH:3.3}
Let $\alpha\in(0,1)$, $\beta\in(1,2)$, and $\sigma = 1-\frac{\alpha}{2}$. The finite difference scheme \eqref{eq:twodimen} is unconditionally stable and the numerical solution satisfies the following priori estimate,
\begin{equation*}
    \Vert \mathbf{U}^{n} \Vert_{L_h^2}^2  \leq 
    \Vert \mathbf{U}^0 \Vert_{L_h^2}^2 +\frac{8\pi^{2+\beta}\Gamma(1-\alpha)}{C\tilde{a}^{\alpha}2^{\frac{3\beta}{2}}} 
    \max\limits_{0\leq k \leq n-1} \left\{\left((k+\sigma)\tau\right)^{\alpha}\Vert \mathbf{F}^{k+\sigma} \Vert_{L_h^2}^2\right\},\,\,1\leq n \leq N, 
\end{equation*}
provided that the temporal stepsize $\tau$ and spatial stepsize $h$ are sufficiently small. Here $C$ is a positive constant. 
\begin{proof} 
Taking the two-dimensional discrete inner product of Eq. \eqref{eq:twodimen2} with $\sigma \mathbf{U}^{n+1}+(1-\sigma)\mathbf{U}^{n}$ gives
\begin{equation}\label{eq:mathbfU1}
\begin{split}
    &\left(\sum\limits_{i=1}\limits^{n+1} c^{(\alpha,\sigma)}_{i,n}\left(\mathbf{U}^i-\mathbf{U}^{i-1}\right),\sigma \mathbf{U}^{n+1}+(1-\sigma)\mathbf{U}^{n} \right)_{h^2} 
    \\
    =& - \left(\frac{1}{h^{\beta}}\mathbf{A}(\sigma\mathbf{U}^{n+1}+(1-\sigma)\mathbf{U}^n), \sigma \mathbf{U}^{n+1}+(1-\sigma)\mathbf{U}^{n}\right)_{h^2} 
    \\
    & + \left(\mathbf{F}^{n+\sigma}, \sigma \mathbf{U}^{n+1}+(1-\sigma)\mathbf{U}^{n}\right)_{h^2},\,\,0\leq n \leq N-1.
\end{split}
\end{equation}
Applying Lemmas \ref{lemma:psi} and \ref{lemma:cikv}, and Cauchy-Schwarz inequality to Eq. \eqref{eq:mathbfU1} yields
\begin{equation*}
\begin{split}
    \frac{1}{2}\sum_{i=1}^{n+1}c^{(\alpha,\sigma)}_{i,n}&\left(\Vert \mathbf{U}^i \Vert_{L_h^2}^2  
    -\Vert \mathbf{U}^{i-1} \Vert_{L_h^2}^2  \right)
    \leq
    - C\left(\frac{2}{\pi}\right)^{\beta}\frac{2^{\frac{\beta-2}{2}}}{4\pi^2} \Vert \sigma \mathbf{U}^{n+1}+(1-\sigma)\mathbf{U}^{n} \Vert_{L_h^2}^2
    \\
    &+C\left(\frac{2}{\pi}\right)^{\beta}\frac{2^{\frac{\beta-2}{2}}}{4\pi^2} \Vert \sigma \mathbf{U}^{n+1}+(1-\sigma)\mathbf{U}^{n} \Vert_{L_h^2}^2 + \frac{2\pi^{2+\beta}}{C2^{\frac{3\beta}{2}}} \Vert \mathbf{F}^{n+\sigma} \Vert_{L_h^2}^2, 
\end{split}
\end{equation*} 
which implies that  
\begin{equation*}
\begin{split}
    &c^{(\alpha,\sigma)}_{n+1,n}\Vert \mathbf{U}^{n+1} \Vert_{L_h^2}^2 
    \leq \sum_{i=1}^{n}(c^{(\alpha,\sigma)}_{i+1,n}-c^{(\alpha,\sigma)}_{i,n})\Vert \mathbf{U}^i \Vert_{L_h^2}^2  + c^{(\alpha,\sigma)}_{1,n}\left(\Vert \mathbf{U}^0 \Vert_{h}^2 + \frac{2}{c^{(\alpha,\sigma)}_{1,n}}\frac{2\pi^{2+\beta}}{C2^{\frac{3\beta}{2}}} \Vert \mathbf{F}^{n+\sigma} \Vert_{L_h^2}^2\right)
    \\
    \leq& \sum_{i=1}^{n}(c^{(\alpha,\sigma)}_{i+1,n}-c^{(\alpha,\sigma)}_{i,n})\Vert \mathbf{U}^i \Vert_{L_h^2}^2 
    +c^{(\alpha,\sigma)}_{1,n}\left(\Vert \mathbf{U}^0 \Vert_{L_h^2}^2 +\frac{8\pi^{2+\beta}\Gamma(1-\alpha)}{C\tilde{a}^{\alpha}2^{\frac{3\beta}{2}}} \left((n+\sigma)\tau\right)^{\alpha}\Vert \mathbf{F}^{n+\sigma} \Vert_{L_h^2}^2\right). 
\end{split}
\end{equation*}
Following the same manner of proving Theorem \ref{TH:3.1}, it can be obtained by using the mathematical induction that 
\begin{equation*}
    \Vert \mathbf{U}^{n} \Vert_{L_h^2}  \leq 
    \Vert \mathbf{U}^0 \Vert_{L_h^2} +\frac{8\pi^{2+\beta}\Gamma(1-\alpha)}{C\tilde{a}^{\alpha}2^{\frac{3\beta}{2}}} 
    \max\limits_{0\leq k \leq n-1} \left\{\left((k+\sigma)\tau\right)^{\alpha}\Vert \mathbf{F}^{k+\sigma} \Vert_{L_h^2}^2\right\} 
\end{equation*}
holds for arbitrary $n = 1, \cdots, N$. 
\end{proof}
\end{theorem}

\begin{theorem}\label{TH:3.4}
Assume that $\alpha\in(0,1)$, $\beta\in(1,2)$, $\sigma = 1-\frac{\alpha}{2}$, and $u(x,y,t) \in C^3([\tilde{a},T],\mathcal{B}^{2+\beta}(\widetilde{\Omega}))$ is the solution to Eq. \eqref{eq:FAE2}. Let $u_{jk}^n$ be the exact solution to Eq. \eqref{eq:FAE2} at $(x_j, y_k, t_n)$ and $U_{jk}^n$ be the numerical solution determined by scheme \eqref{eq:twodimen}. Then the error function given by
\begin{equation*}
    e_{jk}^n = U_{jk}^n - u_{jk}^n, 
    \ 
    0 \leq j, k \leq M,
    \ 
    0 \leq n \leq N,
\end{equation*}
has the following estimate when the temporal stepsize $\tau$ and the spatial stepsize h are sufficiently small, 
\begin{equation}\label{eq:varepsC2D}
    \Vert e^n \Vert _{L_h^2} \leq {C} \left(\tau^{2}+h^2\right) , 
    \ 
    0\leq n \leq N, 
\end{equation}
with ${C}$ being a positive constant independent of $h$ and $\tau$. 
\begin{proof}
The case with $n=0$ is obvious. It remains to show the case with $n\geq 1$. It follows from Eqs. \eqref{eq:FAE2_preprocess} and \eqref{eq:twodimen} that
\begin{equation}\label{eq:etwodimen}
    \left\{
    \begin{array}{ll} 
    \sum\limits_{i=1}\limits^{n+1} c^{(\alpha,\sigma)}_{i,n}\left(e_{jk}^i-e_{jk}^{i-1}\right)+\left(-\Delta_h\right)^{\frac{\beta}{2}}\left(\sigma e^{n+1}_{jk}-(1-\sigma)e^{n}_{jk}\right)=R^{n+\sigma}_{jk},
    \\ 
    \hfill 
    0\leq n\leq N-1,\,1\leq j, k \leq M-1,
    \\
    e^{0}_{jk}=0,
    \hfill 
    1\leq j, k \leq M-1,
    \\
    e^{n}_{jk}=0,
    \hfill 
    (x_j,y_k)\in \partial\Omega_h, 0\leq n\leq N.
    \end{array}
    \right.
\end{equation}
By Eq. \eqref{eq:L1_CH}, Lemma \ref{lemma:ErrorRL}, and Lemma \ref{lemma:tk+sigma}, it is evident that for $0\leq n\leq N$ and $1\leq j, k \leq M-1$, there holds
\begin{equation}
    \left\lvert R_{jk}^{n+\sigma}\right\rvert \leq \widetilde{C} (\tau^{2}+h^2). 
\end{equation}
In view of Theorem \ref{TH:3.3}, we have
\begin{equation*}
\begin{split}
     \Vert e^{n} \Vert_{L_h^2}^2
     &\leq \Vert e^0 \Vert_{L_h^2}^2+\frac{8\pi^{2+\beta}\Gamma(1-\alpha)}{C\tilde{a}^{\alpha}2^{\frac{3\beta}{2}}} 
    \max\limits_{0\leq k \leq n-1} \left\{\left((k+\sigma)\tau\right)^{\alpha}\Vert R^{k+1} \Vert_{L_h^2}^2\right\}
     \\
     &\leq \frac{8\pi^{2+\beta}\Gamma(1-\alpha)}{C\tilde{a}^{\alpha}2^{\frac{3\beta}{2}}}
     T^{\alpha}\left(\widetilde C\left(\tau^2+h^2\right)\right)^2,
     \,\, 1\leq n \leq N, 
\end{split}
\end{equation*}
which yields that 
\begin{equation*}
    \Vert e^n \Vert _{L_h^2} \leq {C}\left(\tau^{2}+h^2\right), 
    \ 
    1\leq n \leq N.
\end{equation*}
This finishes the proof.
\end{proof}
\end{theorem}

\section{Fully discrete schemes with L1-2 formula in temporal discretization}
\label{sec:Scheme2}
In the above section, we discretize the Caputo-Hadamard derivative at the non-integer grid in time direction. Although the accuracy is the $(3-\alpha)$-th order for L2-1$_{\sigma}$ formula, the temporal accuracy is 2-nd order when it is used to solve Eqs. \eqref{eq:FAE} and  \eqref{eq:FAE2}. In order to increase the accuracy in time direction, we use the L1-2 formula to discretize the time fractional derivative in Eqs. \eqref{eq:FAE} and  \eqref{eq:FAE2}. In this section, we focus on this issue.


\subsection{Fully discrete scheme for Eq. \eqref{eq:FAE}}
Denote the exact solution and the numerical solution at the grid point $(x_j,t_n)$ by $u_j^k$ and $U_j^n$, respectively. Let also $f_j^n=f(x_j,t_n)$. Substituting Eqs. \eqref{eq:L12_CH} and \eqref{eq:RZ3} into Eq. \eqref{eq:FAE}, we obtain the following finite difference scheme after omitting the high-order terms 
\begin{equation}\label{eq:CH_RZ_2}
    \left\{
    \begin{array}{ll}
    c^{(\alpha)}_{1,1}U^{1}_{j}-c^{(\alpha)}_{1,1}U^{0}_{j}+h^{-\beta}\sum\limits_{k=0}\limits^{M}r_{j-k}^{(\beta)}U^{1}_{k}=f^{1}_{j},
    \hfill   
    1\leq j \leq M-1,
    \\
    c^{(\alpha)}_{n,n}U^{n}_{j}+\sum\limits_{i=1}\limits^{n-1}(c^{(\alpha)}_{i,n}-c^{(\alpha)}_{i+1,n})U^{i}_{j}-c^{(\alpha)}_{1,n}U^{0}_{j}
    +h^{-\beta}\sum\limits_{k=0}\limits^{M}r_{j-k}^{(\beta)}U^{n}_{k}
    =f^{n}_{j},
    \\ 
    \hfill 
    2\leq n\leq N,\,1\leq j \leq M-1,
    \\
    U^{0}_{j}=u_0(x_j),
    \hfill
    1\leq j \leq M-1,
    \\
    U^{n}_{0}=U^{n}_{M}=0,
    \hfill 
    0\leq n\leq N. 
    \end{array}
    \right.
\end{equation} 

Numerical stability analysis and error estimate for the fully discrete scheme \eqref{eq:CH_RZ_2} are given as follows.

\begin{theorem}\label{TH:4.1}
Let $\alpha\in(0,0.3738)$ and $\beta\in(1,2)$. If the temporal stepsize $\tau$ and the spatial stepsize $h$ are sufficiently small, then the fully discrete scheme \eqref{eq:CH_RZ_2} is numerically stable and the numerical solution satisfies the following priori estimate,  
\begin{equation}\label{eq:Un_2}
    \Vert U^{n} \Vert_{h} \leq \Vert U^0 \Vert_{h}+C^{(\alpha)}  \max\limits_{1\leq k \leq n} \left\{(k\tau)^{\alpha}\Vert f^{k} \Vert_{h}\right\}, \,\, 1\leq n\leq N,
\end{equation}
with $C^{(\alpha)} =  \max\left\{\displaystyle\frac{6\Gamma(1-\alpha)}{\tilde{a}^{\alpha}}, \frac{3\Gamma(2-\alpha)}{\alpha\tilde{a}^{\alpha}}\right\}$.  
\begin{proof}
We prove \eqref{eq:Un_2} by the mathematical induction. When $n=1$, multiplying  $hU_j^{1}$ on the both sides of the first equation of \eqref{eq:CH_RZ_2} and summing $j$ from 1 to $M-1$ yield that
\begin{equation*}
    c^{(\alpha)}_{1,1}h\sum\limits_{j=1}\limits^{M-1}({U}^{1}_{j})^2 + 
    h^{-\beta}h\sum\limits_{j=1}\limits^{M-1}\left(\sum\limits_{k=0}\limits^{M}r_{j-k}^{(\beta)}{U}^{1}_{k}\right){U}^{1}_{j} = c^{(\alpha)}_{1,1}h\sum\limits_{j=1}\limits^{M-1}{U}^{0}_{j}{U}^{1}_{j} +h\sum\limits_{j=1}\limits^{M-1}{f}^{1}_{j}{U}^{1}_{j}.
\end{equation*}
According to Lemma \ref{lemma:RZh>0} and Cauchy-Schwarz inequality, we have  
\begin{equation*}
    c^{(\alpha)}_{1,1}\Vert U^1 \Vert_h^2 \leq
    c^{(\alpha)}_{1,1}\Vert U^0 \Vert_h\Vert U^1 \Vert_h + \Vert f^1 \Vert_h\Vert U^1 \Vert_h.
\end{equation*}
It follows from Lemma \ref{lemma:1/c_kk_alpha} that 
\begin{equation*}
    \Vert U^1 \Vert _h \leq
    \Vert U^0 \Vert _h + \frac{1}{c^{(\alpha)}_{1,1}}\Vert f^1 \Vert_h
    \leq
    \Vert U^0 \Vert _h + C^{(\alpha)}\tau^{\alpha}\Vert f^1 \Vert_h. 
\end{equation*} 

Assume that \eqref{eq:Un_2} holds for $i=2, 3,\cdots, n-1$. Multiplying $hU_j^n$ on the both sides of the second equation of \eqref{eq:CH_RZ_2} and summing $j$ from 1 to $M-1$ give 
\begin{equation*}
\begin{split}
    &c^{(\alpha)}_{n,n}h\sum\limits_{j=1}\limits^{M-1}({U}^{n}_{j})^2 
    +h^{-\beta}h\sum\limits_{j=1}\limits^{M-1}\left(\sum\limits_{k=0}\limits^{M}r_{j-k}^{(\beta)}{U}^{n}_{k}\right){U}^{n}_{j} 
    \\
    &=h\sum\limits_{j=1}\limits^{M-1}\sum\limits_{i=1}\limits^{n-1}(c^{(\alpha)}_{i+1,n}-c^{(\alpha)}_{i,n}){U}^{i}_j{U}^{n}_{j}
    +c^{(\alpha)}_{1,n}h\sum\limits_{j=1}\limits^{M-1}{U}^{0}_j{U}^{n}_{j}
    +h\sum\limits_{j=1}\limits^{M-1}{f}^{n}_{j}{U}^{n}_{j} . 
\end{split}
\end{equation*}
Using Lemmas \ref{lemma:lc_kk_2}, \ref{lemma:1/c_kk_alpha}, \ref{lemma:RZh>0}, and Cauchy-Schwarz inequality, one has
\begin{equation*}
\begin{split}
    c^{(\alpha)}_{n,n}\Vert U^n \Vert _{h}^2 &\leq
    \sum\limits_{i=1}\limits^{n-1}(c^{(\alpha)}_{i+1,n}-c^{(\alpha)}_{i,n})\Vert U^i \Vert _{h}\Vert U^n \Vert _{h} + c^{(\alpha)}_{1,n}\left(\Vert U^0 \Vert _{h}+\frac{\Vert f^n \Vert _{h}}{c^{(\alpha)}_{1,n}} \right)\Vert U^n \Vert _{h}
    \\
    &\leq  \sum\limits_{i=1}\limits^{n-1}(c^{(\alpha)}_{i+1,n}-c^{(\alpha)}_{i,n})\left(\Vert U^0 \Vert_{h}+C^{(\alpha)}  \max\limits_{1\leq k \leq i} \left\{(k\tau)^{\alpha}\Vert f^{k} \Vert_{h}\right\}\right)
    \Vert U^n \Vert _{h}
    \\
    &+ c^{(\alpha)}_{1,n} \left(\Vert U^0 \Vert_{h}+C^{(\alpha)}  \max\limits_{1\leq k \leq n} \left\{(k\tau)^{\alpha}\Vert f^{k} \Vert_{h}\right\}\right)\Vert U^n \Vert _{h}
    \\
    &\leq\left[\sum\limits_{i=1}\limits^{n-1}(c^{(\alpha)}_{i+1,n}-c^{(\alpha)}_{i,n})+ c^{(\alpha)}_{1,n} \right] \left(\Vert U^0 \Vert_{h}+C^{(\alpha)}  \max\limits_{1\leq k \leq n} \left\{(k\tau)^{\alpha}\Vert f^{k} \Vert_{h}\right\}\right)\Vert U^n \Vert _{h}
    \\
    &=c^{(\alpha)}_{n,n}\left(\Vert U^0 \Vert_{h}+C^{(\alpha)}  \max\limits_{1\leq k \leq n} \left\{(k\tau)^{\alpha}\Vert f^{k} \Vert_{h}\right\}\right)\Vert U^n \Vert _{h}.
\end{split}
\end{equation*}
Therefore,
\begin{equation*}
    \Vert U^n \Vert _{h}\leq \Vert U^0 \Vert_{h}+C^{(\alpha)}  \max\limits_{1\leq k \leq n} \left\{(k\tau)^{\alpha}\Vert f^{k} \Vert_{h}\right\}.
\end{equation*}
This completes the proof.
\end{proof}
\end{theorem}

Next, we discuss the error estimate. We show the proof after displaying the theorem.
\begin{theorem}
Let $\alpha\in(0,0.3738)$ and $\beta\in(1,2)$. Assume that $u(x,t)$, $_{RL}{\rm D}_{a,x}^\beta u(x,t)$, $_{RL}{\rm D}_{x,b}^\beta u(x,t)$ and their Fourier transforms belong to $C^3\left([\tilde{a},T], L^1(\mathbb{R})\right)$. Let $u_j^n$ be the exact solution to Eq. \eqref{eq:FAE} at $(x_j,t_n)$ and $U_j^n$ be the numerical solution given by scheme \eqref{eq:CH_RZ_2}. Then the grid function given by  
\begin{equation*}
    \varepsilon_j^n = U_j^n-u_j^n, 
    \ 
    0\leq j\leq M, 
    \ 
    0\leq n \leq N, 
\end{equation*}
has the following estimate for sufficiently small temporal stepsize $\tau$ and spatial stepsize $h$,
\begin{equation}\label{eq:varepsC2}
    \Vert \varepsilon^n \Vert _h \leq C\left(\tau^{3-\alpha}+h^2\right), 
    \ 
    0\leq n \leq N, 
\end{equation}
with $C > 0$ being a constant.
\begin{proof} 
The case with $n=0$ is trivial. For $1\leq n\leq N$, it follows from Eqs. \eqref{eq:FAE} and \eqref{eq:CH_RZ_2} that
\begin{equation}
    \left\{
    \begin{array}{ll}
    c^{(\alpha)}_{1,1}\varepsilon^{1}_{j}
    -c^{(\alpha)}_{1,1}\varepsilon^{0}_{j} 
    +h^{-\beta}\sum\limits_{k=0}\limits^{M}r_{j-k}^{(\beta)}\varepsilon^{1}_{k}
    =R_j^1,
    \hfill 
    1 \leq j \leq M-1,
    \\
    c^{(\alpha)}_{n,n}\varepsilon^{n}_{j}
    -\sum\limits_{i=1}\limits^{n-1}(c^{(\alpha)}_{i+1,n}-c^{(\alpha)}_{i,n})\varepsilon^{i}_{j}
    -c^{(\alpha)}_{1,n}\varepsilon^{0}_{j} 
    +h^{-\beta}\sum\limits_{k=0}\limits^{M}r_{j-k}^{(\beta)}\varepsilon^{n}_{k}=
    R_j^n,
    \\
    \hfill 2\leq n\leq N,\ 1 \leq j \leq M-1,
    \\
    \varepsilon^{0}_{j}=0,
    \hfill
    1\leq j \leq M-1,
    \\
    \varepsilon^{n}_{0}=\varepsilon^{n}_{M}=0,
    \hfill 
    0\leq n\leq N. 
    \end{array}
    \right.
\end{equation}
By \eqref{eq:L12_CH} and \eqref{eq:RZ1}, the truncation error $R_j^n$ satisfies
\begin{equation*}
\begin{split}
    &\left\lvert R_j^n\right\rvert \leq 
    \widetilde C\left(\tau^{3-\alpha}+h^2\right),
    \ 
    1\leq n\leq N, 
    \, 
    1\leq j\leq M-1.
\end{split}
\end{equation*}

In view of Theorem \ref{TH:4.1}, there holds 
\begin{equation*}
\begin{split}
     \Vert \varepsilon^{n} \Vert_{h}
     &\leq \Vert \varepsilon^0 \Vert_{h}+C^{(\alpha)}  \max\limits_{1\leq k \leq n} \left\{(k\tau)^{\alpha}\Vert R^{k} \Vert_{h}\right\}.
     \\
     &\leq C^{(\alpha)}T^{\alpha}\widetilde C\left(\tau^{3-\alpha}+h^2\right),
     \, 1\leq n \leq N, 
\end{split}
\end{equation*}
which implies that \eqref{eq:varepsC2} is valid. All this ends the proof. 
\end{proof}
\end{theorem}

Next, we consider model \eqref{eq:FAE2} in two space dimensions.
 
\subsection{Fully discrete scheme for Eq. \eqref{eq:FAE2}}
In the present subsection, we consider Eq. \eqref{eq:FAE2}. Adopting the L1-2 approximation given in \eqref{eq:L12_CH} for the temporal Caputo-Hadamard derivative and the fractional centred difference formula in Eq. \eqref{eq:Deltah} for the fractional Laplacian, we obtain the following fully discrete scheme after ignoring the high-order terms,
\begin{equation}\label{eq:twodimen_2}
    \left\{
    \begin{array}{ll}
    c^{(\alpha)}_{1,1}U^{1}_{jk}-c^{(\alpha)}_{1,1}U^{0}_{jk}+\left(-\Delta_h\right)^{\frac{\beta}{2}}U^{1}_{jk} 
    =f^{1}_{jk},
    \hfill 
    1\leq j, k \leq M-1, 
    \\
    c^{(\alpha)}_{n,n}U^{n}_{jk}+\sum\limits_{i=1}\limits^{n-1}(c^{(\alpha)}_{i,n}-c^{(\alpha)}_{i+1,n})U^{i}_{jk}-c^{(\alpha)}_{1,n}U^{0}_{jk}+\left(-\Delta_h\right)^{\frac{\beta}{2}}U^{n}_{jk}=f^{n}_{jk},
    \\ 
    \hfill 
    2\leq n\leq N,\,1\leq j, k \leq M-1,
    \\
    U^{0}_{jk}=u_0(x_j,y_k),
    \hfill 
    1\leq j, k \leq M-1,
    \\
    U^{n}_{jk}=0,
    \hfill 
    (x_j,y_k)\in \partial\Omega_h, 0\leq n\leq N.
    \end{array}
    \right.
\end{equation}
Here $U_{jk}^n$ is the approximation of $u(x,y,t)$ at $(x_j,y_k,t_n)$. With the notations introduced in Sec. \ref{Sec:L2_1Sigma_FAE2}, this finite difference system can be written into the matrix form as follows, 
\begin{equation}\label{eq:twodimen2_2}
    \left\{
    \begin{array}{ll}
    c^{(\alpha)}_{1,1}\mathbf{U}^{1}-c^{(\alpha)}_{1,1}\mathbf{U}^{0}+\frac{1}{h^{\beta}}\mathbf{A}\mathbf{U}^1
    =\mathbf{F}^{1},
    &
    n=1,
    \\
    c^{(\alpha)}_{n,n}\mathbf{U}^{n}+\sum\limits_{i=1}\limits^{n-1}(c^{(\alpha)}_{i,n}-c^{(\alpha)}_{i+1,n})\mathbf{U}^{i}-c^{(\alpha)}_{1,n}\mathbf{U}^{0}+\frac{1}{h^{\beta}}\mathbf{A}\mathbf{U}^{n}=\mathbf{F}^{n}, 
    &
    2\leq n\leq N. 
    \end{array}
    \right.
\end{equation} 

Next, we show the numerical stability and the error estimate for the fully discrete scheme \eqref{eq:twodimen_2}.

\begin{theorem}\label{TH:4.3}
Let $\alpha\in(0,0.3738)$ and $\beta\in(1,2)$. Assume that the temporal stepsize $\tau$ and the spatial stepsize $h$ are sufficiently small. The finite difference scheme \eqref{eq:twodimen_2} is numerically stable and the numerical solution satisfies the following priori estimate,
\begin{equation}\label{eq:TH43}
    \Vert \mathbf{U}^{n} \Vert_{L_h^2}  \leq 
    \Vert \mathbf{U}^0 \Vert_{L_h^2} 
    +C^{(\alpha)}  \max\limits_{1\leq k \leq n}\left\{(k\tau)^{\alpha} \Vert \mathbf{F}^{k} \Vert_{L_h^2}\right\},
    \,
    1\leq n \leq N, 
\end{equation}
where $C^{(\alpha)} =  \max\left\{\displaystyle\frac{6\Gamma(1-\alpha)}{\tilde{a}^{\alpha}}, \frac{3\Gamma(2-\alpha)}{\alpha\tilde{a}^{\alpha}}\right\}$. 
\begin{proof}
We first consider the case with $n=1$. Taking the two-dimensional discrete inner product of the first equation in \eqref{eq:twodimen2_2} with $\mathbf{U}^{1}$ yields 
\begin{equation*}
\begin{aligned}
    c_{1,1}^{(\alpha)}\Vert \mathbf{U}^{1} \Vert_{L_h^2}^2
    =&c_{1,1}^{(\alpha)}\left(\mathbf{U}^{0}, \mathbf{U}^{1}\right)_{h^2}
    -\frac{1}{h^{\beta}}\left(\mathbf{A}\mathbf{U}^{1}, \mathbf{U}^{1}\right)_{h^2}
    +\left(\mathbf{F}^{1}, \mathbf{U}^{1}\right)_{h^2}
    \\
    \leq
    &c_{1,1}^{(\alpha)}\Vert\mathbf{U}^{0}\Vert_{L_h^2} \Vert\mathbf{U}^{1}\Vert_{L_h^2} 
    +\Vert\mathbf{F}^{1}\Vert_{L_h^2}\Vert\mathbf{U}^{1}\Vert_{L_h^2},
\end{aligned}
\end{equation*}
i.e.,
\begin{equation*} 
    \Vert \mathbf{U}^{1} \Vert_{L_h^2}  
    \leq
    \Vert\mathbf{U}^{0}\Vert_{L_h^2} 
    +\frac{1}{c_{1,1}^{(\alpha)}}\Vert\mathbf{F}^{1}\Vert_{L_h^2} 
    \leq     \Vert\mathbf{U}^{0}\Vert_{L_h^2}+C^{(\alpha)}\tau^{\alpha}\Vert\mathbf{F}^{1}\Vert_{L_h^2}, 
\end{equation*}
where Lemma \ref{lemma:1/c_kk_alpha} is applied. As a result, \eqref{eq:TH43} is valid for $n=1$.  

Assume that $\Vert \mathbf{U}^{i} \Vert_{L_h^2} \leq \Vert \mathbf{U}^0 \Vert_{L_h^2} + C^{(\alpha)}  \max\limits_{1\leq k \leq i}\left\{(k\tau)^{\alpha} \Vert \mathbf{F}^{k} \Vert_{L_h^2}\right\}$ holds for $i=2, 3, \cdots, n-1$. Taking the two-dimensional discrete inner product of the second equation in \eqref{eq:twodimen2_2} with $\mathbf{U}^{n}$ gives
\begin{equation*}
\begin{split}
    &\left(c^{(\alpha)}_{n,n}\mathbf{U}^{n}+\sum\limits_{i=1}\limits^{n-1}(c^{(\alpha)}_{i,n}-c^{(\alpha)}_{i+1,n})\mathbf{U}^{i}-c^{(\alpha)}_{1,n}\mathbf{U}^{0},\mathbf{U}^{n} \right)_{h^2} 
    \\
    =& - \left(\frac{1}{h^{\beta}}\mathbf{A}\mathbf{U}^{n}, \mathbf{U}^{n}\right)_{h^2} 
     + \left(\mathbf{F}^{n}, \mathbf{U}^{n}\right)_{h^2}, 
    \ 
    2\leq n \leq N.
\end{split}
\end{equation*}
Applying Lemmas \ref{lemma:lc_kk_2}, \ref{lemma:1/c_kk_alpha}, \ref{lemma:psi}, and Cauchy-Schwarz inequality to the above equation, we have
\begin{equation*}
\begin{split}
    c^{(\alpha)}_{n,n}\Vert \mathbf{U}^n \Vert _{L_h^2}^2
    &\leq \sum\limits_{i=1}\limits^{n-1}(c^{(\alpha)}_{i+1,n}-c^{(\alpha)}_{i,n})\left(\mathbf{U}^{i}, \mathbf{U}^n\right)_{h^2}  +  c^{(\alpha)}_{1,n}\left(\mathbf{U}^{0}, \mathbf{U}^n\right)_{h^2} + \left(\mathbf{F}^{n}, \mathbf{U}^{n}\right)_{h^2}
    \\
    &\leq  \sum\limits_{i=1}\limits^{n-1}(c^{(\alpha)}_{i+1,n}-c^{(\alpha)}_{i,n})\Vert \mathbf{U}^i \Vert _{L_h^2}\Vert \mathbf{U}^n \Vert _{L_h^2} + c^{(\alpha)}_{1,n}\Vert \mathbf{U}^{0} \Vert _{L_h^2}\Vert \mathbf{U}^n \Vert _{L_h^2}+\Vert \mathbf{F}^{n} \Vert _{L_h^2}\Vert \mathbf{U}^n \Vert _{L_h^2}
    \\
    &\leq \left[\sum\limits_{i=1}\limits^{n-1}(c^{(\alpha)}_{i+1,n}-c^{(\alpha)}_{i,n})\Vert \mathbf{U}^i \Vert _{L_h^2} + 
    c^{(\alpha)}_{1,n}\left(\Vert \mathbf{U}^{0} \Vert _{L_h^2} + C^{(\alpha)} \left(n\tau\right)^{\alpha}\Vert \mathbf{F}^{n} \Vert _{L_h^2}  \right]
    \right)\Vert \mathbf{U}^n \Vert _{L_h^2}
    \\
    &\leq \left[\sum\limits_{i=1}\limits^{n-1}(c^{(\alpha)}_{i+1,n}-c^{(\alpha)}_{i,n})+c^{(\alpha)}_{1,n}\right]
    \left(\Vert \mathbf{U}^0 \Vert_{L_h^2} + C^{(\alpha)}  \max\limits_{1\leq k \leq n}\left\{(k\tau)^{\alpha} \Vert \mathbf{F}^{k} \Vert_{L_h^2}\right\}\right)\Vert \mathbf{U}^n \Vert _{L_h^2}
    \\
    &=c^{(\alpha)}_{n,n}\left(\Vert \mathbf{U}^0 \Vert_{L_h^2} + C^{(\alpha)}  \max\limits_{1\leq k \leq n}\left\{(k\tau)^{\alpha} \Vert \mathbf{F}^{k} \Vert_{L_h^2}\right\}\right)\Vert \mathbf{U}^n \Vert _{L_h^2}, 
\end{split}
\end{equation*}
which impies that 
\begin{equation*}
    \Vert \mathbf{U}^{n} \Vert_{L_h^2}  \leq 
    \Vert \mathbf{U}^0 \Vert_{L_h^2} 
    +C^{(\alpha)}  \max\limits_{1\leq k \leq n}\left\{(k\tau)^{\alpha} \Vert \mathbf{F}^{k} \Vert_{L_h^2}\right\}.  
\end{equation*}
All this ends the proof. 
\end{proof}
\end{theorem}

\begin{theorem}\label{TH:4.4}
Let $\alpha\in(0,0.3738)$ and $\beta\in(1,2)$. Assume that $u(x,y,t) \in C^3([\tilde{a},T],\mathcal{B}^{2+\beta}(\widetilde{\Omega}))$ is the solution of Eq. \eqref{eq:FAE2}. Let $u_{jk}^n$ be the exact solution to Eq. \eqref{eq:FAE2} at $(x_j,y_k,t_n)$ and $U_{jk}^n$ be the numerical solution given by the fully discrete scheme \eqref{eq:twodimen_2}. Then the grid function given by
\begin{equation*}
    e_{jk}^n = U_{jk}^n - u_{jk}^n, 
    \ 
    0 \leq j, k \leq M,
    \ 
    0 \leq n \leq N,
\end{equation*}
has the following estimate when the temporal stepsize $\tau$ and the spatial stepsize $h$ are sufficiently small, 
\begin{equation}\label{eq:varepsC2D_2}
    \Vert e^n \Vert _{L_h^2} \leq {C} \left(\tau^{3-\alpha}+h^2\right), \ 
    0\leq n \leq N,
\end{equation}
with ${C}$ being a positive constant independent of $h$ and $\tau$. 
\begin{proof}
The case with $n=0$ is trivial. For $1\leq n\leq N$, it follows from Eqs. \eqref{eq:FAE2} and \eqref{eq:twodimen_2} that
\begin{equation}\label{eq:etwodimen_2}
    \left\{
    \begin{array}{ll}
    c^{(\alpha)}_{1,1}e^{1}_{jk}-c^{(\alpha)}_{1,1}e^{0}_{jk}+\left(-\Delta_h\right)^{\frac{\beta}{2}}e^{1}_{jk} =R^{1}_{jk},
    \hfill 
    1\leq j \leq M-1,1\leq k \leq M-1,
    \\
    c^{(\alpha)}_{n,n}e^{n}_{jk}+\sum\limits_{i=1}\limits^{n-1}(c^{(\alpha)}_{i,n}-c^{(\alpha)}_{i+1,n})e^{i}_{jk}-c^{(\alpha)}_{1,n}e^{0}_{jk}+\left(-\Delta_h\right)^{\frac{\beta}{2}}e^{n}_{jk}=R^{n}_{jk},
    \\
    \hfill 
    2\leq n\leq N,\,1\leq j \leq M-1,1\leq k \leq M-1,
    \\
    e^{0}_{jk}=0,
    \hfill 
    1\leq j \leq M-1,1\leq k \leq M-1,
    \\
    e^{n}_{jk}=0,
    \hfill 
    (x_j,y_k)\in \partial\Omega_h, 0\leq n\leq N.
    \end{array}
    \right.
\end{equation}
By Eq. \eqref{eq:L12_CH} and Lemma \ref{lemma:ErrorRL}, it is evident that  
\begin{equation}
    \left\lvert R_{jk}^n\right\rvert \leq \widetilde{C} (\tau^{3-\alpha}+h^2), 
    \, 
    1\leq j, k \leq M-1, 
    \, 
    1\leq n\leq N. 
\end{equation}
In view of Theorem \ref{TH:4.3}, we have
\begin{equation*}
\begin{split}
     \Vert e^{n} \Vert_{L_h^2}
     &\leq \Vert e^0 \Vert_{L_h^2}+
    C^{(\alpha)}  \max\limits_{1\leq k \leq n}\left\{(k\tau)^{\alpha} \Vert R^{k} \Vert_{L_h^2}\right\} 
     \\
     &\leq  C^{(\alpha)}
     T^{\alpha}\widetilde C\left(\tau^{3-\alpha}+h^2\right),
     \,\, 1\leq n \leq N, 
\end{split}
\end{equation*}
which yields that 
\begin{equation*}
    \Vert e^n \Vert _{L_h^2} \leq {C}\left(\tau^{3-\alpha}+h^2\right). 
\end{equation*}
This ends the proof.
\end{proof}
\end{theorem}

\section{Numerical experiments}\label{sec:Experiments}
In the section, we demonstrate the aforementioned theoretical results on convergence and numerical stability. 

\begin{example}\label{example1}
Consider the following fractional diffusion equation in one space dimension with $\alpha\in(0,1)$ and $\beta\in(1,2)$, 
\begin{equation*}
    {}_{CH}{\rm D}^{\alpha}_{1,t}u(x,t)
    -{}_{RZ}{\rm D}^{\beta}_{x}u(x,t)=f(x,t),
    \ 
    t \in (1,2],
\end{equation*}
on a finite interval $0 \leq x\leq 1$, with a given force term
\begin{equation*}
\begin{split}
    f(x,t)&=\frac{6}{\Gamma{\left(4-\alpha\right)}}x^4(1-x)^4\left(\log t\right)^{3-\alpha}+\frac{\left(\log t\right)^{3-\alpha}}{2\cos(\pi \beta /2)}
    \\
    &\times \sum_{l=1}^{4}\left(\frac{(-1)^l4!(4+l)!}{l!(4-l)!\Gamma(5+l-\beta)}
    \left(x^{4+l-\beta}+(1-x)^{4+l-\beta}\right)\right).
\end{split}
\end{equation*}
The initial condition and the boundary conditions are given by  $u(x,1)=0$ and $u(0,t)=u(1,t)=0$, respectively. Its exact solution is
\begin{equation*}
    u(x,t)=(\log t)^3x^4(1-x)^4. 
\end{equation*}
\end{example}
The numerical error is defined as
\begin{equation*}
    Error=\sqrt{\left(h\sum_{j=1}^{M-1}\left\lvert U_j^N-u(x_j,t_N)\right\rvert^2\right)}. 
\end{equation*}

Tables \ref{Tab:exmp1} and \ref{Tab:exmp2} show the numerical error and the convergence order given by the fully discrete scheme \eqref{eq:CH_RZ}. We can see that the scheme is stable and the numerical error is $\mathcal{O}\left(\tau^{2}+h^2\right)$. Tables \ref{Tab:exmp5} and \ref{Tab:exmp6} show the numerical error and the convergence orders given by the fully discrete scheme \eqref{eq:CH_RZ_2}. The numerical results indicate that the scheme is stable and the numerical error is $\mathcal{O}\left(\tau^{3-\alpha}+h^2\right)$ for $\alpha\in(0,1)$ even though the theoretical result is for $\alpha \in(0,0.3738)$.

\begin{table}[!htbp] 
\caption{Spatial convergence of scheme \eqref{eq:CH_RZ} with $\tau= 0.001$ for Example \ref{example1}}
\vskip 3pt\label{Tab:exmp1}
\centering 
\begin{tabular}{*{9}{c}}
\toprule
\multirow{2}*{$\beta$}&\multirow{2}*{$h$}&\multicolumn{2}{c}{$\alpha= 0.3$}&\multicolumn{2}{c}{$\alpha=0.6$}
&\multicolumn{2}{c}{$\alpha=0.9$}
\\
&&$Error$&{\rm Order}&$Error$&{\rm Order}&$Error$&{\rm Order}
\\
\midrule\multirow{5}*{1.3}
 & $1/2^3$ & 5.18e-5 & -     & 4.88e-5 & -     & 4.48e-5 & -     \\
  & $1/2^4$ & 1.30e-5 & 1.99 & 1.24e-5 & 1.98 & 1.14e-5 & 1.97 \\
  & $1/2^5$ & 3.32e-6 & 1.97 & 3.15e-6 & 1.97 & 2.93e-6 & 1.97 \\
  & $1/2^6$ & 8.41e-7 & 1.98 & 7.99e-7 & 1.98 & 7.42e-7 & 1.98 \\
  & $1/2^7$ & 2.12e-7 & 1.99 & 2.01e-7 & 1.99 & 1.87e-7 & 1.99 \\
\midrule\multirow{5}*{1.5}
 & $1/2^3$ & 5.15e-5 & -     & 4.95e-5 & -     & 4.68e-5 & -     \\
  & $1/2^4$ & 1.27e-5 & 2.02 & 1.23e-5 & 2.01 & 1.17e-5 & 2.07 \\
  & $1/2^5$ & 3.21e-6 & 1.99 & 3.10e-6 & 1.99 & 2.957e-6 & 1.99 \\
  & $1/2^6$ & 8.11e-7 & 1.99 & 7.83e-7 & 1.99 & 7.457e-7 & 1.99 \\
  & $1/2^7$ & 2.04e-7 & 1.99 & 1.97e-7 & 1.99 & 1.877e-7 & 1.99 \\
\midrule\multirow{5}*{1.7}
& $1/2^3$ & 4.72e-5 & -     & 4.60e-5 & -     & 4.44e-5 & -     \\
  & $1/2^4$ & 1.14e-5 & 2.05 & 1.12e-5 & 2.04 & 1.08e-5 & 2.04 \\
  & $1/2^5$ & 2.86e-6 & 2.05 & 2.79e-6 & 2.00 & 2.70e-6 & 2.00 \\
  & $1/2^6$ & 7.18e-7 & 2.00 & 7.01e-7 & 1.99 & 6.78e-7 & 1.99 \\
  & $1/2^7$ & 1.80e-7 & 2.00 & 1.76e-7 & 2.00 & 1.70e-7 & 2.00 \\ 
\bottomrule
\end{tabular}
\end{table} 

\begin{table}[!htbp]
\caption{Temporal convergence  of  scheme \eqref{eq:CH_RZ} with  $h= 0.001$ for Example \ref{example1}}
\vskip 3pt\label{Tab:exmp2}
\centering 
\begin{tabular}{*{9}{c}}
\toprule
\multirow{2}*{$\beta$}&\multirow{2}*{$\tau$}&\multicolumn{2}{c}{$\alpha= 0.3$}&\multicolumn{2}{c}{$\alpha=0.6$}
&\multicolumn{2}{c}{$\alpha=0.9$}
\\
&&$Error$&{\rm Order}&$Error$&{\rm Order}&$Error$&{\rm Order}
\\
\midrule\multirow{5}*{1.3}
& $1/2^3$ & 1.26e-6 & -     & 2.20e-6 & -     & 2.70e-6 & -     \\
  & $1/2^4$ & 3.10e-7 & 2.03 & 5.44e-7 & 2.02 & 6.74e-7 & 2.00 \\
  & $1/2^5$ & 7.52e-8 & 2.04 & 1.34e-7 & 2.02 & 1.67e-7 & 2.01 \\
  & $1/2^6$ & 1.73e-8 & 2.12 & 3.21e-8 & 2.06 & 4.07e-8 & 2.04 \\
  & $1/2^7$ & 3.86e-9 & 2.16 & 7.09e-9 & 2.18 & 9.29e-9 & 2.13 \\
\midrule\multirow{5}*{1.5}
& $1/2^3$ & 1.31e-6 & -     & 2.27e-6 & -     & 2.80e-6 & -     \\
  & $1/2^4$ & 3.21e-7 & 2.03 & 5.61e-7 & 2.01 & 7.00e-7 & 2.00 \\
  & $1/2^5$ & 7.82e-8 & 2.04 & 1.38e-7 & 2.02 & 1.73e-7 & 2.01 \\
  & $1/2^6$ & 1.81e-8 & 2.11 & 3.32e-8 & 2.06 & 4.23e-8 & 2.03 \\
  & $1/2^7$ & 4.01e-9 & 2.18 & 7.40e-9 & 2.17 & 9.71e-9 & 2.12 \\
\midrule\multirow{5}*{1.7}
& $1/2^3$ & 1.35e-6 & -     & 2.32e-6 & -     & 2.86e-6 & -     \\
  & $1/2^4$ & 3.31e-7 & 2.02 & 5.75e-7 & 2.01 & 7.13e-7 & 2.00 \\
  & $1/2^5$ & 8.10e-8 & 2.03 & 1.42e-7 & 2.02 & 1.77e-7 & 2.01 \\
  & $1/2^6$ & 1.90e-8 & 2.09 & 3.43e-8 & 2.05 & 4.34e-8 & 2.03 \\
  & $1/2^7$ & 4.14e-9 & 2.20 & 7.70e-9 & 2.15 & 1.00e-8 & 2.11 \\ 
\bottomrule
\end{tabular}
\end{table}

\begin{table}[!htbp]
\caption{Spatial convergence of scheme \eqref{eq:CH_RZ_2} with $\tau= 0.001$ for Example \ref{example1}}
\vskip 3pt\label{Tab:exmp5}
\centering 
\begin{tabular}{*{9}{c}}
\toprule
\multirow{2}*{$\alpha$}&\multirow{2}*{$h$}&\multicolumn{2}{c}{$\beta= 1.3$}&\multicolumn{2}{c}{$\beta= 1.5$}
&\multicolumn{2}{c}{$\beta= 1.7$}
\\
&&$Error$&{\rm Order}&$Error$&{\rm Order}&$Error$&{\rm Order}
\\
\midrule\multirow{5}*{0.1}
  & $1/2^3$ & 5.34e-5 & -    & 5.25e-5 & -    & 4.78e-5 & -    \\
  & $1/2^4$ & 1.34e-5 & 2.00 & 1.30e-5 & 2.02 & 1.16e-5 & 2.05 \\
  & $1/2^5$ & 3.40e-6 & 1.98 & 3.27e-6 & 1.99 & 2.89e-6 & 2.00 \\
  & $1/2^6$ & 8.63e-7 & 1.98 & 8.25e-7 & 1.99 & 7.27e-7 & 1.99 \\
  & $1/2^7$ & 2.17e-7 & 1.99 & 2.08e-7 & 1.99 & 1.82e-7 & 2.00 \\
\midrule\multirow{5}*{0.2}
  & $1/2^3$ & 5.26e-5 & -    & 5.21e-5 & -    & 4.75e-5 & -     \\
  & $1/2^4$ & 1.32e-5 & 2.00 & 1.29e-5 & 2.02 & 1.15e-5 & 2.05 \\
  & $1/2^5$ & 3.36e-6 & 1.98 & 3.24e-6 & 1.99 & 2.88e-6 & 2.00 \\
  & $1/2^6$ & 8.52e-7 & 1.98 & 8.19e-7 & 1.99 & 7.23e-7 & 2.00 \\
  & $1/2^7$ & 2.15e-7 & 1.99 & 2.06e-7 & 1.99 & 1.81e-7 & 2.00 \\
\midrule\multirow{5}*{0.3}
  & $1/2^3$ & 5.18e-5 & -    & 5.15e-5 & -    & 4.72e-5 & -     \\
  & $1/2^4$ & 1.30e-5 & 1.99 & 1.27e-5 & 2.02 & 1.14e-5 & 2.05 \\
  & $1/2^5$ & 3.32e-6 & 1.97 & 3.21e-6 & 1.99 & 2.86e-6 & 2.00 \\
  & $1/2^6$ & 8.41e-7 & 1.98 & 8.11e-7 & 1.99 & 7.18e-7 & 1.99 \\
  & $1/2^7$ & 2.12e-7 & 1.99 & 2.04e-7 & 1.99 & 1.80e-7 & 2.00 \\ 
\midrule\multirow{5}*{0.4}
  & $1/2^3$ & 5.09e-5 & -    & 5.09e-5 & -    & 4.69e-5 & -    \\
  & $1/2^4$ & 1.28e-5 & 1.99 & 1.26e-5 & 2.02 & 1.14e-5 & 2.05 \\
  & $1/2^5$ & 3.27e-6 & 1.97 & 3.18e-6 & 1.99 & 2.84e-6 & 2.00 \\
  & $1/2^6$ & 8.28e-7 & 1.98 & 8.03e-7 & 1.99 & 7.13e-7 & 1.99 \\
  & $1/2^7$ & 2.09e-7 & 1.99 & 2.02e-7 & 1.99 & 1.79e-7 & 2.00 \\
\midrule\multirow{5}*{0.6}
  & $1/2^3$ & 4.88e-5 & -    & 4.95e-5 & -    & 4.60e-5 & -     \\
  & $1/2^4$ & 1.24e-5 & 1.98 & 1.23e-5 & 2.01 & 1.12e-5 & 2.04 \\
  & $1/2^5$ & 3.15e-6 & 1.97 & 3.10e-6 & 1.98 & 2.79e-6 & 2.00 \\
  & $1/2^6$ & 7.99e-7 & 1.98 & 7.83e-7 & 1.99 & 7.01e-7 & 1.99 \\
  & $1/2^7$ & 2.01e-7 & 1.99 & 1.97e-7 & 1.99 & 1.76e-7 & 2.00 \\ 
\midrule\multirow{5}*{0.8}
  & $1/2^3$ & 4.62e-5 & -    & 4.78e-5 & -    & 4.50e-5 & -     \\
  & $1/2^4$ & 1.18e-5 & 1.97 & 1.19e-5 & 2.01 & 1.09e-5 & 2.04 \\
  & $1/2^5$ & 3.01e-6 & 1.97 & 3.01e-6 & 1.98 & 2.73e-6 & 2.00 \\
  & $1/2^6$ & 7.63e-7 & 1.98 & 7.59e-7 & 1.99 & 6.86e-7 & 1.99 \\
  & $1/2^7$ & 1.92e-7 & 1.99 & 1.91e-7 & 1.99 & 1.72e-7 & 2.00 \\ 
\bottomrule
\end{tabular}
\end{table}

\begin{table}[!htbp]
\caption{Temporal convergence  of  scheme \eqref{eq:CH_RZ_2} with  $h= 0.00005$ for Example \ref{example1}}
\vskip 3pt\label{Tab:exmp6}
\centering 
\begin{tabular}{*{9}{c}}
\toprule
\multirow{2}*{$\alpha$}&\multirow{2}*{$\tau$}&\multicolumn{2}{c}{$\beta= 1.3$}&\multicolumn{2}{c}{$\beta= 1.5$}
&\multicolumn{2}{c}{$\beta= 1.7$}
\\
&&$Error$&{\rm Order}&$Error$&{\rm Order}&$Error$&{\rm Order}
\\
\midrule\multirow{5}*{0.1}
  & $1/2^3$ & 1.01e-7 & -    & 8.17e-8 & -    & 6.45e-8 & -     \\
  & $1/2^4$ & 1.45e-8 & 2.79 & 1.18e-8 & 2.79 & 9.37e-9 & 2.79 \\
  & $1/2^5$ & 2.00e-9 & 2.86 & 1.63e-9 & 2.86 & 1.30e-9 & 2.85 \\
  & $1/2^6$ & 2.68e-10& 2.90 & 2.20e-10& 2.89 & 1.76e-10& 2.89 \\
  & $1/2^7$ & 3.89e-11& 2.79 & 3.16e-10& 2.80 & 2.50e-11& 2.82 \\
\midrule\multirow{5}*{0.2}
  & $1/2^3$ & 2.58e-7 & -    & 2.12e-7 & -    & 1.69e-7 & -     \\
  & $1/2^4$ & 3.87e-8 & 2.73 & 3.21e-8 & 2.72 & 2.58e-8 & 2.71 \\
  & $1/2^5$ & 5.55e-9 & 2.80 & 4.64e-9 & 2.79 & 3.76e-9 & 2.78 \\
  & $1/2^6$ & 7.93e-10& 2.81 & 6.52e-10& 2.83 & 5.33e-9 & 2.82 \\
  & $1/2^7$ & 1.21e-10& 2.71 & 9.95e-11& 2.71 & 7.94e-10& 2.75 \\
\midrule\multirow{5}*{0.3}
  & $1/2^3$ & 4.91e-7 & -    & 4.09e-7 & -    & 3.30e-7 & -     \\
  & $1/2^4$ & 7.64e-8 & 2.68 & 6.44e-8 & 2.67 & 5.26e-8 & 2.65 \\
  & $1/2^5$ & 1.13e-8 & 2.76 & 9.62e-9 & 2.74 & 7.97e-9 & 2.72 \\
  & $1/2^6$ & 1.78e-9 & 2.66 & 1.48e-9 & 2.70 & 1.19e-9 & 2.75 \\
  & $1/2^7$ & 2.85e-10& 2.65 & 2.36e-10& 2.64 & 1.91e-10& 2.64 \\ 
\midrule\multirow{5}*{0.4}
 & $1/2^3$ & 8.25e-7 & -     & 6.95e-7 & -    & 5.68e-7 & -     \\
  & $1/2^4$ & 1.34e-7 & 2.63 & 1.13e-7 & 2.62 & 9.40e-8 & 2.60 \\
  & $1/2^5$ & 2.15e-8 & 2.64 & 1.82e-8 & 2.64 & 1.49e-8 & 2.66 \\
  & $1/2^6$ & 3.57e-9 & 2.59 & 3.00e-9 & 2.60 & 2.43e-9 & 2.61 \\
  & $1/2^7$ & 6.01e-10& 2.57 & 5.04e-10& 2.57 & 4.10e-10& 2.57 \\
\midrule\multirow{5}*{0.6}
 & $1/2^3$ & 2.07e-6 & -     & 1.71e-6 & -    & 1.40e-6 & -     \\
  & $1/2^4$ & 3.79e-7 & 2.45 & 3.24e-7 & 2.40 & 2.67e-7 & 2.38 \\
  & $1/2^5$ & 6.78e-8 & 2.48 & 5.75e-8 & 2.49 & 4.77e-8 & 2.49 \\
  & $1/2^6$ & 1.26e-8 & 2.43 & 1.07e-8 & 2.43 & 8.78e-9 & 2.44 \\
  & $1/2^7$ & 2.37e-9 & 2.41 & 2.01e-9 & 2.41 & 1.66e-9 & 2.41 \\
\midrule\multirow{5}*{0.8}
  & $1/2^3$ & 4.95e-6 & -    & 4.18e-6 & -    & 3.39e-6 & -     \\
  & $1/2^4$ & 1.00e-6 & 2.31 & 8.55e-7 & 2.29 & 7.08e-7 & 2.26 \\
  & $1/2^5$ & 2.02e-7 & 2.31 & 1.71e-7 & 2.32 & 1.43e-7 & 2.31 \\
  & $1/2^6$ & 4.21e-8 & 2.26 & 3.58e-8 & 2.26 & 2.96e-8 & 2.27 \\
  & $1/2^7$ & 8.96e-9 & 2.23 & 7.61e-9 & 2.23 & 6.31e-9 & 2.23 \\ 
\bottomrule
\end{tabular}
\end{table} 

\begin{example}\label{example2}
Consider the following fractional diffusion equation in two space dimensions with $\alpha\in(0,1)$ and $ \beta\in(1,2)$, 
\begin{equation*}
    \left\{
    \begin{array}{ll}
    {}_{CH}{\rm D}^{\alpha}_{1,t}u(x,y,t)+\left(-\Delta\right)^{\frac{\beta}{2}}u(x,y,t)=f(x,y,t),
     &(x,y)\in\widetilde{\Omega},t \in (1,2],
    \\
    u(x,y,t) = 0, 
    &(x,y)\in \mathbb{R}^2\backslash \widetilde{\Omega},t \in (1,2],
    \\
    u(x,y,1) = 0, 
    &(x,y)\in \widetilde{\Omega},
    \end{array}
    \right.
\end{equation*}
where $\widetilde{\Omega}=(-1,1)^2$ and the exact solution is set as $u(x,y,t) = (\log t)^3(1-x^2)^4(1-y^2)^4 $. The source term $f$ is not explicitly known and we use very fine mesh to compute it. 
\end{example}
In the present simulation, we evaluate the source term by $f\approx f_h={}_{CH}{\rm D}^{\alpha}_{1,t}u(x,y,t)+\left(-\Delta_h\right)^{\frac{\beta}{2}}u(x,y,t)$ with $h = 2^{-8}$. The numerical error in the spatial direction is
\begin{equation*}
    E(h) = 
    \sqrt{h^2 \sum_{j=0}^{M}
    \sum_{k=0}^{M}\left\lvert U_{jk}^N(h,\tau) - U_{2j,2k}^N(h/2,\tau) \right\rvert^2},
\end{equation*}
where $\tau$ is small enough. The numerical error in the temporal direction is
\begin{equation*}
    F(\tau) = 
    \sqrt{h^2 \sum_{j=0}^{M}
    \sum_{k=0}^{M}\left\lvert U_{jk}^N(h,\tau) - U_{jk}^{2N}(h,\tau/2) \right\rvert^2},
\end{equation*}
where $h$ is small enough.

Tables \ref{Tab:exmp3} and \ref{Tab:exmp4} display the numerical error and convergence order for the finite difference scheme (\ref{eq:twodimen}). We can observe that the numerical results are stable and of 2-nd order accuracy in both time and space. Tables \ref{Tab:exmp7} and \ref{Tab:exmp8} show the numerical error and convergence order given by the fully discrete scheme \eqref{eq:twodimen_2}. The numerical results indicate that the scheme is numerically stable and the numerical error is of $\mathcal{O}\left(\tau^{3-\alpha}+h^2\right)$ for $\alpha\in(0,1)$,  even though the theoretical results are for $\alpha \in(0,0.3738)$.

From the above numerical simulations, one can see the following facts: 1) Although L2-1$_{\sigma}$ formula has $(3-\alpha)$-th order accuracy, the established schemes \eqref{eq:CH_RZ} and \eqref{eq:twodimen}, respectively for Eqs. \eqref{eq:FAE} and \eqref{eq:FAE2} (or Eqs. \eqref{eq:FAE_preprocess} and \eqref{eq:FAE2_preprocess}), only have 2-nd order accuracy in time direction. The main reason is that the temporal discretization at the non-integer time grid $t_{k+\sigma}(\sigma = 1-\frac{\alpha}{2})$ uses the value at the next level, i.e., $u^{k+1}$. From the theoretical analysis and numerical simulations, the numerical schemes \eqref{eq:CH_RZ} and \eqref{eq:twodimen} are both numerically stable for all $\alpha \in (0,1)$. 2) Although the derived schemes \eqref{eq:CH_RZ_2} and \eqref{eq:twodimen_2}, respectively for Eqs. \eqref{eq:FAE} and \eqref{eq:FAE2}, have $(3-\alpha)$-th order accuracy in time direction, $\alpha$ must be in $(0,0.3738)$ such that the derived schemes are stable. This limitation is caused by the fact that $c_{k-1,k}^{(\alpha)}$ in L1-2 formula (see \eqref{eq:L12_CH}) can be sign-changing with respect to $\alpha$ for any fixed $k$. However, numerical examples and simulations seem to be good for all $\alpha$ in $(0,1)$. So when to choose L2-1$_{\sigma}$ or L1-2 formula depends upon the needs of the research questions.

\begin{table}[!htbp]
\caption{Spatial convergence of scheme \eqref{eq:twodimen} with $\tau= 1/2^7$ for Example \ref{example2}}
\vskip 3pt\label{Tab:exmp3}
\centering 
\begin{tabular}{*{9}{c}}
\toprule
\multirow{2}*{$\beta$}&\multirow{2}*{$h$}&\multicolumn{2}{c}{$\alpha= 0.3$}&\multicolumn{2}{c}{$\alpha=0.6$}
&\multicolumn{2}{c}{$\alpha=0.9$}
\\
&&$E(h)$&{\rm Order}&$E(h)$&{\rm Order}&$E(h)$&{\rm Order}
\\
\midrule\multirow{5}*{1.3}
  & 1/$2^2$ & -       & -    & -       & -    & -       & -    \\
  & 1/$2^3$ & 4.49e-3 & -    & 4.05e-3 & -    & 3.49e-3 & -    \\
  & 1/$2^4$ & 1.06e-3 & 2.08 & 9.64e-4 & 2.07 & 8.36e-4 & 2.06 \\
  & 1/$2^5$ & 2.63e-4 & 2.02 & 2.38e-4 & 2.02 & 2.07e-4 & 2.01 \\
  & 1/$2^6$ & 6.55e-5 & 2.00 & 5.94e-5 & 2.00 & 5.16e-5 & 2.00 \\
\midrule\multirow{5}*{1.5}
  & 1/$2^2$ & -       & -    & -       & -    & -       & -    \\
  & 1/$2^3$ & 5.52e-3 & -    & 5.10e-3 & -    & 4.55e-3 & -    \\
  & 1/$2^4$ & 1.29e-3 & 2.09 & 1.20e-3 & 2.09 & 1.08e-3 & 2.08 \\
  & 1/$2^5$ & 3.18e-4 & 2.02 & 2.96e-4 & 2.02 & 2.65e-4 & 2.02 \\
  & 1/$2^6$ & 7.93e-5 & 2.01 & 7.36e-5 & 2.01 & 6.61e-5 & 2.01 \\
\midrule\multirow{5}*{1.7}
  & 1/$2^2$ & -       & -    & -       & -    & -       & -    \\
  & 1/$2^3$ & 6.59e-3 & -    & 6.21e-3 & -    & 5.70e-3 & -    \\
  & 1/$2^4$ & 1.52e-3 & 2.11 & 1.44e-3 & 2.11 & 1.33e-3 & 2.10 \\
  & 1/$2^5$ & 3.75e-4 & 2.03 & 3.54e-4 & 2.03 & 3.27e-4 & 2.02 \\
  & 1/$2^6$ & 9.32e-5 & 2.01 & 8.81e-5 & 2.01 & 8.13e-5 & 2.01 \\ 
\bottomrule
\end{tabular}
\end{table} 

\begin{table}[!htbp]
\caption{Temporal convergence  of  scheme \eqref{eq:twodimen} with  $h=1/2^6$ for Example \ref{example2}}
\vskip 3pt\label{Tab:exmp4}
\centering 
\begin{tabular}{*{9}{c}}
\toprule
\multirow{2}*{$\beta$}&\multirow{2}*{$\tau$}&\multicolumn{2}{c}{$\alpha= 0.3$}&\multicolumn{2}{c}{$\alpha=0.6$}
&\multicolumn{2}{c}{$\alpha=0.9$}
\\
&&$F(\tau)$&{\rm Order}&$F(\tau)$&{\rm Order}&$F(\tau)$&{\rm Order}
\\
\midrule\multirow{5}*{1.3}
  & 1/$2^2$ & -       & -    & -       & -    & -       & -     \\
  & 1/$2^3$ & 1.05e-3 & -    & 1.80e-3 & -    & 2.18e-3 & -    \\
  & 1/$2^4$ & 2.57e-4 & 2.04 & 4.46e-4 & 2.02 & 5.42e-4 & 2.01 \\
  & 1/$2^5$ & 6.30e-5 & 2.03 & 1.11e-4 & 2.01 & 1.36e-4 & 2.00 \\
  & 1/$2^6$ & 1.55e-5 & 2.02 & 2.74e-5 & 2.01 & 3.39e-5 & 2.00 \\
\midrule\multirow{5}*{1.5}
  & 1/$2^2$ & -       & -    & -       & -    & -       & -     \\
  & 1/$2^3$ & 1.09e-3 & -    & 1.86e-3 & -    & 2.25e-3 & -     \\
  & 1/$2^4$ & 2.65e-4 & 2.04 & 4.59e-4 & 2.02 & 5.60e-4 & 2.01 \\
  & 1/$2^5$ & 6.50e-5 & 2.03 & 1.14e-4 & 2.01 & 1.40e-4 & 2.00 \\
  & 1/$2^6$ & 1.61e-5 & 2.02 & 2.83e-5 & 2.01 & 3.50e-5 & 2.00 \\
\midrule\multirow{5}*{1.7}
  & 1/$2^2$ & -       & -    & -       & -    & -       & -     \\
  & 1/$2^3$ & 1.12e-3 & -    & 1.90e-3 & -    & 2.32e-3 & -     \\
  & 1/$2^4$ & 2.72e-4 & 2.04 & 4.70e-4 & 2.02 & 5.76e-4 & 2.01 \\
  & 1/$2^5$ & 6.69e-5 & 2.02 & 1.17e-4 & 2.01 & 1.44e-4 & 2.00 \\
  & 1/$2^6$ & 1.66e-5 & 2.02 & 2.90e-5 & 2.01 & 3.60e-5 & 2.00 \\
\bottomrule
\end{tabular}
\end{table}

\begin{table}[!htbp]
\caption{Spatial convergence of scheme \eqref{eq:twodimen_2} with $\tau= 1/2^7$ for Example \ref{example2}}
\vskip 3pt\label{Tab:exmp7}
\centering 
\begin{tabular}{*{9}{c}}
\toprule
\multirow{2}*{$\alpha$}&\multirow{2}*{$h$}&\multicolumn{2}{c}{$\beta= 1.3$}&\multicolumn{2}{c}{$\beta= 1.5$}
&\multicolumn{2}{c}{$\beta= 1.7$}
\\
&&$E(h)$&{\rm Order}&$E(h)$&{\rm Order}&$E(h)$&{\rm Order}
\\
\midrule\multirow{5}*{0.1}
& 1/$2^2$ & -          & -     & -          & -     & -          & -     \\
  & 1/$2^3$ & 4.72e-3 & -    & 5.73e-3 & -    & 6.78e-3 & -     \\
  & 1/$2^4$ & 1.12e-3 & 2.08 & 1.34e-3 & 2.10 & 1.57e-3 & 2.11 \\
  & 1/$2^5$ & 2.76e-4 & 2.02 & 3.30e-4 & 2.02 & 3.85e-4 & 2.03 \\
  & 1/$2^6$ & 6.87e-5 & 2.01 & 8.23e-5 & 2.01 & 9.58e-5 & 2.01 \\
\midrule\multirow{5}*{0.2}
  & 1/$2^2$ & -       & -    & -       & -    & -         & -     \\
  & 1/$2^3$ & 4.61e-3 & -    & 5.63e-3 & -    & 6.69e-3 & -     \\
  & 1/$2^4$ & 1.09e-3 & 2.08 & 1.32e-3 & 2.10 & 1.55e-3 & 2.11 \\
  & 1/$2^5$ & 2.70e-4 & 2.02 & 3.25e-4 & 2.02 & 3.80e-4 & 2.03 \\
  & 1/$2^6$ & 6.72e-5 & 2.01 & 8.09e-5 & 2.01 & 9.46e-5 & 2.01 \\
\midrule\multirow{5}*{0.3}
  & 1/$2^2$ & -       & -    & -       & -    & -       & -     \\
  & 1/$2^3$ & 4.49e-3 & -    & 5.52e-3 & -    & 6.59e-3 & -     \\
  & 1/$2^4$ & 1.06e-3 & 2.08 & 1.29e-3 & 2.09 & 1.52e-3 & 2.11 \\
  & 1/$2^5$ & 2.63e-4 & 2.02 & 3.18e-4 & 2.02 & 3.75e-4 & 2.03 \\
  & 1/$2^6$ & 6.55e-5 & 2.00 & 7.93e-5 & 2.00 & 9.32e-5 & 2.01 \\ 
  \midrule\multirow{5}*{0.4}
  & 1/$2^2$ & -       & -     & -          & -     & -          & -     \\
  & 1/$2^3$ & 4.36e-3 & -    & 5.39e-3 & -    & 6.48e-3 & -     \\
  & 1/$2^4$ & 1.03e-3 & 2.08 & 1.26e-3 & 2.09 & 1.50e-3 & 2.11 \\
  & 1/$2^5$ & 2.55e-4 & 2.02 & 3.12e-4 & 2.02 & 3.68e-4 & 2.03 \\
  & 1/$2^6$ & 6.36e-5 & 2.00 & 7.76e-5 & 2.01 & 9.17e-5 & 2.01 \\
\midrule\multirow{5}*{0.6}
  & 1/$2^2$ & -       & -    & -       & -    & -       & -     \\
  & 1/$2^3$ & 4.05e-3 & -    & 5.10e-3 & -    & 6.21e-3 & -     \\
  & 1/$2^4$ & 9.64e-4 & 2.07 & 1.20e-3 & 2.09 & 1.44e-3 & 2.11 \\
  & 1/$2^5$ & 2.38e-4 & 2.02 & 2.96e-4 & 2.02 & 3.54e-4 & 2.03 \\
  & 1/$2^6$ & 5.94e-5 & 2.00 & 7.36e-5 & 2.01 & 8.81e-5 & 2.01 \\
\midrule\multirow{5}*{0.8}
  & 1/$2^2$ & -       & -    & -       & -    & -       & -     \\
  & 1/$2^3$ & 3.69e-3 & -    & 4.75e-3 & -    & 5.89e-3 & -     \\
  & 1/$2^4$ & 8.82e-4 & 2.07 & 1.12e-3 & 2.08 & 1.37e-3 & 2.10 \\
  & 1/$2^5$ & 2.18e-4 & 2.02 & 2.76e-4 & 2.02 & 3.37e-4 & 2.02 \\
  & 1/$2^6$ & 5.44e-5 & 2.00 & 6.88e-5 & 2.01 & 8.38e-5 & 2.01 \\ 
\bottomrule
\end{tabular}
\end{table}

\begin{table}[!htbp]
\caption{Temporal convergence  of  scheme \eqref{eq:twodimen_2} with  $h=1/2^6$ for Example \ref{example2}}
\vskip 3pt\label{Tab:exmp8}
\centering 
\begin{tabular}{*{9}{c}}
\toprule
\multirow{2}*{$\alpha$}&\multirow{2}*{$\tau$}&\multicolumn{2}{c}{$\beta= 1.3$}&\multicolumn{2}{c}{$\beta= 1.5$}
&\multicolumn{2}{c}{$\beta= 1.7$}
\\
&&$F(\tau)$&{\rm Order}&$F(\tau)$&{\rm Order}&$F(\tau)$&{\rm Order}
\\
\midrule\multirow{5}*{0.1}
  & 1/$2^2$ & -       & -    & -       & -    & -       & -     \\
  & 1/$2^3$ & 8.05e-5 & -    & 7.00e-5 & -    & 6.00e-5 & -     \\
  & 1/$2^4$ & 1.06e-5 & 2.93 & 9.21e-6 & 2.93 & 7.88e-6 & 2.93 \\
  & 1/$2^5$ & 1.43e-6 & 2.90 & 1.24e-6 & 2.90 & 1.06e-6 & 2.90 \\
  & 1/$2^6$ & 1.94e-7 & 2.88 & 1.69e-7 & 2.88 & 1.44e-7 & 2.88 \\
\midrule\multirow{5}*{0.2}
  & 1/$2^2$ & -       & -    & -       & -    & -       & -     \\
  & 1/$2^3$ & 2.02e-4 & -    & 1.75e-4 & -    & 1.49e-4 & -     \\
  & 1/$2^4$ & 2.77e-5 & 2.87 & 2.39e-5 & 2.87 & 2.03e-5 & 2.87 \\
  & 1/$2^5$ & 3.90e-6 & 2.83 & 3.35e-6 & 2.83 & 2.85e-6 & 2.83 \\
  & 1/$2^6$ & 5.55e-7 & 2.81 & 4.78e-7 & 2.81 & 4.05e-6 & 2.81 \\
\midrule\multirow{5}*{0.3}
  & 1/$2^2$ & -       & -    & -       & -    & -       & -     \\
  & 1/$2^3$ & 3.89e-4 & -    & 3.34e-4 & -    & 2.82e-4 & -     \\
  & 1/$2^4$ & 5.57e-5 & 2.80 & 4.77e-5 & 2.81 & 4.03e-5 & 2.81 \\
  & 1/$2^5$ & 8.24e-6 & 2.76 & 7.05e-6 & 2.76 & 5.94e-6 & 2.76 \\
  & 1/$2^6$ & 1.24e-6 & 2.73 & 1.06e-7 & 2.73 & 8.93e-7 & 2.73 \\
\midrule\multirow{5}*{0.4}
& 1/$2^2$ & -          & -     & -     & -     & -   & -   \\
  & 1/$2^3$ & 6.77e-4 & -    & 5.78e-4 & -    & 4.85e-3 & -     \\
  & 1/$2^4$ & 1.02e-4 & 2.73 & 8.69e-5 & 2.73 & 7.27e-5 & 2.74 \\
  & 1/$2^5$ & 1.60e-5 & 2.68 & 1.36e-5 & 2.68 & 1.14e-5 & 2.68 \\
  & 1/$2^6$ & 2.56e-6 & 2.64 & 2.17e-6 & 2.64 & 1.82e-6 & 2.65 \\
\midrule\multirow{5}*{0.6}
 & 1/$2^2$ & -        & -    & -       & -    & -       & -     \\
  & 1/$2^3$ & 1.82e-3 & -    & 1.54e-3 & -    & 1.27e-3 & -     \\
  & 1/$2^4$ & 3.07e-4 & 2.57 & 2.58e-4 & 2.58 & 2.13e-4 & 2.58 \\
  & 1/$2^5$ & 5.43e-5 & 2.50 & 4.57e-5 & 2.50 & 3.77e-5 & 2.50 \\
  & 1/$2^6$ & 9.89e-6 & 2.46 & 8.32e-6 & 2.46 & 6.86e-6 & 2.46 \\
\midrule\multirow{5}*{0.8}
  & 1/$2^2$ & -       & -    & -       & -    & -       & -     \\
  & 1/$2^3$ & 4.56e-3 & -    & 3.85e-3 & -    & 3.17e-3 & -     \\
  & 1/$2^4$ & 8.73e-4 & 2.39 & 7.29e-4 & 2.40 & 5.94e-4 & 2.42 \\
  & 1/$2^5$ & 1.74e-4 & 2.32 & 1.46e-4 & 2.32 & 1.19e-4 & 2.32 \\
  & 1/$2^6$ & 3.60e-5 & 2.28 & 3.03e-5 & 2.27 & 2.48e-5 & 2.27 \\
\bottomrule
\end{tabular}
\end{table}

\section*{Data availability}
No data was used for the research described in the article.




\newpage
\begin{thebibliography}{100}

{
\bibitem{Alikhanov2022}
A.A. Alikhanov, C.M. Huang, 
A class of time-fractional diffusion equations with generalized fractional derivatives, 
J. Comput. Appl. Math. 
414 (2022) 114424.}

{
\bibitem{Bai22021}
Z.Z. Bai, K.Y. Lu, 
Optimal rotated block-diagonal preconditioning for discretized optimal control problems constrained with fractional time-dependent diffusive equations, 
Appl. Numer. Math. 
163 (2021) 126-146.}

\bibitem{Bre2021} 
E. Le Breton, S. Brune, K. Ustaszewski, S. Zahirovic, M. Seton, R.D. Muller, 
Kinematics and extent of the Piemont-Liguria Basin-implications for subduction processes in the Alps, 
Solid Earth 
12 (2021) 885-913.

\bibitem{Den2010}
S.I. Denisov, H. Kantz, 
Continuous-time random walk theory of superslow diffusion, 
Europhys. Lett. 
92 (2010) 30001.

\bibitem{DingL2017}
H.F. Ding, C.P. Li, 
High-order numerical algorithms for Riesz derivatives via constructing new generating functions, 
J. Sci. Comput. 
71 (2017) 759-784.

{
\bibitem{Diethelm2010}
K. Diethelm, 
{\em The Analysis of Fractional Differential Equations},
Springer-Verlag, Berlin, 2010.}
 
\bibitem{Fan2022}
E.Y. Fan, C.P. Li, Z.Q. Li, Numerical approaches to Caputo-Hadamard fractional derivatives with applications to long-term integration of fractional differential systems, Commun. Nonlinear Sci. Numer. Simul. 106 (2022) 106096.

\bibitem{HaoCao2021}
Z.P. Hao, W.R. Cao, S.Y. Li, Numerical correction of finite difference solution for two-dimensional space-fractional diffusion equations with boundary singularity, 
Numer. Algorithms, 86 (2021) 1071-1087.

\bibitem{HaoZhang2021}
Z.P. Hao, Z.Q. Zhang, R. Du, Fractional centered difference scheme for high-dimensional integral fractional Laplacian, J. Comput. Phys. 424 (2021) 109851.

{
\bibitem{Jiao2021}
C.Y. Jiao, A. Khaliq, C.P. Li, H.X. Wang, 
Difference between Riesz derivative and fractional Laplacian on the proper subset of $\mathbb{R}$, 
Fract. Calc. Appl. Anal. 
24 (2021) 1716-1734.}

\bibitem{LiCai2019}
C.P. Li, M. Cai, 
{\em Theory and Numerical Approximations of Fractional
Integrals and Derivatives}, 
SIAM, Philadelphia, USA, 2019.

{
\bibitem{LiL2021}
C.P. Li, Z.Q. Li, 
Stability and logarithmic decay of the solution to Hadamard-type fractional differential equation, 
J. Nonlinear Sci. 
31 (2021) 31.}

\bibitem{LiL2020}
C.P. Li, Z.Q. Li, Z. Wang, 
Mathematical analysis and the local discontinuous Galerkin method for Caputo-Hadamard fractional partial differential equation, 
J. Sci. Comput. 
85 (2020) 41.

\bibitem{LiandZeng2015}
C.P. Li, F.H. Zeng,
{\em Numerical Methods for Fractional Calculus}, 
Chapman and Hall/CRC, Boca Raton, USA, 2015. 

{
\bibitem{Lin2017}
X.L. Lin, M.K. Ng, H.W. Sun, 
A splitting preconditioner for Toeplitz-like linear systems arising from fractional diffusion equations, 
SIAM J. Matrix Anal. Appl. 
38 (2017) 1580-1614.}

\bibitem{Lom2016}
A. Lomin, V. M\'{e}ndez, 
Does ultra-slow diffusion survive in a three dimensional cylindrical comb? 
Chaos Solitons Fractals 
82 (2016) 142-147.

\bibitem{Lom1956} 
C. Lomnitz, 
Creep measurements in igneous rocks, 
J. Geol. 
64 (1956) 473-479.

\bibitem{Lom1957}
C. Lomnitz, 
Linear dissipation in solids, 
J. Appl. Phys. 
28 (1957) 201-205.

\bibitem{Lom1962} 
C. Lomnitz, 
Application of the logarithmic creep law to stress wave attenuation in the solid earth, 
J. Geophys. Res. 
67 (1962) 365-368.

\bibitem{Meerschaert2006}
M.M. Meerschaert, C. Tadjeran, Finite difference approximations for two-sided space-fractional partial differential equations, Appl. Numer. Math. 56 (2006) 80-90.

\bibitem{Ou2022}
C.X. Ou, D.K. Cen, S. Vong, Z.B. Wang, 
Mathematical analysis and numerical methods for Caputo-Hadamard fractional diffusion-wave equations, 
Appl. Numer. Math. 
177 (2022) 34-57.

\bibitem{Podlubny1999}
I. Podlubny, 
{\em Fractional Differential Equations}, 
Academic Press, San Diego, CA, 1999.

\bibitem{Sil2013} 
F.L. da Silva, 
EEG and MEG: Relevance to neuroscience, 
Neuron 
80 (2013) 1112-1128.

\bibitem{Tao2020} 
C.H. Tao, W.E. Seyfried, R.P. Lowell, et al., 
Deep high-temperature hydrothermal circulation in a detachment faulting system on the ultra-slow spreading ridge, 
Nat. Commun. 
11 (2020) 1300.

\bibitem{Tian2015}
W.Y. Tian, H. Zhou, W.H. Deng, A class of second order difference approximations for solving space fractional diffusion equations, Math. Comp. 
84 (2015) 1703-1727.

\bibitem{WangCai2023}
Y. Wang, M. Cai, 
Finite difference schemes for time-space fractional diffusion equations in one- and two-dimensions, 
Commun. Appl. Math. Comput. (2023)  
{DOI : 10.1007/s42967-022-00244-8}. 

\bibitem{Wan2022}
Y.Y. Wang, Z.P. Hao, R. Du, 
A linear finite difference scheme for the two-dimensional nonlinear Schr\"{o}dinger equation with fractional Laplacian, 
J. Sci. Comput. 
90 (2022) 24.

\bibitem{War2016} 
J.M. Warren, 
Global variations in abyssal peridotite compositions,
Lithos 
248 (2016) 193-219.


\bibitem{Yang2022}
Z.W. Yang, 
Numerical approximation and error analysis for Caputo-Hadamard fractional stochastic differential equations, 
Z. Angew. Math. Phys. 
73 (2022) 253.

\bibitem{Zaky2022}
M.A. Zaky, A.S. Hendy, D. Suragan, 
Logarithmic Jacobi collocation method for Caputo-Hadamard fractional differential equations, 
Appl. Numer. Math. 
181 (2022) 326-346.

{
\bibitem{Zheng2021}
X.C. Zheng, 
Logarithmic transformation between (variable-order) Caputo and Caputo-Hadamard fractional problems and applications, 
Appl. Math. Lett. 
121 (2021) 107366.}

\end{thebibliography}
\end{document}